\pdfoutput=1
\RequirePackage{ifpdf}
\ifpdf 
\documentclass[pdftex]{sigma}
\else
\documentclass{sigma}
\fi

\usepackage{multirow}
\usepackage{subfigure}

\numberwithin{equation}{section}
\newtheorem{Theorem}{Theorem}[section]
\newtheorem{Corollary}[Theorem]{Corollary}
\newtheorem{Lemma}[Theorem]{Lemma}
\newtheorem{Proposition}[Theorem]{Proposition}
 { \theoremstyle{definition}

\newtheorem{Remark}[Theorem]{Remark} }

\begin{document}

\newcommand{\arXivNumber}{1412.4655}

\allowdisplaybreaks

\renewcommand{\PaperNumber}{059}

\FirstPageHeading

\ShortArticleName{A Perturbation of the Dunkl Harmonic Oscillator on the Line}

\ArticleName{A Perturbation of the Dunkl Harmonic Oscillator\\ on the Line}

\Author{Jes\'us A.~\'ALVAREZ L\'OPEZ~$^\dag$, Manuel CALAZA~$^\ddag$ and Carlos FRANCO~$^\dag$}

\AuthorNameForHeading{J.A.~\'Alvarez L\'opez, M.~Calaza and C.~Franco}

\Address{$^\dag$~Departamento de Xeometr\'{\i}a e Topolox\'{\i}a,
 Facultade de Matem\'aticas,\\
\hphantom{$^\dag$}~Universidade de Santiago de Compostela,
 15782 Santiago de Compostela, Spain}
\EmailD{\href{jesus.alvarez@usc.es}{jesus.alvarez@usc.es}, \href{carlosluis.franco@usc.es}{carlosluis.franco@usc.es}}

\Address{$^\ddag$~Laboratorio de Investigaci\'on 2 and Rheumatology Unit,\\
\hphantom{$^\ddag$}~Hospital Clinico Universitario
de Santiago, Santiago de Compostela, Spain}
\EmailD{\href{manuel.calaza@usc.es}{manuel.calaza@usc.es}}

\ArticleDates{Received February 19, 2015, in f\/inal form July 20, 2015; Published online July 25, 2015;\newline
Corrected June 28, 2017}

\Abstract{Let $J_\sigma$ be the Dunkl harmonic oscillator on ${\mathbb{R}}$ ($\sigma>-\frac{1}{2}$). For $0<u<1$ and $\xi>0$, it is proved that, if $\sigma>u-\frac{1}{2}$, then the operator $U=J_\sigma+\xi|x|^{-2u}$, with appropriate domain, is essentially self-adjoint in $L^2({\mathbb{R}},|x|^{2\sigma} dx)$, the Schwartz space ${\mathcal{S}}$ is a~core of~$\overline U^{1/2}$, and $\overline U$ has a discrete spectrum, which is estimated in terms of the spectrum of~$\overline{J_\sigma}$. A generalization $J_{\sigma,\tau}$ of $J_\sigma$ is also considered by taking dif\/ferent parameters~$\sigma$ and~$\tau$ on even and odd functions. Then extensions of the above result are proved for $J_{\sigma,\tau}$, where the perturbation has an additional term involving, either the factor $x^{-1}$ on odd functions, or the factor $x$ on even functions. Versions of these results on~${\mathbb{R}}_+$ are derived.}


\Classification{47A55; 47B25; 33C45}


\section{Introduction} \label{s: intro}

The Dunkl operators on ${\mathbb{R}}^n$ were introduced by Dunkl \cite{Dunkl1988,Dunkl1989,Dunkl1991}, and gave rise to what is now called the Dunkl theory \cite{Rosler2003}. They play an important role in physics and stochastic processes (see, e.g., \cite{GenestVinetZhedanov2013,Rosler1998, DiejenVinet2000}). In particular, the Dunkl harmonic oscillators on ${\mathbb{R}}^n$ were studied in \cite{Dunkl2002,NowakStempak2009b,NowakStempak2009a, Rosenblum1994}. We will consider only this operator on ${\mathbb{R}}$, where it is uniquely determined by one parameter. In this case, a conjugation of the Dunkl operator was previously introduced by Yang~\cite{Yang1951} (see also~\cite{Plyushchay2000}).

Let us f\/ix some notation that is used in the whole paper. Let ${\mathcal{S}}={\mathcal{S}}({\mathbb{R}})$ be the Schwartz space on ${\mathbb{R}}$, with its Fr\'echet topology. It decomposes as direct sum of subspaces of even and odd functions, ${\mathcal{S}}={\mathcal{S}}_{\text{\rm ev}}\oplus{\mathcal{S}}_{\text{\rm odd}}$. The even/odd component of a function in ${\mathcal{S}}$ is denoted with the subindex ev/odd. Since ${\mathcal{S}}_{\text{\rm odd}}=x{\mathcal{S}}_{\text{\rm ev}}$, where $x$ is the standard coordinate of ${\mathbb{R}}$, $x^{-1}\phi\in{\mathcal{S}}_{\text{\rm ev}}$ is def\/ined for $\phi\in{\mathcal{S}}_{\text{\rm odd}}$. Let $L^2_\sigma=L^2({\mathbb{R}},|x|^{2\sigma} dx)$ ($\sigma\in{\mathbb{R}}$), whose scalar product and norm are denoted by $\langle\ ,\ \rangle_\sigma$ and $\|\ \|_\sigma$. The above decomposition of ${\mathcal{S}}$ extends to an orthogonal decomposition, $L^2_\sigma=L^2_{\sigma,\text{\rm ev}}\oplus L^2_{\sigma,\text{\rm odd}}$, because the function $|x|^{2\sigma}$ is even. ${\mathcal{S}}$ is a dense subspace of $L^2_\sigma$ if $\sigma>-\frac{1}{2}$, and ${\mathcal{S}}_{\text{\rm odd}}$ is a dense subspace of $L^2_{\tau,\text{\rm odd}}$ if $\tau>-\frac{3}{2}$. Unless otherwise stated, we assume $\sigma>-\frac{1}{2}$ and $\tau>-\frac{3}{2}$. The domain of a (densely def\/ined) operator $P$ in a Hilbert space is denoted by ${\mathsf{D}}(P)$. If $P$ is closable, its closure is denoted by $\overline P$. The domain of a (densely def\/ined) sesquilinear form ${\mathfrak{p}}$ in a Hilbert space is denoted by ${\mathsf{D}}({\mathfrak{p}})$. The quadratic form of ${\mathfrak{p}}$ is also denoted by ${\mathfrak{p}}$. If ${\mathfrak{p}}$ is closable, its closure is denoted by $\bar{\mathfrak{p}}$. For an operator in $L^2_\sigma$ preserving the above decomposition, its restrictions to $L^2_{\sigma,\text{\rm ev/odd}}$ will be indicated with the subindex ev/odd. The operator of multiplication by a continuous function $h$ in $L^2_\sigma$ is also denoted by $h$. The harmonic oscillator is the operator $H=-\frac{d^2}{dx^2}+s^2x^2$ ($s>0$) in $L^2_0$ with ${\mathsf{D}}(H)={\mathcal{S}}$.

The Dunkl operator on ${\mathbb{R}}$ is the operator $T$ in $L^2_\sigma$, with ${\mathsf{D}}(T)={\mathcal{S}}$, determined by $T=\frac{d}{dx}$ on ${\mathcal{S}}_{\text{\rm ev}}$ and $T=\frac{d}{dx}+2\sigma x^{-1}$ on ${\mathcal{S}}_{\text{\rm odd}}$, and the Dunkl harmonic oscillator on ${\mathbb{R}}$ is the operator $J=-T^2+s^2x^2$ in $L^2_\sigma$ with ${\mathsf{D}}(J)={\mathcal{S}}$. Thus $J$ preserves the above decomposition of ${\mathcal{S}}$, being $J_{\text{\rm ev}}=H-2\sigma x^{-1} \frac{d}{dx}$ and $J_{\text{\rm odd}}=H-2\sigma\frac{d}{dx} x^{-1}$. The subindex $\sigma$ is added to $J$ if needed. This $J$ is essentially self-adjoint, and the spectrum of $\overline{J}$ is well known \cite{Rosenblum1994}; in particular, $\overline{J}>0$. In fact, even for $\tau>-\frac{3}{2}$, the operator $J_{\tau,\text{\rm odd}}$ is def\/ined in $L^2_{\tau,\text{\rm odd}}$ with ${\mathsf{D}}(J_{\tau,\text{\rm odd}})={\mathcal{S}}_{\text{\rm odd}}$ because it is a~conjugation of $J_{\tau+1,\text{\rm ev}}$ by a unitary operator (Section~\ref{s: prelim}). Some operators of the form $J+\xi x^{-2}$ ($\xi\in{\mathbb{R}}$) are conjugates of $J$ by powers $|x|^a$ ($a\in{\mathbb{R}}$), and therefore their study can be reduced to the case of~$J$~\cite{AlvCalaza2014}. Our f\/irst theorem analyzes a dif\/ferent perturbation of $J$.

\begin{Theorem}\label{t: UU}
	Let $0<u<1$ and $\xi>0$. If $\sigma>u-\frac{1}{2}$, then there is a positive self-adjoint operator ${\mathcal{U}}$ in $L^2_\sigma$ satisfying the following:
 		\begin{itemize}\itemsep=0pt
 			\item[$(i)$] ${\mathcal{S}}$ is a core of ${\mathcal{U}}^{1/2}$, and, for all $\phi,\psi\in{\mathcal{S}}$,
				\begin{gather}\label{langle UU^1/2 phi, UU^1/2 psi rangle_sigma}
					\big\langle{\mathcal{U}}^{1/2}\phi,{\mathcal{U}}^{1/2}\psi\big\rangle_\sigma
					=\langle J\phi,\psi\rangle_\sigma
					+\xi\big\langle|x|^{-u}\phi,|x|^{-u}\psi\big\rangle_\sigma.
				\end{gather}
 			\item[$(ii)$] ${\mathcal{U}}$ has a discrete spectrum. Let $\lambda_0\le\lambda_1\le\cdots$ be its eigenvalues, repeated according to their multiplicity. There is some $D=D(\sigma,u)>0$, and, for each $\epsilon>0$, there is some $C=C(\epsilon,\sigma,u)>0$ so that, for all $k\in{\mathbb{N}}$,
				\begin{gather}\label{lambda_k, case of UU}
					(2k+1+2\sigma)s+\xi Ds^u(k+1)^{-u}\le\lambda_k
					\le(2k+1+2\sigma)(s+\xi\epsilon s^u)+\xi Cs^u.
				\end{gather}
 		\end{itemize}
\end{Theorem}

\begin{Remark}\label{r: UU}
	In Theorem~\ref{t: UU}, observe the following:
 		\begin{itemize}\itemsep=0pt
		
			\item[(i)] The second term of the right hand side of~\eqref{langle UU^1/2 phi, UU^1/2 psi rangle_sigma} makes sense because $|x|^{-u}{\mathcal{S}}\subset L^2_\sigma$ since $\sigma>u-\frac{1}{2}$.
			
			\item[(ii)] ${\mathcal{U}}=\overline U$, where $U:=J+\xi|x|^{-2u}$ with ${\mathsf{D}}(U)=\bigcap_{m=0}^\infty{\mathsf{D}}({\mathcal{U}}^m)$ (see \cite[Chapter~VI, Section~2.5]{Kato1995}). The more explicit notation $U_\sigma$ will be also used if necessary.
			
			\item[(iii)] The restrictions ${\mathcal{U}}_{\text{ev/odd}}$ are self-adjoint in $L^2_{\sigma,\text{\rm ev/odd}}$ and satisfy~\eqref{langle UU^1/2 phi, UU^1/2 psi rangle_sigma} with $\phi,\psi\in{\mathcal{S}}_{\text{\rm ev/odd}}$ and~\eqref{lambda_k, case of UU} with $k$ even/odd. In fact, by the comments before the statement, ${\mathcal{U}}_{\tau,\text{odd}}$ is def\/ined and satisf\/ies these properties if $\tau>u-\frac{3}{2}$.
			
		\end{itemize}
\end{Remark}

To prove Theorem~\ref{t: UU}, we consider the positive def\/inite symmetric sesquilinear form ${\mathfrak{u}}$ def\/ined by the right hand side of~\eqref{langle UU^1/2 phi, UU^1/2 psi rangle_sigma}. Perturbation theory \cite{Kato1995} is used to show that ${\mathfrak{u}}$ is closable and~$\bar{\mathfrak{u}}$ induces a self-adjoint operator~${\mathcal{U}}$, and to relate the spectra of~${\mathcal{U}}$ and~$\overline J$. Most of the work is devoted to check the conditions to apply this theory so that~\eqref{lambda_k, case of UU} follows; indeed,~\eqref{lambda_k, case of UU} is stronger than a general eigenvalue estimate given by that theory (Remark~\ref{r: type (B)}).

The following generalizations of Theorem~\ref{t: UU} follow with a simple adaptation of the proof. If $\xi<0$, we only have to reverse the inequalities of~\eqref{lambda_k, case of UU}. In~\eqref{langle UU^1/2 phi, UU^1/2 psi rangle_sigma}, we may use a f\/inite sum $\sum_i\xi_i\langle|x|^{-u_i}\phi,|x|^{-u_i}\psi\rangle_\sigma$, where $0<u_i<1$, $\sigma>u_i-\frac{1}{2}$ and $\xi_i>0$; then~\eqref{lambda_k, case of UU} would be modif\/ied by using $\max_iu_i$ and $\min_i\xi_i$ in the left hand side, and $\max_i\xi_i$ in the right hand side. In turn, this can be extended by taking ${\mathbb{R}}^p$-valued functions ($p\in{\mathbb{Z}}_+$), and a f\/inite sum $\sum_i\langle|x|^{-u_i}\Xi_i\phi,|x|^{-u_i}\psi\rangle_\sigma$ in~\eqref{langle UU^1/2 phi, UU^1/2 psi rangle_sigma}, where each $\Xi_i$ is a positive def\/inite self-adjoint endomorphism of ${\mathbb{R}}^p$; then the minimum and maximum eigenvalues of all $\Xi_i$ would be used in~\eqref{lambda_k, case of UU}.

As an open problem, we may ask for a version of Theorem~\ref{t: UU} using Dunkl operators on ${\mathbb{R}}^n$, but we are interested in the following dif\/ferent type of extension. For $\sigma>-\frac{1}{2}$ and $\tau>-\frac{3}{2}$, let $L^2_{\sigma,\tau}=L^2_{\sigma,\text{\rm ev}}\oplus L^2_{\tau,\text{\rm odd}}$, whose scalar product and norm are denoted by $\langle\ ,\ \rangle_{\sigma,\tau}$ and $\|\ \|_{\sigma,\tau}$. Matrix expressions of operators refer to this decomposition. Let $J_{\sigma,\tau}=J_{\sigma,\text{\rm ev}}\oplus J_{\tau,\text{\rm odd}}$ in $L^2_{\sigma,\tau}$, with ${\mathsf{D}}(J_{\sigma,\tau})={\mathcal{S}}$. The hypotheses of the generalization of Theorem~\ref{t: UU} are rather involved to cover enough cases of certain application that will be indicated.

\begin{Theorem}\label{t: VV}
	Let $\xi>0$ and $\eta\in{\mathbb{R}}$, let
		\begin{gather}\label{u, sigma, tau, theta}
			0<u<1,\qquad\sigma>u-\tfrac{1}{2},\qquad\tau>u-\tfrac{3}{2},\qquad\theta>-\tfrac{1}{2},
		\end{gather}
	and set $v=\sigma+\tau-2\theta$. Suppose that the following conditions hold:
		\begin{itemize}\itemsep=0pt
		\item[$(a)$] If $\sigma=\theta\ne\tau$ and $\tau-\sigma\not\in-{\mathbb{N}}$, then
				\begin{gather}\label{VV, a}
					\sigma-1<\tau<\sigma+1,2\sigma+\tfrac{1}{2}.
				\end{gather}
			
			\item[$(b)$] If $\sigma\ne\theta=\tau$ and $\sigma-\tau\not\in-{\mathbb{N}}$, then
				\begin{gather}\label{VV, b}
					-\tau,\tau-1<\sigma<3\tau+1,11\tau+2,\tau+1.
				\end{gather}
			
			\item[$(c)$] If $\sigma\ne\theta=\tau+1$ and $\sigma-\tau-1\not\in-{\mathbb{N}}$, then
				\begin{gather}\label{VV, c}
					\tau+1<\sigma<\tau+3,2\tau+\tfrac{7}{2}.
				\end{gather}
			
			\item[$(d)$] If $\sigma\ne\theta\ne\tau$ and $\sigma-\theta,\tau-\theta\not\in-{\mathbb{N}}$, then
				\begin{gather}\label{VV, d}
\begin{split}
&		\tfrac{\sigma-\tau}{2}-1,\tfrac{\tau-\sigma}{2},\tfrac{\sigma+\tau-1}{4},\tfrac{\sigma+3\tau-2}{14},\tfrac{3\sigma+\tau-4}{14},\tfrac{\sigma+\tau-1}{2}
						<\theta<\tfrac{\sigma+\tau+1}{2},\\
&\tau-1<\sigma<\tau+3.
\end{split}
				\end{gather}
			\end{itemize}
	Then there is a positive self-adjoint operator ${\mathcal{V}}$ in $L^2_{\sigma,\tau}$ satisfying the following:
 		\begin{itemize}\itemsep=0pt
 			\item[$(i)$] ${\mathcal{S}}$ is a core of ${\mathcal{V}}^{1/2}$, and, for all $\phi,\psi\in{\mathcal{S}}$,
				\begin{gather}
					\big\langle{\mathcal{V}}^{1/2}\phi,{\mathcal{V}}^{1/2}\psi\big\rangle_{\sigma,\tau}
					=\langle J_{\sigma,\tau}\phi,\psi\rangle_{\sigma,\tau}
					+\xi\big\langle|x|^{-u}\phi,|x|^{-u}\psi\big\rangle_{\sigma,\tau}\notag\\
					\hphantom{\big\langle{\mathcal{V}}^{1/2}\phi,{\mathcal{V}}^{1/2}\psi\big\rangle_{\sigma,\tau}=}{}
					+\eta\big(\big\langle x^{-1}\phi_{\text{\rm odd}},\psi_{\text{\rm ev}}\big\rangle_\theta
					+\big\langle\phi_{\text{\rm ev}},x^{-1}\psi_{\text{\rm odd}}\big\rangle_\theta\big).
					\label{langle VV^1/2 phi, VV^1/2 psi rangle_sigma,tau}
				\end{gather}

 			\item[$(ii)$] Let $\varsigma_k=\sigma$ if $k$ is even, and $\varsigma_k=\tau$ if $k$ is odd. ${\mathcal{V}}$ has a discrete spectrum. Its eigenvalues form two groups, $\lambda_0\le\lambda_2\le\cdots$ and $\lambda_1\le\lambda_3\le\cdots$, repeated according to their multiplicity, such that there is some $D=D(\sigma,\tau,u)>0$ and, for every $\epsilon>0$, there are some $C=C(\epsilon,\sigma,\tau,u)>0$ and $E=E(\epsilon,\sigma,\tau,\theta)>0$ so that, for all $k\in{\mathbb{N}}$,
				\begin{gather}
					\lambda_k\ge(2k+1+2\varsigma_k)\big(s-2|\eta|\epsilon s^{\frac{v+1}{2}}\big)
					+\xi Ds^u(k+1)^{-u}-2|\eta|Es^{\frac{v+1}{2}},
					\label{lambda_k ge ..., case of VV}\\
					\lambda_k\le(2k+1+2\varsigma_k)\big(s+\epsilon\big(\xi s^u+2|\eta|s^{\frac{v+1}{2}}\big)\big)
					+\xi Cs^u+2|\eta|Es^{\frac{v+1}{2}}.
					\label{lambda_k le ..., case of VV}
				\end{gather}
				
			\item[$(iii)$] Let $\tilde u\in{\mathbb{R}}$ such that
				\begin{gather}\label{tilde u}
					0,v,\tau-2\theta+\tfrac{1}{2},\sigma-2\theta-\tfrac{1}{2}<\tilde u<1,v+1,\sigma+\tfrac{1}{2},\tau+\tfrac{3}{2},
				\end{gather}
			and let $\hat u=\max\{\tilde u,v+1-\tilde u\}$. There is some $D=D(\sigma,\tau,u)>0$ and, for any $\epsilon>0$, there is some $\widetilde C=\widetilde C(\epsilon,\sigma,\tau,u)>0$ so that, for all $k\in{\mathbb{N}}$,
				\begin{gather}\label{lambda_k ge ..., case of VV with tilde u}
					\lambda_k\ge(2k+1+2\varsigma_k)\big(s-|\eta|\epsilon s^{\hat u}\big)+\xi Ds^u(k+1)^{-u}-|\eta|\widetilde Cs^{\hat u}.
				\end{gather}
				
			\item[$(iv)$] If $u=\frac{v+1}{2}$ and $\xi\ge|\eta|$, then there is some $\widetilde D=\widetilde D(\sigma,\tau,u)>0$ so that, for all $k\in{\mathbb{N}}$,
				\begin{gather}\label{lambda_k ge ..., case of VV with u = (v+1)/2, xi ge |eta|}
					\lambda_k\ge(2k+1+2\varsigma_k)s+(\xi-|\eta|)\widetilde Ds^u(k+1)^{-u}.
				\end{gather}
				
			\item[$(v)$] If we add the term $\xi'\langle\phi_{\text{\rm ev}},\psi_{\text{\rm ev}}\rangle_\sigma+\xi''\langle\phi_{\text{\rm odd}},\psi_{\text{\rm odd}}\rangle_\tau$ to the right hand side of~\eqref{langle VV^1/2 phi, VV^1/2 psi rangle_sigma,tau}, for some $\xi',\xi''\in{\mathbb{R}}$, then the result holds as well with the additional term $\max\{\xi',\xi''\}$ in the right hand side of~\eqref{lambda_k le ..., case of VV}, and the additional term, $\xi'$ for $k\in2{\mathbb{N}}$ and $\xi''$ for $k\in2{\mathbb{N}}+1$, in the right hand sides of~\eqref{lambda_k ge ..., case of VV},~\eqref{lambda_k ge ..., case of VV with tilde u} and~\eqref{lambda_k ge ..., case of VV with u = (v+1)/2, xi ge |eta|}.
\end{itemize}
\end{Theorem}
	
\begin{Remark}\label{r: VV}
	Note the following in Theorem~\ref{t: VV}:
		\begin{itemize}\itemsep=0pt
			\item[(i)] Like in Remark~\ref{r: UU}(ii), we have ${\mathcal{V}}=\overline V$, where
				\begin{gather*}
					V=
						\begin{pmatrix}
							U_{\sigma,\text{\rm ev}} & \eta|x|^{2(\theta-\sigma)}x^{-1} \\
							\eta|x|^{2(\theta-\tau)}x^{-1} & U_{\tau,\text{\rm odd}}
						\end{pmatrix},
				\end{gather*}
with ${\mathsf{D}}(V)=\bigcap_{m=0}^\infty{\mathsf{D}}({\mathcal{V}}^m)$. Note that the adjoint of $|x|^{2(\theta-\sigma)}x^{-1}\colon {\mathcal{S}}_{\text{\rm odd}}\to|x|^{2(\theta-\sigma)}{\mathcal{S}}_{\text{\rm ev}}$, as a densely def\/ined operator of $L^2_{\tau,\text{\rm odd}}$ to $L^2_{\sigma,\text{\rm ev}}$, is given by $|x|^{2(\theta-\tau)}x^{-1}$, with the appropriate domain.

			\item[(ii)] Taking $\theta'=\theta-1>-\frac{3}{2}$, since
				\begin{gather*}
					\langle x\phi,\psi\rangle_{\theta'}=\big\langle\phi,x^{-1}\psi\big\rangle_\theta
				\end{gather*}
			for all $\phi\in{\mathcal{S}}_{\text{\rm ev}}$ and $\psi\in{\mathcal{S}}_{\text{\rm odd}}$, we can write~\eqref{langle VV^1/2 phi, VV^1/2 psi rangle_sigma,tau} as
				\begin{gather*}
					\big\langle{\mathcal{V}}^{1/2}\phi,{\mathcal{V}}^{1/2}\psi\big\rangle_{\sigma,\tau}
					=\langle J_{\sigma,\tau}\phi,\psi\rangle_{\sigma,\tau}
					+\xi\big\langle|x|^{-u}\phi,|x|^{-u}\psi\big\rangle_{\sigma,\tau}\\
					\hphantom{\big\langle{\mathcal{V}}^{1/2}\phi,{\mathcal{V}}^{1/2}\psi\big\rangle_{\sigma,\tau}=}{}+\eta\big(\langle\phi_{\text{\rm odd}},x\psi_{\text{\rm ev}}\rangle_{\theta'}
					+\langle x\phi_{\text{\rm ev}},\psi_{\text{\rm odd}}\rangle_{\theta'}\big)
				\end{gather*}
			for all $\phi,\psi\in{\mathcal{S}}$, and, correspondingly,
				\begin{gather*}
					V=
						\begin{pmatrix}
							U_{\sigma,\text{\rm ev}} & \eta|x|^{2(\theta'-\sigma)}x \\
							\eta|x|^{2(\theta'-\tau)}x & U_{\tau,\text{\rm odd}}
						\end{pmatrix}.
				\end{gather*}
			
			\item[(iii)] The conditions~\eqref{VV, a},~\eqref{VV, b} and~\eqref{VV, c} describe three convex open subsets of ${\mathbb{R}}^2$ (Fig.~\ref{fig: a, b, c}). The condition~\eqref{VV, d} describes a convex open subset of ${\mathbb{R}}^3$ (Fig.~\ref{fig: VV, d}), which is symmetric with respect to the plane def\/ined by $\sigma=\tau+1$. It is a ``semi-inf\/inite bar'' with 4 lateral faces, and 5 faces at the ``bounded end.''
			
			\item[(iv)] In Theorem~\ref{t: VV}(iii), the condition~\eqref{tilde u} means that~\eqref{u, sigma, tau, theta} also holds with $\tilde u$ and $v+1-\tilde u$ instead of $u$. There exists $\tilde u$ satisfying~\eqref{tilde u} just when
				\begin{gather}\label{exists tilde u}
					0,v,\tau-2\theta+\tfrac{1}{2},\sigma-2\theta-\tfrac{1}{2}<1,v+1,\sigma+\tfrac{1}{2},\tau+\tfrac{3}{2}.
				\end{gather}
			This property is satisf\/ied in the cases~(b) and~(d) by~\eqref{u, sigma, tau, theta},~\eqref{VV, b} and~\eqref{VV, d}; in particular, we can take $\tilde u=\frac{v+1}{2}$. In the case~(a), if $\tau<3\sigma$, then~\eqref{exists tilde u} holds by~\eqref{u, sigma, tau, theta} and~\eqref{VV, a}. In the case~(c), if $\sigma<3\tau+4$, then~\eqref{exists tilde u} holds by~\eqref{u, sigma, tau, theta} and~\eqref{VV, c}.
\end{itemize}
\end{Remark}

\begin{figure}[th]
\centering
\subfigure[Set def\/ined by~\eqref{VV, a}.]{
\includegraphics[width=3.8cm]{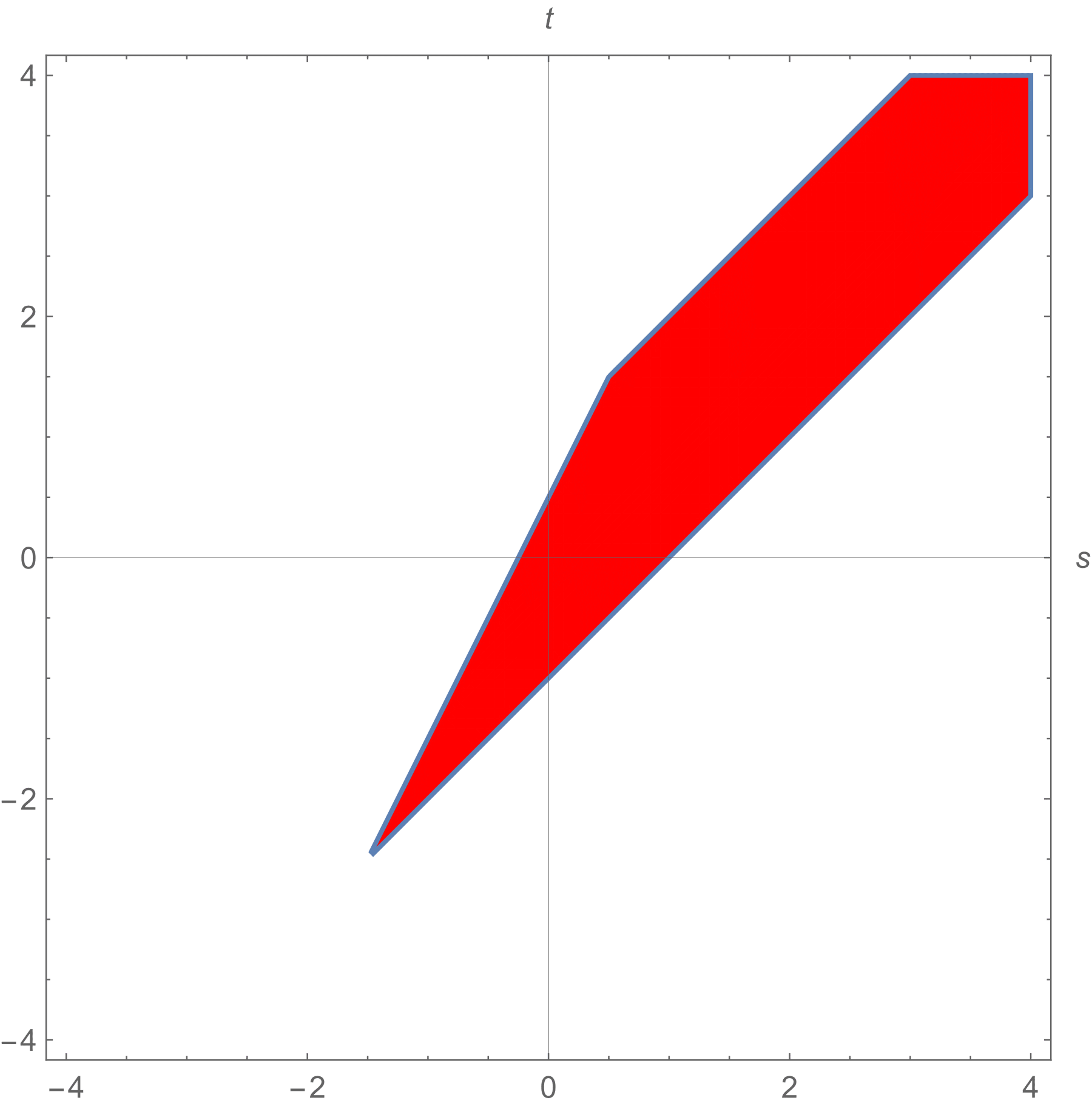}}
\subfigure[Set def\/ined by~\eqref{VV, b}.]{
\includegraphics[width=3.8cm]{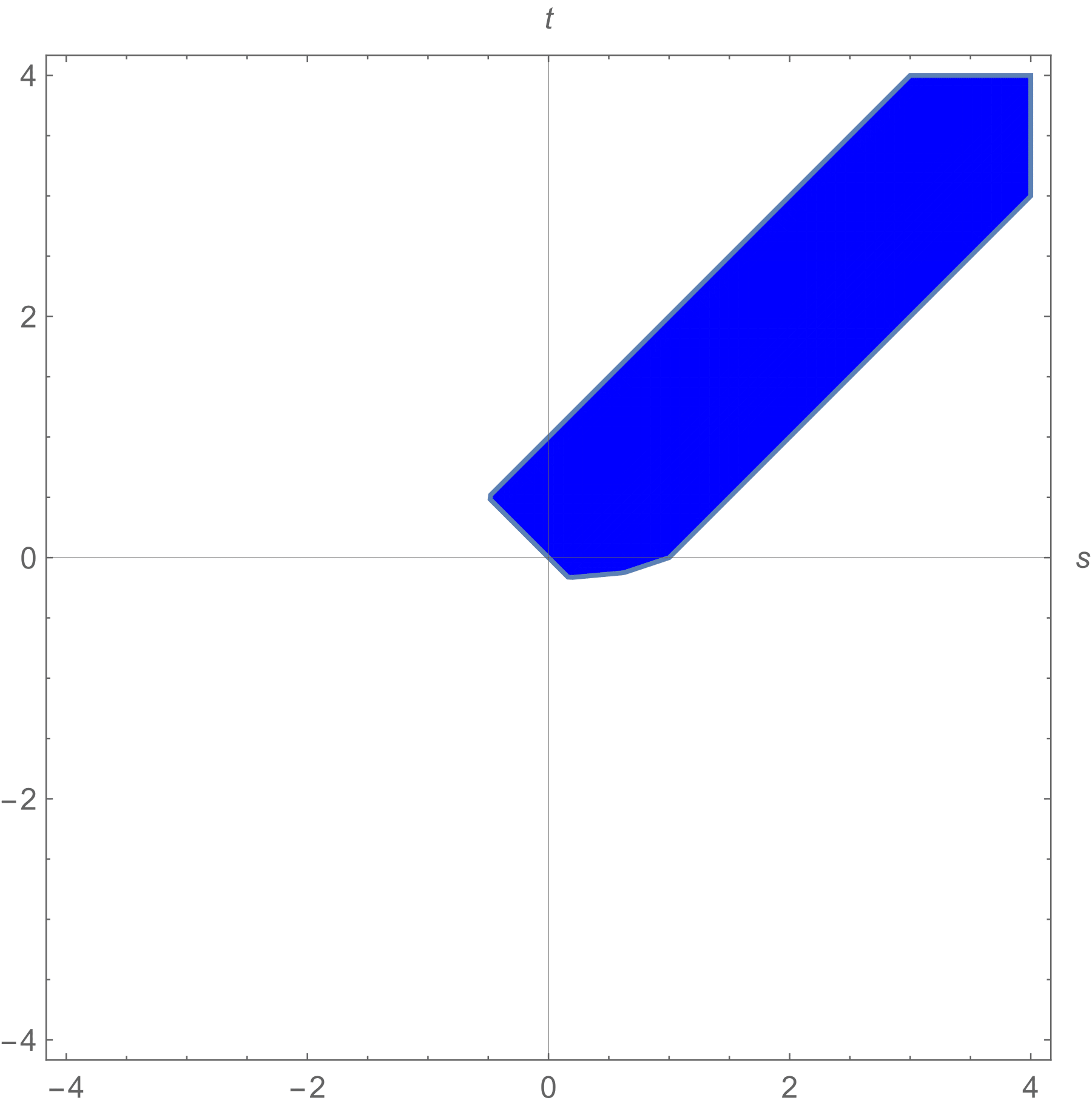}}
\subfigure[Set def\/ined by~\eqref{VV, c}.]{
\includegraphics[width=3.8cm]{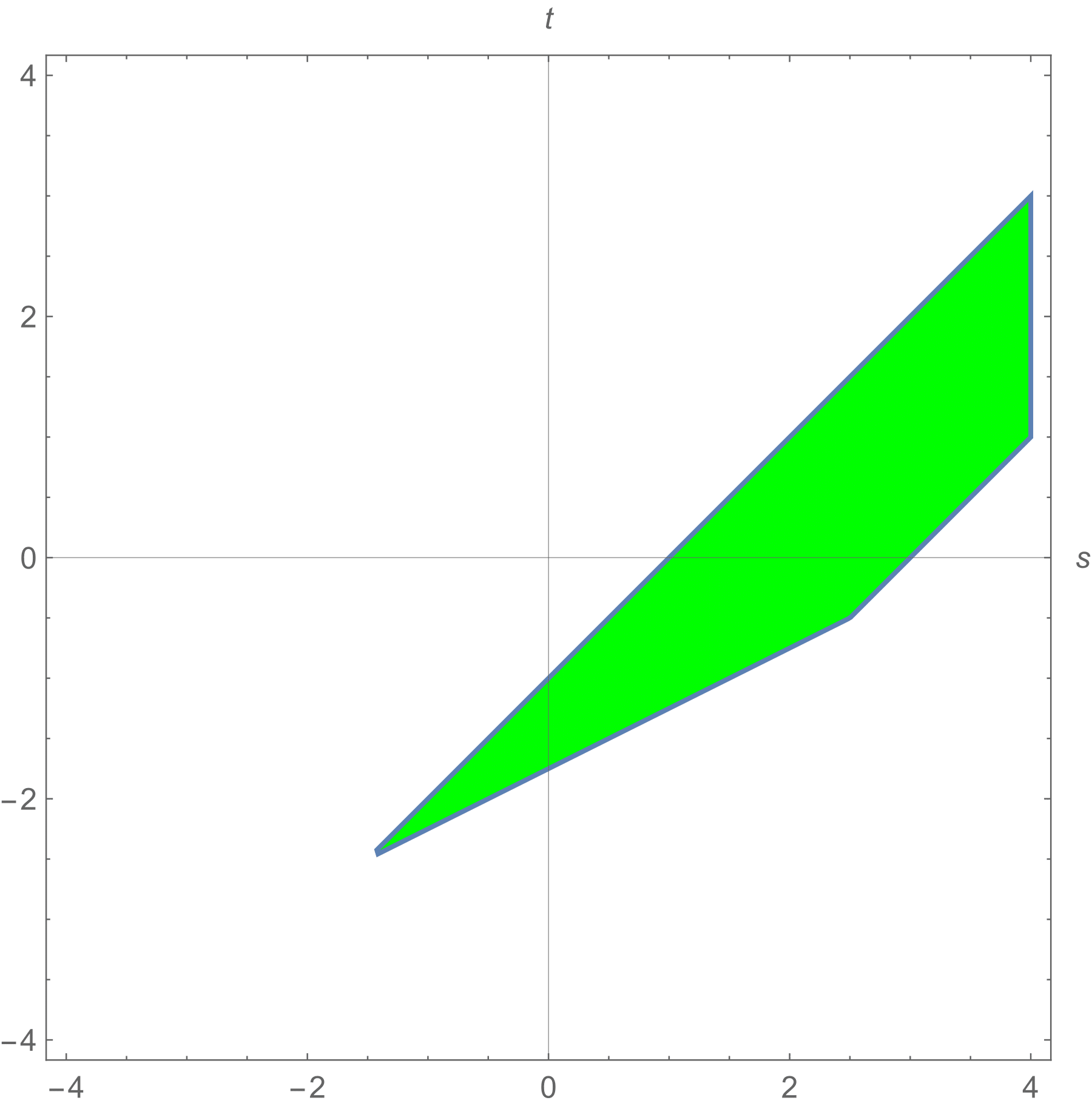}}
\vspace{-2mm}

\caption{Sets in Theorem~\ref{t: VV}(a),(b),(c).}\label{fig: a, b, c}
\end{figure}

\begin{figure}[th]\centering
\includegraphics[height=5.5cm]{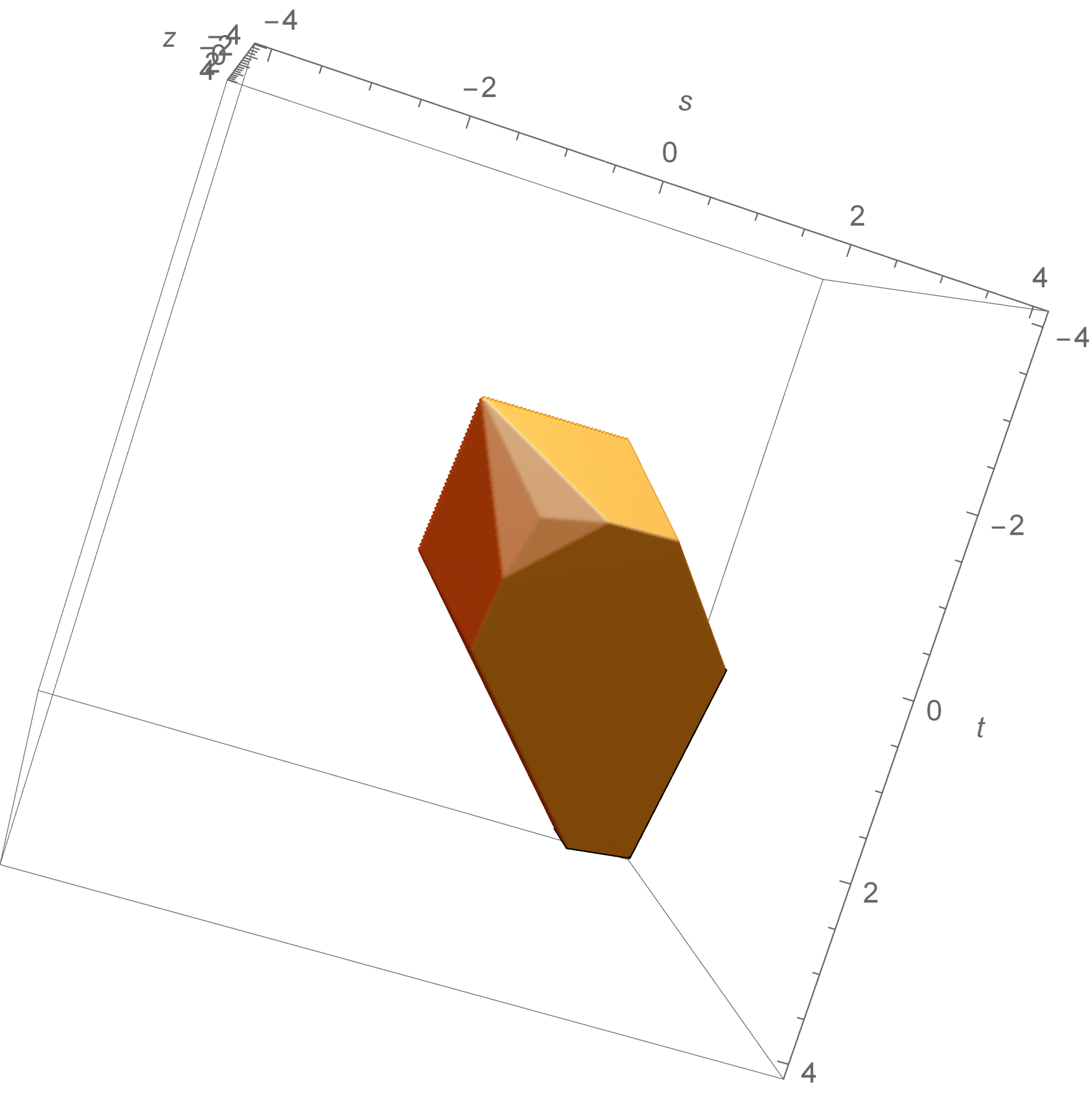}
\vspace{-2mm}

\caption{Set def\/ined by~\eqref{VV, d} in Theorem~\ref{t: VV}(d).}\label{fig: VV, d}
\end{figure}

The main arguments of the proofs of Theorems~\ref{t: UU} and~\ref{t: VV} are given in Sections~\ref{s: ft}--\ref{s: ft'}. But some needed estimates are postponed to Sections~\ref{s: preliminary estimate} and~\ref{s: main estimates} because they are of rather independent nature, and with rather long and tedious proofs.

Versions of these results on ${\mathbb{R}}_+$ are also derived in Section~\ref{s: R_+} (Corollaries~\ref{c: PP},~\ref{c: QQ} and~\ref{c: WW}). In~\cite{AlvCalazaFranco:Witten-general}, these corollaries are used to study a version of the Witten's perturbation $\Delta_s$ of the Laplacian on strata with the general adapted metrics of \cite{BrasseletHectorSaralegi1992,Nagase1983,Nagase1986}. This gives rise to an analytic proof of Morse inequalities in strata involving intersection homology of arbitrary perversity, which was our original motivation. The simplest case of adapted metrics, corresponding to the lower middle perversity, was treated in~\cite{AlvCalaza:Witten} using an operator induced by~$J$ on~${\mathbb{R}}_+$. The perturbations of $J$ studied here show up in the local models of~$\Delta_s$ when general adapted metrics are considered. Some details of this application are given in Section~\ref{s: Witten}.

\section{Preliminaries}\label{s: prelim}

The {\em Dunkl annihilation} and {\em creation} operators are $B=sx+T$ and $B'=sx-T$ ($s>0$). Like~$J$, the operators $B$ and $B'$ are considered in $L^2_\sigma$ with domain ${\mathcal{S}}$. They are perturbations of the usual annihilation and creation operators. The operators $T$, $B$, $B'$ and $J$ are continuous on ${\mathcal{S}}$. The following properties hold~\cite{AlvCalaza2014, Rosenblum1994}:
	\begin{itemize}\itemsep=0pt
	
		\item $B'$ is adjoint of $B$, and $J$ is essentially self-adjoint.
		
		\item The spectrum of $\overline J$ consists of the eigenvalues\footnote{It is assumed that $0\in{\mathbb{N}}$.} $(2k+1+2\sigma)s$, $k\in{\mathbb{N}}$, of multiplicity one.
		
		\item The corresponding normalized eigenfunctions $\phi_k$ are inductively def\/ined by
 			\begin{gather}
 				\phi_0 =s^{(2\sigma+1)/4}\Gamma\big(\sigma+\tfrac{1}{2}\big)^{-\frac{1}{2}}e^{-sx^2/2} ,\label{phi_0}\\
 				\phi_k =
 					\begin{cases}
 						(2ks)^{-\frac{1}{2}}B'\phi_{k-1} & \text{if $k$ is even},\\
 						(2(k+2\sigma)s)^{-\frac{1}{2}}B'\phi_{k-1} & \text{if $k$ is odd},
 					\end{cases}
 				\qquad k\ge1.\label{phi_k}
 			\end{gather}
		\item The eigenfunctions $\phi_k$ also satisfy
 			\begin{gather}
 				B\phi_0 =0 ,\label{B phi_0}\\
 				B\phi_k =
 					\begin{cases}
 					(2ks)^{\frac{1}{2}}\phi_{k-1} & \text{if $k$ is even},\\
 					(2(k+2\sigma)s)^{\frac{1}{2}}\phi_{k-1} & \text{if $k$ is odd},
 					\end{cases}
 				\qquad k\ge1 .\label{B phi_k}
 			\end{gather}
		\item $\bigcap_{m=0}^\infty{\mathsf{D}}\big(\overline J^m\big)={\mathcal{S}}$.
	\end{itemize}
By~\eqref{phi_0} and~\eqref{phi_k}, we get $\phi_k=p_ke^{-sx^2/2}$, where $p_k$ is the sequence of polynomials inductively given by $p_0=s^{(2\sigma+1)/4}\Gamma(\sigma+\frac{1}{2})^{-\frac{1}{2}}$ and
 	\begin{gather*}
 		p_k=
 			\begin{cases}
 				(2ks)^{-\frac{1}{2}}(2sxp_{k-1}-Tp_{k-1}) & \text{if $k$ is even},\\
 				(2(k+2\sigma)s)^{-\frac{1}{2}}(2sxp_{k-1}-Tp_{k-1}) & \text{if $k$ is odd},
 			\end{cases}
 		\qquad k\ge1 .
 	\end{gather*}
Up to normalization, $p_k$ is the sequence of generalized Hermite polynomials \cite[p.~380, Problem~25]{Szego1975}, and $\phi_k$ is the sequence of generalized Hermite functions. Each $p_k$ is of degree $k$, even/odd if $k$ is even/odd, and with positive leading coef\/f\/icient. They satisfy the recursion formula \cite[equation~(13)]{AlvCalaza2014}
 	\begin{gather}\label{recursion}
 		p_k=
 			\begin{cases}
 			k^{-\frac{1}{2}}\big((2s)^{\frac{1}{2}}xp_{k-1}-(k-1+2\sigma)^{\frac{1}{2}}p_{k-2}\big) & \text{if $k$ is even},\\
 			(k+2\sigma)^{-\frac{1}{2}}\big((2s)^{\frac{1}{2}}xp_{k-1}-(k-1)^{\frac{1}{2}}p_{k-2}\big) & \text{if $k$ is odd} .
 			\end{cases}
 	\end{gather}
When $k=2m+1$ ($m\in{\mathbb{N}}$), we have \cite[equation~(14)]{AlvCalaza2014}
 \begin{gather} \label{x^-1 p_k}
 x^{-1}p_k=\sum_{i=0}^m(-1)^{m-i}
 \sqrt{\frac{m! \Gamma(i+\frac{1}{2}+\sigma)s}{i! \Gamma(m+\frac{3}{2}+\sigma)}} p_{2i} .
 \end{gather}
The Pochhammer symbol could be used to simplify this expression, as well as many other expressions in Sections~\ref{s: ft} and~\ref{s: mixed}. However there are quotients of gamma functions in Sections~\ref{s: mixed} and~\ref{s: ft'} that can not be simplif\/ied in this way (see e.g.\ Proposition~\ref{p: hat c_k,ell, case where sigma = theta ne tau}). Thus, for the sake of uniformity, we use gamma functions in all quotients of this type.

Let ${\mathfrak{j}}$ be the positive def\/inite symmetric sesquilinear form in $L^2_\sigma$, with ${\mathsf{D}}({\mathfrak{j}})={\mathcal{S}}$, given by ${\mathfrak{j}}(\phi,\psi)=\langle J\phi,\psi\rangle_\sigma$. Like in the case of~$J$, the subindex $\sigma$ will be added to the notation $T$, $B$, $B'$ and $\phi_k$ and ${\mathfrak{j}}$ if necessary. Observe that
	\begin{gather}
		B_\sigma =
			\begin{cases}
				B_\tau & \text{on ${\mathcal{S}}_{\text{\rm ev}}$},\\
				B_\tau+2(\sigma-\tau)x^{-1} & \text{on ${\mathcal{S}}_{\text{\rm odd}}$} ,
			\end{cases}
		\label{B_sigma}\\
		B'_\sigma =
			\begin{cases}
				B'_\tau & \text{on ${\mathcal{S}}_{\text{\rm ev}}$},\\
				B'_\tau+2(\tau-\sigma)x^{-1} & \text{on ${\mathcal{S}}_{\text{\rm odd}}$} .
			\end{cases}
		\label{B'_sigma}
	\end{gather}
	
The operator $x\colon {\mathcal{S}}_{\text{\rm ev}}\to{\mathcal{S}}_{\text{\rm odd}}$ is a homeomorphism \cite{AlvCalaza2014}, which extends to a unitary operator $x\colon L^2_{\sigma,\text{\rm ev}}\to L^2_{\sigma-1,\text{\rm odd}}$. We get $x J_{\sigma,\text{\rm ev}} x^{-1}=J_{\sigma-1,\text{\rm odd}}$ because $x \big[\frac{d^2}{dx^2},x^{-1}\big]=-2\frac{d}{dx} x^{-1}$. Thus, even for any $\tau>-\frac{3}{2}$, the operator $J_{\tau,\text{\rm odd}}$ is densely def\/ined in $L^2_{\tau,\text{\rm odd}}$, with ${\mathsf{D}}(J_{\tau,\text{\rm odd}})={\mathcal{S}}_{\text{\rm odd}}$, and has the same spectral properties as $J_{\tau+1,\text{\rm ev}}$; in particular, the eigenvalues of $\overline{J_{\tau,\text{\rm odd}}}$ are $(2k+1+2\tau)s$ ($k\in2{\mathbb{N}}+1$), and $\phi_{\tau,k}=x\phi_{\tau+1,k-1}$.

 To prove the results of the paper, alternative arguments could be given by using the expression of the generalized Hermite polynomials in terms of the Laguerre ones (see, e.g., \cite[p.~525]{Rosler1998} or \cite[p.~23]{Rosler2003}). In particular, certain asymptotic estimates of Laguerre functions \cite{Erdelyi1960,Muckenhoupt1970b} (see also \cite{AskeyWainger1965,Muckenhoupt1970a-II}), yield the following asymptotic estimates of the generalized Hermite functions \cite[Section~2.4]{AlvCalaza:generalizedHermite}: there are some $C,c>0$, depending only on $\sigma$, such that
 \begin{gather}\label{|phi_k(x)| |x|^sigma le ...}
 |\phi_k(x)x^\sigma|\le
 \begin{cases}
 Cs^{\frac{\bar\sigma}{2}+\frac{1}{4}} x^{\bar\sigma}\nu^{\frac{\bar\sigma}{2}-\frac{1}{4}} & \text{if $0<x\le\sqrt{\frac{1}{s\nu}}$,}\\
 Cs^{\frac{1}{4}}\nu^{-\frac{1}{4}} & \text{if $\sqrt{\frac{1}{s\nu}}<x\le\sqrt{\frac{\nu}{2s}}$,}\\
 Cs^{\frac{1}{4}}(\nu^{\frac{1}{3}}+|sx^2-\nu|)^{-\frac{1}{4}} & \text{if $\sqrt{\frac{\nu}{2s}}<x\le\sqrt{\frac{3\nu}{2s}}$,}\\
 C(sx)^{\frac{1}{2}}e^{-csx^2} & \text{if $\sqrt{\frac{3\nu}{2s}}<x$,}
 \end{cases}
\end{gather}
where $\bar\sigma=\bar\sigma_k=\sigma+\frac{1-(-1)^k}{2}$ and $\nu=\nu_k=2k+1+2\sigma$, with the proviso that we must take $\nu=2$ if $k=0$ and $\sigma<\frac{1}{2}$.

\section[The sesquilinear form ${\mathfrak{t}}$]{The sesquilinear form $\boldsymbol{{\mathfrak{t}}}$}\label{s: ft}

Let $0<u<1$ such that $\sigma>u-\frac{1}{2}$. Then $|x|^{-u}{\mathcal{S}}\subset L^2_\sigma$, and therefore a positive def\/inite symmetric sesquilinear form ${\mathfrak{t}}$ in $L^2_\sigma$, with ${\mathsf{D}}({\mathfrak{t}})={\mathcal{S}}$, is def\/ined by
	\begin{gather*}
		{\mathfrak{t}}(\phi,\psi)=\big\langle|x|^{-u}\phi,|x|^{-u}\psi\big\rangle_\sigma=\langle\phi,\psi\rangle_{\sigma-u} .
	\end{gather*}
The notation ${\mathfrak{t}}_\sigma$ may be also used. The goal of this section is to study~${\mathfrak{t}}$ and apply it to prove Theorem~\ref{t: UU}. Precisely, an estimation of the values ${\mathfrak{t}}(\phi_k,\phi_\ell)$ is needed.

\begin{Lemma}\label{l: ft'(B' phi,psi), ft'(phi,B' psi)}
	For all $\phi\in{\mathcal{S}}_{\text{\rm odd}}$ and $\psi\in{\mathcal{S}}_{\text{\rm ev}}$,
		\begin{gather*}
{\mathfrak{t}}(B'\phi,\psi)-{\mathfrak{t}}(\phi,B\psi)={\mathfrak{t}}(\phi,B'\psi)-{\mathfrak{t}}(B\phi,\psi)=-2u {\mathfrak{t}}\big(x^{-1}\phi,\psi\big).
		\end{gather*}
\end{Lemma}

\begin{proof}
	By~\eqref{B_sigma} and~\eqref{B'_sigma}, for all $\phi\in{\mathcal{S}}_{\text{\rm odd}}$ and $\psi\in{\mathcal{S}}_{\text{\rm ev}}$,
		\begin{gather*}
			{\mathfrak{t}}(B'_\sigma\phi,\psi)-{\mathfrak{t}}(\phi,B_\sigma\psi)
			 =\langle B'_{\sigma-u}\phi,\psi\rangle_{\sigma-u}-2u\big\langle x^{-1}\phi,\psi\big\rangle_{\sigma-u}
			-\langle\phi,B_{\sigma-u}\psi\rangle_{\sigma-u}\\
			\hphantom{{\mathfrak{t}}(B'_\sigma\phi,\psi)-{\mathfrak{t}}(\phi,B_\sigma\psi)}{}
=-2u {\mathfrak{t}}\big(x^{-1}\phi,\psi\big) ,\\
			{\mathfrak{t}}(\phi,B'_\sigma\psi)-{\mathfrak{t}}(B_\sigma\phi,\psi)
			 =\langle\phi,B'_{\sigma-u}\psi\rangle_{\sigma-u}-\langle B_{\sigma-u}\phi,\psi\rangle_{\sigma-u}
			-2u\big\langle x^{-1}\phi,\psi\big\rangle_{\sigma-u}\\
			\hphantom{{\mathfrak{t}}(\phi,B'_\sigma\psi)-{\mathfrak{t}}(B_\sigma\phi,\psi)}{}
			 =-2u {\mathfrak{t}}\big(x^{-1}\phi,\psi\big) .\tag*{\qed}
		\end{gather*}
\renewcommand{\qed}{}
\end{proof}

In the whole of this section, $k$, $\ell$, $m$, $n$, $i$, $j$, $p$ and $q$ will be natural numbers. Let $c_{k,\ell}={\mathfrak{t}}(\phi_k,\phi_\ell)$ and $d_{k,\ell}=c_{k,\ell}/c_{0,0}$. Thus $d_{k,\ell}=d_{\ell,k}$, and $d_{k,\ell}=0$ when $k+\ell$ is odd. Since
	\begin{gather}\label{int_-infty^infty e^-sx^2 |x|^2 kappa dx}
		\int_{-\infty}^\infty e^{-sx^2}|x|^{2\kappa} dx=s^{-\frac{2\kappa+1}{2}}\Gamma\big(\kappa+\tfrac{1}{2}\big)
	\end{gather}
for $\kappa>-\frac{1}{2}$, we get
	\begin{gather}\label{c_0,0}
		c_{0,0}=\Gamma\big(\sigma-u+\tfrac{1}{2}\big)\Gamma\big(\sigma+\tfrac{1}{2}\big)^{-1}s^u.
	\end{gather}

\begin{Lemma}\label{l: d_k,0, k=2n>0}
	If $k=2m>0$, then
		\begin{gather*}
			d_{k,0}=\frac{u}{\sqrt{m}}\sum_{j=0}^{m-1}(-1)^{m-j}
			\sqrt{\frac{(m-1)! \Gamma(j+\frac{1}{2}+\sigma)}{j! \Gamma(m+\frac{1}{2}+\sigma)}} d_{2j,0} .
		\end{gather*}
\end{Lemma}

\begin{proof}
	By~\eqref{phi_k},~\eqref{B phi_0},~\eqref{x^-1 p_k} and Lemma~\ref{l: ft'(B' phi,psi), ft'(phi,B' psi)},
		\begin{gather*}
			c_{k,0}=\frac{1}{\sqrt{2sk}} {\mathfrak{t}}(B'\phi_{k-1},\phi_0)
		=\frac{1}{\sqrt{2sk}} {\mathfrak{t}}(\phi_{k-1},B\phi_0)
			-\frac{2u}{\sqrt{2sk}} {\mathfrak{t}}\big(x^{-1}\phi_{k-1},\phi_0\big)\\
	\hphantom{c_{k,0}}{}
			=-\frac{2u}{\sqrt{2sk}} {\mathfrak{t}}\big(x^{-1}\phi_{k-1},\phi_0\big)
			=\frac{u}{\sqrt{m}}\sum_{j=0}^{m-1}(-1)^{m-j}\sqrt{\frac{(m-1)!
			\Gamma(j+\frac{1}{2}+\sigma)}{j! \Gamma(m+\frac{1}{2}+\sigma)}} c_{2j,0} .\tag*{\qed}
		\end{gather*}
\renewcommand{\qed}{}
\end{proof}

\begin{Lemma}\label{l: d_k,ell, k=2m>0, ell = 2n>0}
	If $k=2m>0$ and $\ell=2n>0$, then
		\begin{gather*}
			d_{k,\ell}=\sqrt{\frac{m}{n}} d_{k-1,\ell-1}+\frac{u}{\sqrt{n}}\sum_{j=0}^{n-1}(-1)^{n-j}
			\sqrt{\frac{(n-1)! \Gamma(j+\frac{1}{2}+\sigma)}{j! \Gamma(n+\frac{1}{2}+\sigma)}} d_{k,2j} .
		\end{gather*}
\end{Lemma}

\begin{proof}
	By~\eqref{phi_k},~\eqref{B phi_k},~\eqref{x^-1 p_k} and Lemma~\ref{l: ft'(B' phi,psi), ft'(phi,B' psi)},
		\begin{gather*}
			c_{k,\ell}=\frac{1}{\sqrt{2s\ell}} {\mathfrak{t}}(\phi_k,B'\phi_{\ell-1})
			=\frac{1}{\sqrt{2\ell s}} {\mathfrak{t}}(B\phi_k,\phi_{\ell-1})
			-\frac{2u}{\sqrt{2\ell s}} {\mathfrak{t}}\big(\phi_k,x^{-1}\phi_{\ell-1}\big)\\
\hphantom{c_{k,\ell}}{}
			=\sqrt{\frac{m}{n}} c_{k-1,\ell-1}
			+\frac{u}{\sqrt{n}} \sum_{j=0}^{n-1}(-1)^{n-j}
			\sqrt{\frac{(n-1)! \Gamma(j+\frac{1}{2}+\sigma)}
			{j! \Gamma(n+\frac{1}{2}+\sigma)}} c_{k,2j} .\tag*{\qed}
		\end{gather*}
\renewcommand{\qed}{}
\end{proof}

\begin{Lemma}\label{l: d_k,ell, k=2m+1, ell = 2n+1}
	If $k=2m+1$ and $\ell=2n+1$, then
		\begin{gather*}
			d_{k,\ell}=\sqrt{\frac{n+\frac{1}{2}+\sigma}{m+\frac{1}{2}+\sigma}} d_{k-1,\ell-1}
			-\frac{u}{\sqrt{m+\frac{1}{2}+\sigma}} \sum_{j=0}^n(-1)^{n-j}
			\sqrt{\frac{n! \Gamma(j+\frac{1}{2}+\sigma)}{j! \Gamma(n+\frac{3}{2}+\sigma)}} d_{k-1,2j} .
		\end{gather*}
\end{Lemma}

\begin{proof}
	By~\eqref{phi_k},~\eqref{B phi_k},~\eqref{x^-1 p_k} and Lemma~\ref{l: ft'(B' phi,psi), ft'(phi,B' psi)},
		\begin{gather*}
			c_{k,\ell} =\frac{1}{\sqrt{2(k+2\sigma)s}} {\mathfrak{t}}(B'\phi_{k-1},\phi_\ell) \\
\hphantom{c_{k,\ell}}{}
			 =\frac{1}{\sqrt{2(k+2\sigma)s}} {\mathfrak{t}}(\phi_{k-1},B\phi_\ell)
			-\frac{2u}{\sqrt{2(k+2\sigma)s}} {\mathfrak{t}}\big(\phi_{k-1},x^{-1}\phi_\ell\big)\\
\hphantom{c_{k,\ell}}{}
			=\sqrt{\frac{n+\frac{1}{2}+\sigma}{m+\frac{1}{2}+\sigma}} c_{k-1,\ell-1}
 -\frac{u}{\sqrt{m+\frac{1}{2}+\sigma}}\sum_{j=0}^n(-1)^{n-j}
			\sqrt{\frac{n! \Gamma(j+\frac{1}{2}+\sigma)}
			{j!\Gamma(n+\frac{3}{2}+\sigma)}} c_{k-1,2j} .\tag*{\qed}
		\end{gather*}
\renewcommand{\qed}{}
\end{proof}

The following def\/initions are given for $k\ge\ell$ with $k+\ell$ even. Let	
	\begin{gather}\label{Pi_k,ell, case k = 2m ge ell = 2n}
		\Pi_{k,\ell}=\sqrt{\frac{m! \Gamma(n+\frac{1}{2}+\sigma)}{n! \Gamma(m+\frac{1}{2}+\sigma)}}
	\end{gather}
if $k=2m\ge\ell=2n$, and
	\begin{gather}\label{Pi_k,ell, case k = 2m+1 ge ell = 2n+1}
		\Pi_{k,\ell}=\sqrt{\frac{m! \Gamma(n+\frac{3}{2}+\sigma)}{n! \Gamma(m+\frac{3}{2}+\sigma)}}
	\end{gather}
if $k=2m+1\ge\ell=2n+1$. Let $\Sigma_{k,\ell}$ be inductively def\/ined as follows\footnote{We use the convention that a product of an empty set of factors is~$1$. Such empty products are possible in~\eqref{Sigma_k,0} (when $m=0$), in Lemma~\ref{l: Weierstrass} and its proof, and in the proofs of Lemma~\ref{l: Sigma_k,ell le C' m^-u(1-u) (m-n+1)^-(1-u)^2} and Remark~\ref{r: ft(phi_k) ge D s^u/2 (k+1)^-u}. Consistently, the sum of an empty set of terms is $0$. Such empty sums are possible in Lemma~\ref{l: hat d_k,ell, case k = 2m+1, ell = 2n+1, and sigma = theta ne tau} and its proof, and in the proof of Proposition~\ref{p: hat c_k,ell, case where sigma = theta ne tau}.}:
	\begin{gather}\label{Sigma_k,0}
		\Sigma_{k,0}=\prod_{i=1}^m\left(1-\frac{1-u}{i}\right)=\frac{\Gamma(m+u)}{m!\Gamma(u)}
	\end{gather}
if $k=2m$;
	\begin{gather}\label{Sigma_k,ell, k = 2m ge ell = 2n}
		\Sigma_{k,\ell}=\Sigma_{k-1,\ell-1}
		+u\sum_{j=0}^{n-1}\frac{(n-1)! \Gamma(j+\frac{1}{2}+\sigma)}{j! \Gamma(n+\frac{1}{2}+\sigma)}
		\Sigma_{k,2j}
	\end{gather}
if $k=2m\ge\ell=2n>0$; and
	\begin{gather}
		\Sigma_{k,\ell} =\Sigma_{k-1,\ell-1}
		-u\sum_{j=0}^n\frac{n! \Gamma(j+\frac{1}{2}+\sigma)}{j! \Gamma(n+\frac{3}{2}+\sigma)} \Sigma_{k-1,2j}
		\label{Sigma_k,ell, k = 2m+1 ge ell = 2n+1}\\
		\hphantom{\Sigma_{k,\ell}}{} =\left(1-\frac{u}{n+\frac{1}{2}+\sigma}\right)\Sigma_{k-1,\ell-1}
		 -\frac{nu}{n+\frac{1}{2}+\sigma} \sum_{j=0}^{n-1}
		\frac{(n-1)! \Gamma(j+\frac{1}{2}+\sigma)}{j! \Gamma(n+\frac{1}{2}+\sigma)} \Sigma_{k-1,2j}
		\label{Sigma_k,ell, k = 2m+1, ell = 2n+1, 2nd}
	\end{gather}
if $k=2m+1\ge\ell=2n+1$. Thus $\Sigma_{0,0}=1$, $\Sigma_{2,0}=u$, $\Sigma_{4,0}=\frac{1}{2}u(1+u)$, and
	\begin{gather}\label{Sigma_k,1}
		\Sigma_{k,1}=\left(1-\frac{u}{\frac{1}{2}+\sigma}\right)\Sigma_{k-1,0}
	\end{gather}
if $k$ is odd. From~\eqref{Sigma_k,0} and using induction on $m$, it easily follows that
	\begin{gather}\label{Sigma_k,0 = frac u m sum_j=0^m-1 Sigma_2j,0}
		\Sigma_{k,0}=\frac{u}{m}\sum_{j=0}^{m-1}\Sigma_{2j,0}
	\end{gather}
for $k=2m>0$. Combining~\eqref{Sigma_k,ell, k = 2m ge ell = 2n} with~\eqref{Sigma_k,ell, k = 2m+1 ge ell = 2n+1}, and~\eqref{Sigma_k,ell, k = 2m+1, ell = 2n+1, 2nd} with~\eqref{Sigma_k,ell, k = 2m ge ell = 2n}, we get
	\begin{gather}\label{Sigma_k,ell = Sigma_k-2,ell-2 - ..., k = 2m ge ell = 2n > 0}
		\Sigma_{k,\ell}=\Sigma_{k-2,\ell-2}
		-u\sum_{j=0}^{n-1}\frac{(n-1)!\Gamma(j+\frac{1}{2}+\sigma)}{j!\Gamma(n+\frac{1}{2}+\sigma)}(\Sigma_{k-2,2j}-\Sigma_{k,2j})
	\end{gather}
if $k=2m\ge\ell=2n>0$; and
	\begin{gather}
		\Sigma_{k,\ell}=\left(1-\frac{u}{n+\frac{1}{2}+\sigma}\right)\Sigma_{k-2,\ell-2}\nonumber\\
\hphantom{\Sigma_{k,\ell}=}{}
		+\left(1-\frac{u+n}{n+\frac{1}{2}+\sigma}\right)u\sum_{j=0}^{n-1}
		\frac{(n-1)!\Gamma(j+\frac{1}{2}+\sigma)}{j!\Gamma(n+\frac{1}{2}+\sigma)}\Sigma_{k-1,2j}\label{Sigma_k,ell = Sigma_k-2,ell-2 + ..., k = 2m+1 ge ell = 2n+1 > 1}
	\end{gather}
if $k=2m+1\ge\ell=2n+1>1$.

\begin{Proposition}\label{p: d_k,ell}
	$d_{k,\ell}=(-1)^{m+n}\Pi_{k,\ell}\Sigma_{k,\ell}$ if $k=2m\ge\ell=2n$, or if $k=2m+1\ge\ell=2n+1$.
\end{Proposition}

\begin{proof}
	We proceed by induction on $k$ and $l$. The statement is obvious for $k=\ell=0$ because $d_{0,0}=\Pi_{0,0}=\Sigma_{0,0}=1$.
	
	Let $k=2m>0$, and assume that the result is true for all $d_{2j,0}$ with $j<m$. Then, by Lemma~\ref{l: d_k,0, k=2n>0},~\eqref{Pi_k,ell, case k = 2m ge ell = 2n} and~\eqref{Sigma_k,0 = frac u m sum_j=0^m-1 Sigma_2j,0},
		\begin{gather*}
			d_{k,0} =\frac{u}{\sqrt{m}}\sum_{j=0}^{m-1}(-1)^{m-j}
			\sqrt{\frac{(m-1)! \Gamma(j+\frac{1}{2}+\sigma)}{j! \Gamma(m+\frac{1}{2}+\sigma)}}
			(-1)^j\Pi_{2j,0}\Sigma_{2j,0}\\
\hphantom{d_{k,0}}{}
=(-1)^m \frac{u}{\sqrt{m}}\sum_{j=0}^{m-1}
			\sqrt{\frac{(m-1)! \Gamma(j+\frac{1}{2}+\sigma)}{j! \Gamma(m+\frac{1}{2}+\sigma)}}
			\sqrt{\frac{j! \Gamma(\frac{1}{2}+\sigma)}{\Gamma(j+\frac{1}{2}+\sigma)}}\Sigma_{2j,0}\\
\hphantom{d_{k,0}}{} =(-1)^m\Pi_{k,0} \frac{u}{m}\sum_{j=0}^{m-1}\Sigma_{2j,0}
			=(-1)^m\Pi_{k,0}\Sigma_{k,0} .
		\end{gather*}
	
	Now, take $k=2m\ge\ell=2n>0$ so that the equality of the statement holds for $d_{k-1,\ell-1}$ and all $d_{k,2j}$ with $j<n$. Then, by Lemma~\ref{l: d_k,ell, k=2m>0, ell = 2n>0},
		\begin{gather*}
			d_{k,\ell} =\sqrt{\frac{m}{n}} (-1)^{m+n}\Pi_{k-1,\ell-1}\Sigma_{k-1,\ell-1}\\
\phantom{d_{k,\ell}=}{}+\frac{u}{\sqrt{n}}\sum_{j=0}^{n-1}(-1)^{n-j}
			\sqrt{\frac{(n-1)! \Gamma(j+\frac{1}{2}+\sigma)}{j! \Gamma(n+\frac{1}{2}+\sigma)}}
			(-1)^{m+j}\Pi_{k,2j}\Sigma_{k,2j} .
		\end{gather*}
	Here, by~\eqref{Pi_k,ell, case k = 2m ge ell = 2n} and~\eqref{Pi_k,ell, case k = 2m+1 ge ell = 2n+1}, $\sqrt{m/n} \Pi_{k-1,\ell-1}=\Pi_{k,\ell}$, and{\samepage
		\begin{gather*}
			\frac{1}{\sqrt{n}}
			\sqrt{\frac{(n-1)! \Gamma(j+\frac{1}{2}+\sigma)}{j! \Gamma(n+\frac{1}{2}+\sigma)}} \Pi_{k,2j}
					 =\frac{1}{\sqrt{n}}
					\sqrt{\frac{m! \Gamma(n+\frac{1}{2}+\sigma)}{(n-1)!
					\Gamma(m+\frac{1}{2}+\sigma)}}
					\frac{(n-1)! \Gamma(j+\frac{1}{2}+\sigma)}{j! \Gamma(n+\frac{1}{2}+\sigma)}\\
\hphantom{\frac{1}{\sqrt{n}}
			\sqrt{\frac{(n-1)! \Gamma(j+\frac{1}{2}+\sigma)}{j! \Gamma(n+\frac{1}{2}+\sigma)}} \Pi_{k,2j}}{}
=\Pi_{k,\ell} \frac{(n-1)!
					\Gamma(j+\frac{1}{2}+\sigma)}{j! \Gamma(n+\frac{1}{2}+\sigma)} .
\end{gather*}
	Thus, by~\eqref{Sigma_k,ell, k = 2m ge ell = 2n}, $d_{k,\ell}=(-1)^{m+n}\Pi_{k,\ell}\Sigma_{k,\ell}$.}
		
	Finally, take $k=2m+1\ge\ell=2n+1$ such that the equality of the statement holds for all $d_{k-1,2j}$ with $j\le n$. Then, by Lemma~\ref{l: d_k,ell, k=2m+1, ell = 2n+1},
		\begin{gather*}
			d_{k,\ell}=\sqrt{\frac{n+\frac{1}{2}+\sigma}{m+\frac{1}{2}+\sigma}} (-1)^{m+n}
			\Pi_{k-1,\ell-1}\Sigma_{k-1,\ell-1}\\
\hphantom{d_{k,\ell}=}{}-\frac{u}{\sqrt{m+\frac{1}{2}+\sigma}} \sum_{j=0}^n(-1)^{n-j}
			\sqrt{\frac{n! \Gamma(j+\frac{1}{2}+\sigma)}{j! \Gamma(n+\frac{3}{2}+\sigma)}}
			(-1)^{m+j}\Pi_{k-1,2j}\Sigma_{k-1,2j} .
		\end{gather*}
	Here, by~\eqref{Pi_k,ell, case k = 2m ge ell = 2n} and~\eqref{Pi_k,ell, case k = 2m+1 ge ell = 2n+1},
		\begin{gather*}
			\sqrt{\frac{n+\frac{1}{2}+\sigma}{m+\frac{1}{2}+\sigma}} \Pi_{k-1,\ell-1}=\Pi_{k,\ell} ,
		\end{gather*}
	and
		\begin{gather*}
			\frac{1}{\sqrt{m+\frac{1}{2}+\sigma}}
			\sqrt{\frac{n! \Gamma(j+\frac{1}{2}+\sigma)}{j! \Gamma(n+\frac{3}{2}+\sigma)}} \Pi_{k-1,2j}\\
			\qquad{} =\frac{1}{\sqrt{m+\frac{1}{2}+\sigma}}
			\sqrt{\frac{m! \Gamma(n+\frac{3}{2}+\sigma)}{n! \Gamma(m+\frac{1}{2}+\sigma)}}
			\frac{n! \Gamma(j+\frac{1}{2}+\sigma)}{j! \Gamma(n+\frac{3}{2}+\sigma)}
			=\Pi_{k,\ell} \frac{n!
			\Gamma(j+\frac{1}{2}+\sigma)}{j! \Gamma(n+\frac{3}{2}+\sigma)} .
		\end{gather*}
	Thus, by~\eqref{Sigma_k,ell, k = 2m+1 ge ell = 2n+1}, $d_{k,\ell}=(-1)^{m+n}\Pi_{k,\ell}\Sigma_{k,\ell}$.	
\end{proof}

\begin{Lemma}\label{l: Sigma_k,ell > 0}
	$\Sigma_{k,\ell}>0$ for all $k$ and $\ell$.
\end{Lemma}

\begin{proof}
	We proceed by induction on $\ell$. For $\ell\in\{0,1\}$, this is true by~\eqref{Sigma_k,0} and~\eqref{Sigma_k,1} because $\sigma>u-\frac{1}{2}$. If $\ell>1$ and the results holds for $\Sigma_{k',\ell'}$ with $\ell'<\ell$, then $\Sigma_{k,\ell}>0$ by~\eqref{Sigma_k,ell, k = 2m ge ell = 2n} and~\eqref{Sigma_k,ell = Sigma_k-2,ell-2 + ..., k = 2m+1 ge ell = 2n+1 > 1} since $\sigma>u-\frac{1}{2}$.
\end{proof}

\begin{Lemma}\label{l: Sigma_k,ell le (1 - frac 1-u m) Sigma_k-2,ell}
	If $k=2m>\ell=2n$ or $k=2m+1>\ell=2n+1$, then
		\begin{gather*}
			\Sigma_{k,\ell}\le\left(1-\frac{1-u}{m}\right)\Sigma_{k-2,\ell}.
		\end{gather*}
\end{Lemma}

\begin{proof}
	We proceed by induction on $\ell$. This is true for $\ell\in\{0,1\}$ by~\eqref{Sigma_k,0} and~\eqref{Sigma_k,1}.
	
	Now, suppose that the result is satisf\/ied by $\Sigma_{k',\ell'}$ with $\ell'<\ell$. If $k=2m>\ell=2n>0$, then, by~\eqref{Sigma_k,ell, k = 2m ge ell = 2n} and Lemma~\ref{l: Sigma_k,ell > 0},
		\begin{gather*}
		\Sigma_{k,\ell} \le\left(1-\frac{1-u}{m-1}\right)\Sigma_{k-3,\ell-1}
		 +u\sum_{j=0}^{n-1}
		\frac{(n-1)! \Gamma(j+\frac{1}{2}+\sigma)}{j! \Gamma(n+\frac{1}{2}+\sigma)}
		\left(1-\frac{1-u}{m}\right)\Sigma_{k-2,2j}\\
\hphantom{\Sigma_{k,\ell}}{}
\le\left(1-\frac{1-u}{m}\right)\Sigma_{k-2,\ell} .
		\end{gather*}
	If $k=2m+1>\ell=2n+1>1$, then, by~\eqref{Sigma_k,ell = Sigma_k-2,ell-2 + ..., k = 2m+1 ge ell = 2n+1 > 1} and Lemma~\ref{l: Sigma_k,ell > 0}, and since $\sigma>u-\frac{1}{2}$,
	\begin{gather*}
		\Sigma_{k,\ell}\le\left(1-\frac{u}{n+\frac{1}{2}+\sigma}\right)
		\left(1-\frac{1-u}{m-1}\right)\Sigma_{k-4,\ell-2}\\
		 \phantom{\Sigma_{k,\ell} =}{}+\left(1-\frac{u+n}{n+\frac{1}{2}+\sigma}\right)u\sum_{j=0}^{n-1}
		\frac{(n-1)!\Gamma(j+\frac{1}{2}+\sigma)}{j!\Gamma(n+\frac{1}{2}+\sigma)}
		\left(1-\frac{1-u}{m}\right)\Sigma_{k-3,2j}\\
	\phantom{\Sigma_{k,\ell}}{} <\left(1-\frac{1-u}{m}\right)\Sigma_{k-2,\ell} .\tag*{\qed}
	\end{gather*}
\renewcommand{\qed}{}
\end{proof}

\begin{Corollary}\label{c: Sigma_k,ell < (1 - frac u(1-u) m)Sigma_k-2,ell-2}
	If $k=2m\ge\ell=2n>0$, then
		\begin{gather*}
			\Sigma_{k-1,\ell-1}<\Sigma_{k,\ell}\le\left(1-\frac{u(1-u)}{m}\right)\Sigma_{k-2,\ell-2} .
		\end{gather*}
\end{Corollary}

\begin{proof}
	The f\/irst inequality is a direct consequence of~\eqref{Sigma_k,ell, k = 2m ge ell = 2n}, and Lemma~\ref{l: Sigma_k,ell > 0}. On the other hand, by~\eqref{Sigma_k,ell = Sigma_k-2,ell-2 - ..., k = 2m ge ell = 2n > 0}, and Lemmas~\ref{l: Sigma_k,ell > 0} and~\ref{l: Sigma_k,ell le (1 - frac 1-u m) Sigma_k-2,ell},
		\begin{gather*}
			\Sigma_{k,\ell}
			 \le\Sigma_{k-2,\ell-2}
			-\frac{u(1-u)}{m}\sum_{j=0}^{n-1}
			\frac{(n-1)!\Gamma(j+\frac{1}{2}+\sigma)}{j!\Gamma(n+\frac{1}{2}+\sigma)}\Sigma_{k-2,2j}\\
\hphantom{\Sigma_{k,\ell}}{}
			 =\left(1-\frac{u(1-u)}{m}\right)\Sigma_{k-2,\ell-2}
			-\frac{u(1-u)}{m}\sum_{j=0}^{n-2}
			\frac{(n-1)!\Gamma(j+\frac{1}{2}+\sigma)}{j!\Gamma(n+\frac{1}{2}+\sigma)}\Sigma_{k-2,2j}\\
\hphantom{\Sigma_{k,\ell}}{}
 \le\left(1-\frac{u(1-u)}{m}\right)\Sigma_{k-2,\ell-2} .\tag*{\qed}
		\end{gather*}
\renewcommand{\qed}{}	
\end{proof}

\begin{Corollary}\label{c: Sigma_k,ell < (1 - frac u n + frac 1 2 + sigma) Sigma_k-1,ell-1}
	If $k=2m+1\ge\ell=2n+1$, then
		\begin{gather*}
			\left(1-\frac{u}{n+\frac{1}{2}+\sigma}\right)\Sigma_{k-2,\ell-2}<\Sigma_{k,\ell}
			\le\left(1-\frac{u}{n+\frac{1}{2}+\sigma}\right)\Sigma_{k-1,\ell-1} .
		\end{gather*}
\end{Corollary}

\begin{proof}
	This follows from~\eqref{Sigma_k,ell, k = 2m+1, ell = 2n+1, 2nd},~\eqref{Sigma_k,ell = Sigma_k-2,ell-2 + ..., k = 2m+1 ge ell = 2n+1 > 1} and Lemma~\ref{l: Sigma_k,ell > 0} because $\sigma>u-\frac{1}{2}$.
\end{proof}

\begin{Lemma}\label{l: Weierstrass}
	For $0<t<1$, there is some $C_0=C_0(t)\ge1$ such that, for all $p$,
			\begin{gather*}
				C_0^{-1}(p+1)^{-t}\le\prod_{i=1}^p\left(1-\frac{t}{i}\right)\le C_0(p+1)^{-t} .
			\end{gather*}
\end{Lemma}

\begin{proof}
	 For each $t>0$, by the Weierstrass def\/inition of the gamma function,
		\begin{gather*}
			\Gamma(t)=\frac{e^{-\gamma t}}{t}
			\prod_{i=1}^\infty\left(1+\frac{t}{i}\right)^{-1}e^{t/i} ,
		\end{gather*}
	where $\gamma=\lim\limits_{j\to\infty}\big(\sum\limits_{i=1}^j\frac{1}{i}-\ln j\big)$ (the Euler--Mascheroni constant), there is some $K_0\ge1$ such that, for all $p\in{\mathbb{Z}}_+$,
			\begin{gather}\label{Weierstrass}
				K_0^{-1}\prod_{i=1}^pe^{-t/i}\le\prod_{i=1}^p\left(1+\frac{t}{i}\right)^{-1}
				\le K_0\prod_{i=1}^pe^{-t/i} .
			\end{gather}
	
	Now, assume that $0<t<1$, and observe that
		\begin{gather*}
			\prod_{i=1}^p\left(1-\frac{t}{i}\right)=\prod_{i=1}^p\left(1+\frac{t}{i-t}\right)^{-1} .
		\end{gather*}
	By the second inequality of~\eqref{Weierstrass}, for $p\ge1$,
		\begin{gather*}
			\prod_{i=1}^p\left(1+\frac{t}{i-t}\right)^{-1}
			<\prod_{i=1}^p\left(1+\frac{t}{i}\right)^{-1}
			\le K_0\prod_{i=1}^pe^{-t/i}=K_0\exp\left({-}t\sum_{i=1}^p\frac{1}{i}\right)\\
\hphantom{\prod_{i=1}^p\left(1+\frac{t}{i-t}\right)^{-1}}{}
			\le K_0\exp\left({-}t\int_1^{p+1}\frac{dx}{x}\right)
			=K_0(p+1)^{-t} .
		\end{gather*}
	On the other hand, by the f\/irst inequality of~\eqref{Weierstrass}, for $p\ge2$,
		\begin{gather*}
			\prod_{i=1}^p\left(1+\frac{t}{i-t}\right)^{-1}
			\ge(1-t)\prod_{i=1}^{p-1}\left(1+\frac{t}{i}\right)^{-1}
			\ge(1-t)K_0^{-1}\prod_{i=1}^{p-1}e^{-t/i}\\
\hphantom{\prod_{i=1}^p\left(1+\frac{t}{i-t}\right)^{-1}}{}
			=(1-t)K_0^{-1}\exp\left({-}t\sum_{i=1}^{p-1}\frac{1}{i}\right)
			\ge(1-t)K_0^{-1}\exp\left({-}t\left(1+\int_1^{p-1}\frac{dx}{x}\right)\right)\\
\hphantom{\prod_{i=1}^p\left(1+\frac{t}{i-t}\right)^{-1}}{}
			=(1-t)K_0^{-1}e^{-t}(p-1)^{-t}
			>(1-t)K_0^{-1}e^{-t}(p+1)^{-t} .\tag*{\qed}
		\end{gather*}
\renewcommand{\qed}{}
\end{proof}

\begin{Lemma}\label{l: Sigma_k,ell le C' m^-u(1-u) (m-n+1)^-(1-u)^2}
	There is some $C'=C'(u)>0$ such that
		\begin{gather*}
			\Sigma_{k,\ell}\le C'(m+1)^{-u(1-u)}(m-n+1)^{-(1-u)^2}
		\end{gather*}
	for $k=2m\ge\ell=2n$ or $k=2m+1\ge\ell=2n+1$.
\end{Lemma}

\begin{proof}
	 Suppose f\/irst that $k=2m\ge\ell=2n$. By Lemma~\ref{l: Sigma_k,ell le (1 - frac 1-u m) Sigma_k-2,ell} and Corollary~\ref{c: Sigma_k,ell < (1 - frac u(1-u) m)Sigma_k-2,ell-2}, we get
	 	\begin{gather*}
			\Sigma_{k,\ell}
			 \le\prod_{i=m-n+1}^m\left(1-\frac{u(1-u)}{i}\right)\prod_{i=1}^{m-n}\left(1-\frac{1-u}{i}\right)\\
\hphantom{\Sigma_{k,\ell}}{}
			 =\prod_{i=1}^m\left(1-\frac{u(1-u)}{i}\right)
			\prod_{i=1}^{m-n}\left(1-\frac{u(1-u)}{i}\right)^{-1}
			\prod_{i=1}^{m-n}\left(1-\frac{1-u}{i}\right).
		\end{gather*}
Then the result follows in this case from Lemma~\ref{l: Weierstrass}.
	
	When $k=2m+1\ge\ell=2n+1$, the result follows from the above case and Corollary~\ref{c: Sigma_k,ell < (1 - frac u n + frac 1 2 + sigma) Sigma_k-1,ell-1}.
\end{proof}

\begin{Lemma}\label{l: Gautschi}
	For each $t\in{\mathbb{R}}\setminus(-{\mathbb{N}})$, there is some $C_1=C_1(t)\ge1$ such that, for all $p$,
		\begin{gather*}
			C_1^{-1}(p+1)^{1-t}\le\frac{\Gamma(p+1)}{|\Gamma(p+t)|}\le C_1(p+1)^{1-t} .
		\end{gather*}
\end{Lemma}

\begin{proof}
	We can assume that $p\ge1$. Write $t=q+r$, where $q=\lfloor t\rfloor$. If $q=0$, then $0<r<1$ and the result follows from the Gautschi's inequality, stating that
		\begin{gather}\label{Gautschi}
			x^{1-r}\le\frac{\Gamma(x+1)}{\Gamma(x+r)}\le(x+1)^{1-r}
		\end{gather}
	for $0<r<1$ and $x>0$, because $x^{1-r}\ge2^{r-1}(x+1)^{1-r}$ for $x\ge1$.
	
	If $q\ge1$ and $r=0$, then
		\begin{gather*}
			\frac{\Gamma(p+1)}{\Gamma(p+t)}=\frac{p!}{(p+q-1)!}\le\frac{1}{(p+1)^{q-1}}=(p+1)^{1-t} ,\\
			\frac{\Gamma(p+1)}{\Gamma(p+t)}=\frac{p!}{(p+q-1)!}\ge\frac{1}{(p+q-1)^{q-1}}
			\ge\frac{1}{(qp)^{q-1}}\ge t^{1-t}(p+1)^{1-t} .
		\end{gather*}
		
	If $q\ge1$ and $r>0$, then, by~\eqref{Gautschi},
		\begin{gather*}
			\frac{\Gamma(p+1)}{\Gamma(p+t)}
			\le\frac{\Gamma(p+1)}{(p+1)^{q-1}(p+r)\Gamma(p+r)}
			\le\frac{(p+1)^{2-q-r}}{p+r}\le2(p+1)^{1-t},\\
			\frac{\Gamma(p+1)}{\Gamma(p+t)}\ge\frac{\Gamma(p+1)}{(p+t-1)^q\Gamma(p+r)}
			\ge\frac{p^{1-r}}{(p+t-1)^q}\ge\frac{(p+1)^{1-r}}{2^{1-r}(p+t-1)^q}\\
\hphantom{\frac{\Gamma(p+1)}{\Gamma(p+t)}}{}
\ge\min\{1,(t-1)^{-q}\}2^{r-1}(p+1)^{1-t},
		\end{gather*}
	because
		\begin{gather*}
			(p+t-1)^{-q}\ge
				\begin{cases}
					(p+1)^{-q} & \text{if $0<t\le2$}, \\
					(t-1)^{-q}(p+1)^{-q} & \text{if $t>2$}.
				\end{cases}
		\end{gather*}
	
	In the case $q<0$ ($t<0$), apply reverse induction on $q$: with $C_1=C_1(t+1)$, we get
		\begin{gather*}
			\frac{\Gamma(p+1)}{|\Gamma(p+t)|}=\frac{|p+t|\Gamma(p+1)}{|\Gamma(p+t+1)|}
			\le|p+t|C_1(p+1)^{-t}\le C_1|q|(p+1)^{1-t},\\
			\frac{\Gamma(p+1)}{|\Gamma(p+t)|}=\frac{|p+t|\Gamma(p+1)}{|\Gamma(p+t+1)|}
			\ge|p+t|C_1^{-1}(p+1)^{-t}=\frac{|p+t|}{p+1}C_1^{-1}(p+1)^{1-t},
		\end{gather*}
	where $|p+t|/(p+1)$ is bounded uniformly on $p$.
\end{proof}

\begin{Corollary}\label{c: Pi_k,ell le C''(frac n+1 m+1)^sigma/2 - 1/4}
	There is some $C''=C''(\sigma)>0$ such that
		\begin{gather*}
			\Pi_{k,\ell}\le
				\begin{cases}
					\displaystyle C''\left(\frac{n+1}{m+1}\right)^{\frac{\sigma}{2}-\frac{1}{4}} & \text{if $k=2m\ge\ell=2n$}, \vspace{1mm}\\
					\displaystyle C''\left(\frac{n+1}{m+1}\right)^{\frac{\sigma}{2}+\frac{1}{4}} & \text{if $k=2m+1\ge\ell=2n+1$} .
				\end{cases}
		\end{gather*}
\end{Corollary}

\begin{proof}
	This follows from~\eqref{Pi_k,ell, case k = 2m ge ell = 2n},~\eqref{Pi_k,ell, case k = 2m+1 ge ell = 2n+1} and Lemma~\ref{l: Gautschi}.
\end{proof}

For the sake of simplicity, let us use the following notation. For real valued functions~$f$ and~$g$ of $(m,n)$, for $(m,n)$ in some subset of ${\mathbb{N}}\times{\mathbb{N}}$, write $f\preccurlyeq g$ if there is some $C>0$ such that $f(m,n)\le C g(m,n)$ for all $(m,n)$. The same notation is used for functions depending also on other variables, $s,\sigma,u,\dots$, taking $C$ independent of $m$, $n$ and $s$, but possibly depending on the rest of variables.

\begin{Lemma}\label{l: (m+1)^alpha (n+1)^beta (m-n+1)^gamma}
	For $\alpha,\beta,\gamma\in{\mathbb{R}}$, if $\alpha+\beta,\alpha+\gamma,\alpha+\beta+\gamma<0$, then there is some $\omega>0$ such that, for all naturals $m\ge n$,
		\begin{gather*}
			(m+1)^\alpha(n+1)^\beta(m-n+1)^\gamma\preccurlyeq(m+1)^{-\omega}(n+1)^{-\omega} .
		\end{gather*}
\end{Lemma}

\begin{proof}
	We consider the following cases:

1. If $\alpha,\beta,\gamma<0$, then
				\begin{gather*}
					(m+1)^\alpha(n+1)^\beta(m-n+1)^\gamma
					\le(m+1)^\alpha(n+1)^\beta .
				\end{gather*}
			
2. If $\beta\ge0$ and $\gamma<0$, then
				\begin{gather*}
					(m+1)^\alpha(n+1)^\beta(m-n+1)^\gamma \le(m+1)^{\alpha+\beta}
					 \le(m+1)^{\frac{\alpha+\beta}{2}}(n+1)^{\frac{\alpha+\beta}{2}} .
				\end{gather*}
				
3. If $\alpha\ge0$ and $m+1\le2(n+1)$, then $\beta,\gamma<0$ and
				\begin{gather*}
					(m+1)^\alpha(n+1)^\beta(m-n+1)^\gamma \le2^{-\beta}(m+1)^{\alpha+\beta}
					 \le2^{-\beta}(m+1)^{\frac{\alpha+\beta}{2}}(n+1)^{\frac{\alpha+\beta}{2}} .
				\end{gather*}
			
4. If $\alpha\ge0$ and $m+1>2(n+1)$, then $\beta,\gamma<0$ and $m-n+1>(m+1)/2$, and therefore
				\begin{gather*}
					(m+1)^\alpha(n+1)^\beta(m-n+1)^\gamma
					\le2^{-\gamma}(m+1)^{\alpha+\gamma}(n+1)^\beta .
				\end{gather*}
			
5. If $\beta<0$ and $\gamma\ge0$, then
				\begin{gather*}
					(m+1)^\alpha(n+1)^\beta(m-n+1)^\gamma
					\le(m+1)^{\alpha+\gamma}(n+1)^\beta .
				\end{gather*}
			
6. If $\beta\ge0$ and $\gamma\ge0$, then
				\begin{gather*}
					(m+1)^\alpha(n+1)^\beta(m-n+1)^\gamma
					 \le(m+1)^{\alpha+\beta+\gamma}
					 \le(m+1)^{\frac{\alpha+\beta+\gamma}{2}}(n+1)^{\frac{\alpha+\beta+\gamma}{2}} .\tag*{\qed}
				\end{gather*}
\renewcommand{\qed}{}
\end{proof}

\begin{Proposition}\label{p: |d_k,ell|}
	There is some $\omega=\omega(\sigma,u)>0$ such that
		\begin{gather*}
			|d_{k,\ell}|\preccurlyeq(m+1)^{-\omega}(n+1)^{-\omega}
		\end{gather*}
	 for $k=2m$ and $\ell=2n$, or for $k=2m+1$ and $\ell=2n+1$.
\end{Proposition}

\begin{proof}
	We can assume $k\ge\ell$ because $d_{k,\ell}=d_{\ell,k}$.
	
	If $k=2m+1\ge\ell=2n+1$, then, according to Proposition~\ref{p: d_k,ell}, Lemma~\ref{l: Sigma_k,ell le C' m^-u(1-u) (m-n+1)^-(1-u)^2} and Corollary~\ref{c: Pi_k,ell le C''(frac n+1 m+1)^sigma/2 - 1/4},
		\begin{gather*}
			|d_{k,\ell}|\preccurlyeq(m+1)^{-\frac{\sigma}{2}-\frac{1}{4}-u(1-u)}(n+1)^{\frac{\sigma}{2}+\frac{1}{4}}(m-n+1)^{-(1-u)^2} .
		\end{gather*}
	Thus the result follows by Lemma~\ref{l: (m+1)^alpha (n+1)^beta (m-n+1)^gamma} since
		\begin{gather*}
			-\tfrac{\sigma}{2}-\tfrac{1}{4}-u(1-u)-(1-u)^2=-\tfrac{\sigma}{2}+u-\tfrac{5}{4}<\tfrac{u}{2}-1<0 .
		\end{gather*}
			
If $k=2m\ge\ell=2n$, then, according to Proposition~\ref{p: d_k,ell}, Lemma~\ref{l: Sigma_k,ell le C' m^-u(1-u) (m-n+1)^-(1-u)^2} and Corollary~\ref{c: Pi_k,ell le C''(frac n+1 m+1)^sigma/2 - 1/4},
		\begin{gather*}
			|d_{k,\ell}|\preccurlyeq(m+1)^{-\frac{\sigma}{2}+\frac{1}{4}-u(1-u)}(n+1)^{\frac{\sigma}{2}-\frac{1}{4}}(m-n+1)^{-(1-u)^2} .
		\end{gather*}
	Thus the result follows by Lemma~\ref{l: (m+1)^alpha (n+1)^beta (m-n+1)^gamma} since
		\begin{gather*}
			-\tfrac{\sigma}{2}+\tfrac{1}{4}-u(1-u)-(1-u)^2=-\tfrac{\sigma}{2}+u-\tfrac{3}{4}<\tfrac{u}{2}-\tfrac{1}{2}<0 .\tag*{\qed}
		\end{gather*}
\renewcommand{\qed}{}
\end{proof}

\begin{Corollary}\label{c: |c_k,ell|}
	There is some $\omega=\omega(\sigma,u)>0$ such that, for $k=2m$ and $\ell=2n$, or for $k=2m+1$ and $\ell=2n+1$,
		\begin{gather*}
			|c_{k,\ell}|\preccurlyeq s^u(m+1)^{-\omega}(n+1)^{-\omega} .
		\end{gather*}
\end{Corollary}

\begin{proof}
	This follows from Proposition~\ref{p: |d_k,ell|} and~\eqref{c_0,0}.
\end{proof}

\begin{Proposition}\label{p: ft le epsilon s^u-1 fl(phi)+C s^u |phi|_sigma^2}
	For any $\epsilon>0$, there is some $C=C(\epsilon,\sigma,u)>0$ such that, for all $\phi\in{\mathcal{S}}$,
		\begin{gather*}
			{\mathfrak{t}}(\phi)\le\epsilon s^{u-1} {\mathfrak{j}}(\phi)+Cs^u \|\phi\|_\sigma^2 .
		\end{gather*}
\end{Proposition}

\begin{proof}
	For each $k$, let $\nu_k=2k+1+2\sigma$. By Corollary~\ref{c: |c_k,ell|}, there are $K_0=K_0(\sigma,u)>0$ and $\omega=\omega(\sigma,u)>0$ such that
		\begin{gather}\label{|c_k,ell|}
			|c_{k,\ell}|\le K_0s^u\nu_k^{-\omega}\nu_\ell^{-\omega}
		\end{gather}
	for all $k$ and $\ell$. Since $S=S(\sigma,u):=\sum_k\nu_k^{-1-2\omega}<\infty$, given $\epsilon>0$, there is some $k_0=k_0(\epsilon,\sigma,u)$ so that
		\begin{gather*}
			S_0=S_0(\epsilon,\sigma,u):=\sum_{k>k_0}\nu_k^{-1-2\omega}<\frac{\epsilon^2}{4K_0^2S} .
		\end{gather*}
	Let $S_1=S_1(\epsilon,\sigma,u)=\sum\limits_{k\le k_0}\nu_k^{-\omega}$. For $\phi=\sum_kt_k\phi_k\in{\mathcal{S}}$, by~\eqref{|c_k,ell|} and the Schwartz inequality, we have
		\begin{gather*}
			{\mathfrak{t}}(\phi) =\sum_{k,\ell}t_k\overline{t_\ell} c_{k,\ell}\le\sum_{k,\ell}|t_k||t_\ell||c_{k,\ell}|\\
		\hphantom{{\mathfrak{t}}(\phi)}{}	 \le K_0s^{u-\frac{1}{2}}\sum_{k\le k_0}\frac{|t_k|}{\nu_k^\omega}
			\sum_\ell\frac{|t_\ell|(\nu_\ell s)^{\frac{1}{2}}}{\nu_\ell^{\frac{1}{2}+\omega}}
			 +K_0s^{u-1}
			\sum_{k>k_0}\frac{|t_k|(\nu_ks)^{\frac{1}{2}}}{\nu_k^{\frac{1}{2}+\omega}}
			\sum_\ell\frac{|t_\ell|(\nu_\ell s)^{\frac{1}{2}}}{\nu_\ell^{\frac{1}{2}+\omega}}\\
		\hphantom{{\mathfrak{t}}(\phi)}{}
			 \le K_0S_1S^{\frac{1}{2}}s^{u-\frac{1}{2}} \|\phi\|_\sigma {\mathfrak{j}}(\phi)^{\frac{1}{2}}
			+K_0S_0^{\frac{1}{2}}S^{\frac{1}{2}}s^{u-1} {\mathfrak{j}}(\phi)\\
		\hphantom{{\mathfrak{t}}(\phi)}{}
			\le K_0S_1S^{\frac{1}{2}}s^{u-\frac{1}{2}}\|\phi\|_\sigma{\mathfrak{j}}(\phi)^{\frac{1}{2}}+\frac{\epsilon s^{u-1}}{2}{\mathfrak{j}}(\phi)
			\le\frac{K_0^2S_1^2Ss^u}{2\epsilon}\|\phi\|_\sigma^2+\epsilon s^{u-1}{\mathfrak{j}}(\phi).\tag*{\qed}
		\end{gather*}
\renewcommand{\qed}{}
\end{proof}

\begin{Proposition}\label{p: ft(phi) ge D s^u/2 (k+1)^-u |phi|_sigma^2}
	There is some $D=D(\sigma,u)>0$ such that, for all $k\in{\mathbb{N}}$ and $\phi$ in the linear span of $\phi_0,\dots,\phi_k$,
		\begin{gather*}
			{\mathfrak{t}}(\phi)\ge Ds^u(k+1)^{-u}\|\phi\|_\sigma^2.
		\end{gather*}
\end{Proposition}

\begin{proof}
	Let $\phi=\sum\limits_{i=0}^kt_i\phi_i$ ($t_i\in{\mathbb{C}}$) and $\nu=\nu_k=2k+1+2\sigma$. Let $K\ge3$, which will be f\/ixed later. By~\eqref{|phi_k(x)| |x|^sigma le ...},
		\begin{gather*}
			\int_{|x|\ge\sqrt{\frac{K\nu}{2s}}}|\phi(x)|^2|x|^{2\sigma} dx
			=\sum_{i,j=0}^kt_i\overline{t_j}\int_{|x|\ge\sqrt{\frac{K\nu}{2s}}}\phi_i(x)\overline{\phi_j(x)} |x|^{2\sigma} dx\\
			\hphantom{\int_{|x|\ge\sqrt{\frac{K\nu}{2s}}}|\phi(x)|^2|x|^{2\sigma} dx}{}
			\le2\sum_{i,j=0}^k|t_i||t_j|\int_{\sqrt{\frac{K\nu}{2s}}}^\infty|\phi_i(x)||\phi_j(x)|x^{2\sigma} dx\\
			\hphantom{\int_{|x|\ge\sqrt{\frac{K\nu}{2s}}}|\phi(x)|^2|x|^{2\sigma} dx}{}
			\le\sum_{i,j=0}^k\big(|t_i|^2+|t_j|^2\big)C^2s\int_{\sqrt{\frac{K\nu}{2s}}}^\infty xe^{-2csx^2} dx\\
			\hphantom{\int_{|x|\ge\sqrt{\frac{K\nu}{2s}}}|\phi(x)|^2|x|^{2\sigma} dx}{}
			=2(k+1)\|\phi\|_\sigma^2C^2s\int_{\sqrt{\frac{K\nu}{2s}}}^\infty xe^{-2csx^2} dx\\
			\hphantom{\int_{|x|\ge\sqrt{\frac{K\nu}{2s}}}|\phi(x)|^2|x|^{2\sigma} dx}{}
			=\frac{C^2(k+1)}{2c}e^{-Kc\nu}\|\phi\|_\sigma^2,
		\end{gather*}
	 where $C,c>0$ depend only on $\sigma$. We can choose $K=K(\sigma)\ge3$ and $D=D(\sigma,u)>0$ such that
		\begin{gather*}
			\left(\frac{2s}{K\nu}\right)^u\left(1-\frac{C^2(k+1)}{2c}e^{-Kc\nu}\right)\ge Ds^u(k+1)^{-u}
		\end{gather*}
	for all $s>0$ and $k\in{\mathbb{N}}$, obtaining
		\begin{gather*}
			{\mathfrak{t}}(\phi)\ge\int_{|x|\le\sqrt{\frac{K\nu}{2s}}}|\phi(x)|^2|x|^{2\sigma-2u} dx
			\ge\left(\frac{2s}{K\nu}\right)^u\int_{|x|\le\sqrt{\frac{K\nu}{2s}}}|\phi(x)|^2|x|^{2\sigma} dx\\
			\hphantom{{\mathfrak{t}}(\phi)}{}
			\ge\left(\frac{2s}{K\nu}\right)^u\|\phi\|_\sigma^2\left(1-\frac{C^2(k+1)}{2c}e^{-Kc\nu}\right)
			\ge Ds^u(k+1)^{-u}\|\phi\|_\sigma^2.\tag*{\qed}
		\end{gather*}
\renewcommand{\qed}{}
\end{proof}

\begin{Remark}\label{r: ft(phi_k) ge D s^u/2 (k+1)^-u}
	For $\phi=\phi_k$, we can also use the following argument. By Proposition~\ref{p: d_k,ell} and~\eqref{c_0,0}, and since $\Pi_{k,k}=1$, it is enough to prove that there is some $D_0=D_0(\sigma,u)>0$ so that $\Sigma_{k,k}\ge D_0(k+1)^{-u}$. Moreover we can assume that $k=2m+1$ by Corollary~\ref{c: Sigma_k,ell < (1 - frac u(1-u) m)Sigma_k-2,ell-2}. We have $p_0:=\lfloor\frac{1}{2}+\sigma\rfloor\ge0$ because $\frac{1}{2}+\sigma>u$. According to Corollary~\ref{c: Sigma_k,ell < (1 - frac u n + frac 1 2 + sigma) Sigma_k-1,ell-1}, Lemma~\ref{l: Weierstrass} and~\eqref{Sigma_k,1}, there is some $C_0=C_0(u)\ge1$ such that
		\begin{gather*}
			\Sigma_{k,k}\ge\prod_{i=0}^m\left(1-\frac{u}{i+\frac{1}{2}+\sigma}\right)
			\ge\left(1-\frac{u}{\frac{1}{2}+\sigma}\right)\prod_{p=1+p_0}^{m+p_0}\left(1-\frac{u}{p}\right)\\
			\hphantom{\Sigma_{k,k}}{}
			=\left(1-\frac{u}{\frac{1}{2}+\sigma}\right)\prod_{p=1}^{m+p_0}\left(1-\frac{u}{p}\right)
			\prod_{p=1}^{p_0}\left(1-\frac{u}{p}\right)^{-1}\\
			\hphantom{\Sigma_{k,k}}{}
			\ge\left(1-\frac{u}{\frac{1}{2}+\sigma}\right)C_0^{-2}(m+p_0+1)^{-u}(p_0+1)^u
			\ge\left(1-\frac{u}{\frac{1}{2}+\sigma}\right)C_0^{-2}(k+1)^{-u}.
		\end{gather*}
\end{Remark}

\begin{Remark}
	If $0<u<\frac{1}{2}$, then $\lim\limits_m{\mathfrak{t}}(\phi_{2m+1})=0$. To check it, we use that there is some $K=K(\sigma,s)>0$ so that $|x|^{2\sigma}\phi_k^2(x)\le Kk^{-\frac{1}{6}}$ for all $x\in{\mathbb{R}}$ and all odd $k\in{\mathbb{N}}$ \cite[Theorem~1.1(ii)]{AlvCalaza:generalizedHermite} (this also follows from~\eqref{|phi_k(x)| |x|^sigma le ...}). For any $\epsilon>0$, take some $x_0>0$ and $k_0\in{\mathbb{N}}$ such that $x_0^{-2u}<\epsilon/2$ and $Kk_0^{-\frac{1}{6}}x_0^{1-2u}<\epsilon(1-2u)/4$. Then, for all odd natural $k\ge k_0$,
		\begin{gather*}
			{\mathfrak{t}}(\phi_k)
			=2\int_0^{x_0}\phi_k^2(x) x^{2(\sigma-u)} dx+2\int_{x_0}^\infty\phi_k^2(x) x^{2(\sigma-u)} dx\\
\hphantom{{\mathfrak{t}}(\phi_k)}{}
			\le2Kk^{-\frac{1}{6}}\int_0^{x_0}x^{-2u} dx+2x_0^{-2u}\int_{x_0}^\infty\phi_k^2(x) x^{2\sigma} dx
			\le2Kk^{-\frac{1}{6}}\frac{x_0^{1-2u}}{1-2u}+x_0^{-2u}<\epsilon,
		\end{gather*}
	because $1-2u>0$ and $\|\phi_k\|_\sigma=1$. In the case where $\sigma\ge0$, this argument is also valid when $k$ is even. We do not know if $\inf_k{\mathfrak{t}}(\phi_k)>0$ when $\frac{1}{2}\le u<1$.
\end{Remark}

\begin{proof}[Proof of Theorem~\ref{t: UU}]
	The positive def\/inite sesquilinear form ${\mathfrak{j}}$ of Section~\ref{s: prelim} is closable by \cite[Theorems~VI-2.1 and~VI-2.7]{Kato1995}. Then, taking $\epsilon>0$ so that $\xi\epsilon s^{u-1}<1$, it follows from \cite[Theorem~VI-1.33]{Kato1995} and Proposition~\ref{p: ft le epsilon s^u-1 fl(phi)+C s^u |phi|_sigma^2} that the positive def\/inite sesquilinear form ${\mathfrak{u}}:={\mathfrak{j}}+\xi{\mathfrak{t}}$ is also closable, and ${\mathsf{D}}(\bar{\mathfrak{u}})={\mathsf{D}}(\overline{{\mathfrak{j}}})$. By \cite[Theorems~VI-2.1,~VI-2.6 and~VI-2.7]{Kato1995}, there is a unique positive def\/inite self-adjoint operator ${\mathcal{U}}$ such that ${\mathsf{D}}({\mathcal{U}})$ is a core of ${\mathsf{D}}(\bar{\mathfrak{u}})$, which consists of the elements $\phi\in{\mathsf{D}}(\bar{\mathfrak{u}})$ so that, for some $\chi\in L^2_\sigma$, we have $\bar{\mathfrak{u}}(\phi,\psi)=\langle\chi,\psi\rangle_\sigma$ for all $\psi$ in some core of~$\bar{\mathfrak{u}}$ (in this case, ${\mathcal{U}}(\phi)=\chi$). By~\cite[Theorem~VI-2.23]{Kato1995}, we have ${\mathsf{D}}({\mathcal{U}}^{1/2})={\mathsf{D}}(\bar{\mathfrak{u}})$, ${\mathcal{S}}$ is a core of~${\mathcal{U}}^{1/2}$ (since it is a core of ${\mathfrak{u}}$), and~\eqref{langle UU^1/2 phi, UU^1/2 psi rangle_sigma} is satisf\/ied. By Proposition~\ref{p: ft(phi) ge D s^u/2 (k+1)^-u |phi|_sigma^2}, there is some $D(\sigma,u)$ so that, for all $s>0$ and $k\in{\mathbb{N}}$, and every $\phi$ is in the linear span of $\phi_0,\dots,\phi_k$, we have ${\mathfrak{t}}(\phi)\ge Ds^u(k+1)^{-u}\|\phi\|_\sigma^2$. Moreover we can assume that the sequence $(2k+1+2\sigma)s+\xi Ds^u(k+1)^{-u}$ is strictly increasing after reducing $D$ if necessary. So
		\begin{gather*}
			{\mathfrak{u}}(\phi)\ge\big((2k+1+2\sigma)s+\xi Ds^u(k+1)^{-u}\big)\|\phi\|_\sigma^2
		\end{gather*}
	 if $\phi\in{\mathcal{S}}$ is orthogonal in $L^2_\sigma$ to the linear span of $\phi_0,\dots,\phi_{k-1}$ (assuming that this span is $0$ when $k=0$). Therefore ${\mathcal{U}}$ has a discrete spectrum satisfying the f\/irst inequality of~\eqref{lambda_k, case of UU} by the form version of the min-max principle \cite[Theorem~XIII.2]{ReedSimon1978}. The second inequality of~\eqref{lambda_k, case of UU} holds because
		\begin{gather*}
			\bar{\mathfrak{u}}(\phi)\le\big(1+\xi\epsilon s^{u-1}\big)\bar{\mathfrak{j}}(\phi)+\xi Cs^u\|\phi\|_\sigma^2
		\end{gather*}
	for all $\phi\in{\mathsf{D}}(\bar{\mathfrak{u}})$ by Proposition~\ref{p: ft le epsilon s^u-1 fl(phi)+C s^u |phi|_sigma^2} and \cite[Theorem~VI-1.18]{Kato1995}, since ${\mathcal{S}}$ is a core of $\bar{\mathfrak{u}}$ and $\bar{\mathfrak{j}}$.
\end{proof}

\begin{Remark}\label{r: bar fp = bar fl + xi bar ft}
	In the above proof, note that $\bar{\mathfrak{u}}=\bar{\mathfrak{j}}+\xi\bar{\mathfrak{t}}$ and ${\mathsf{D}}(\bar{\mathfrak{j}})={\mathsf{D}}\big(\overline J^{1/2}\big)$. Thus~\eqref{langle UU^1/2 phi, UU^1/2 psi rangle_sigma} can be extended to $\phi,\psi\in{\mathsf{D}}\big({\mathcal{U}}^{1/2}\big)$ using $\big\langle\overline J^{1/2}\phi,\overline J^{1/2}\psi\big\rangle_\sigma$ instead of $\langle J\phi,\psi\rangle_\sigma$.
\end{Remark}

\begin{Remark}\label{r: type (B)}
	 Extend the def\/inition of the above forms and operators to the case of $\xi\in{\mathbb{C}}$. Then $|\bar{\mathfrak{t}}(\phi)|\le\epsilon s^{u-1}\,\Re\bar{\mathfrak{j}}(\phi)+Cs^u \|\phi\|_\sigma^2$ for all $\phi\in{\mathsf{D}}(\bar{\mathfrak{j}})$, like in the proof of Theorem~\ref{t: UU}. Thus the family $\bar{\mathfrak{u}}=\bar{\mathfrak{u}}(\xi)$ becomes holomorphic of type~(a) by Remark~\ref{r: bar fp = bar fl + xi bar ft} and \cite[Theorem~VII-4.8]{Kato1995}, and therefore ${\mathcal{U}}={\mathcal{U}}(\xi)$ is a self-adjoint holomorphic family of type~(B). So the functions $\lambda_k=\lambda_k(\xi)$ ($\xi\in{\mathbb{R}}$) are continuous and piecewise holomorphic \cite[Remark~VII-4.22, Theorem~VII-3.9, and~VII-\S~3.4]{Kato1995}, with $\lambda_k(0)=(2k+1+2\sigma)s$. Moreover \cite[Theorem~VII-4.21]{Kato1995} gives an exponential estimate of $|\lambda_k(\xi)-\lambda_k(0)|$ in terms of $\xi$. But~\eqref{lambda_k, case of UU} is a~better estimate.
\end{Remark}

\section{Scalar products of mixed generalized Hermite functions}\label{s: mixed}

Let $\sigma,\tau,\theta>-\frac{1}{2}$, and write $v=\sigma+\tau-2\theta$. This section is devoted to describe the scalar products
	\begin{gather*}
		\hat c_{k,\ell}=\hat c_{\sigma,\tau,\theta,k,\ell}=\langle\phi_{\sigma,k},\phi_{\tau,\ell}\rangle_\theta ,
	\end{gather*}
which will be needed to prove Theorem~\ref{t: VV}. Note that $\hat c_{k,\ell}=0$ if $k+\ell$ is odd, and
	\begin{gather}\label{hat c_sigma,tau,xi,k,ell = hat c_tau,sigma,xi,ell,k}
		\hat c_{\sigma,\tau,\theta,k,\ell}=\hat c_{\tau,\sigma,\theta,\ell,k}
	\end{gather}
for all $k$ and $\ell$. Of course, $\hat c_{k,\ell}=\delta_{k,\ell}$ if $\sigma=\tau=\theta$.

According to Section~\ref{s: prelim}, if $k$ and $\ell$ are odd, then $\hat c_{\sigma,\tau,\theta,k,\ell}$ is also def\/ined when $\sigma,\tau,\theta>-\frac{3}{2}$, and we have
		\begin{gather}\label{hat c_sigma,tau,xi,k,ell = hat c_sigma+1,tau+1,xi+1,k-1,ell-1}
			\hat c_{\sigma,\tau,\theta,k,\ell}=\langle x\phi_{\sigma+1,k-1},x\phi_{\tau+1,\ell-1}\rangle_\theta
			=\hat c_{\sigma+1,\tau+1,\theta+1,k-1,\ell-1} .
		\end{gather}

\subsection[Case where $\sigma=\theta\ne\tau$ and $\tau-\sigma\not\in-{\mathbb{N}}$]{Case where $\boldsymbol{\sigma=\theta\ne\tau}$ and $\boldsymbol{\tau-\sigma\not\in-{\mathbb{N}}}$}\label{ss: mixed, case where sigma = theta ne tau}

In this case, we have $v=\tau-\sigma$. By~\eqref{phi_0} and~\eqref{int_-infty^infty e^-sx^2 |x|^2 kappa dx},
	\begin{gather}\label{hat c_0,0, case sigma = theta ne tau}
		\hat c_{0,0}=s^{\frac{v}{2}} \Gamma\big(\sigma+\tfrac{1}{2}\big)^{\frac{1}{2}} \Gamma\big(\tau+\tfrac{1}{2}\big)^{-\frac{1}{2}}.
	\end{gather}

\begin{Lemma}\label{l: hat c_k,0, case where sigma = theta ne tau}
	If $k>0$ is even, then $\hat c_{k,0}=0$.
\end{Lemma}

\begin{proof}
	By~\eqref{phi_k},~\eqref{B phi_0} and~\eqref{B_sigma},
		\begin{gather*}
			\hat c_{k,0}=\frac{1}{\sqrt{2ks}}\langle B'_\sigma\phi_{\sigma,k-1},\phi_{\tau,0}\rangle_\sigma
			=\frac{1}{\sqrt{2ks}}\langle\phi_{\sigma,k-1},B_\tau\phi_{\tau,0}\rangle_\sigma
			=0 .\tag*{\qed}
		\end{gather*}
\renewcommand{\qed}{}
\end{proof}

\begin{Lemma}\label{l: hat c_0,ell, 1st expression, case where sigma = theta ne tau}
	If $\ell=2n>0$, then
		\begin{gather*}
			\hat c_{0,\ell}=\frac{v}{\sqrt{n}}
			\sum_{j=0}^{n-1}(-1)^{n-j}\sqrt{\frac{(n-1)! \Gamma(j+\frac{1}{2}+\tau)}
			{j! \Gamma(n+\frac{1}{2}+\tau)}} \hat c_{0,2j} .
		\end{gather*}
\end{Lemma}

\begin{proof}
	By~\eqref{phi_k},~\eqref{B phi_0},~\eqref{x^-1 p_k} and~\eqref{B'_sigma},
		\begin{gather*}
			\hat c_{0,\ell}
			=\frac{1}{\sqrt{2\ell s}}
			\langle\phi_{\sigma,0},B'_\tau\phi_{\tau,\ell-1}\rangle_\sigma
			=\frac{1}{\sqrt{2\ell s}}
			\langle\phi_{\sigma,0},(B'_\sigma-2vx^{-1})\phi_{\tau,\ell-1}\rangle_\sigma\\
\hphantom{\hat c_{0,\ell}}{}
			=\frac{1}{\sqrt{2\ell s}}
			\langle B_\sigma\phi_{\sigma,0},\phi_{\tau,\ell-1}\rangle_\sigma
			-\frac{2v}{\sqrt{2\ell}}
			\sum_{j=0}^{n-1}(-1)^{n-1-j}\sqrt{\frac{(n-1)! \Gamma(j+\frac{1}{2}+\tau)}
			{j! \Gamma(n+\frac{1}{2}+\tau)}} \hat c_{0,2j}\\
\hphantom{\hat c_{0,\ell}}{}
			=\frac{v}{\sqrt{n}}
			\sum_{j=0}^{n-1}(-1)^{n-j}\sqrt{\frac{(n-1)! \Gamma(j+\frac{1}{2}+\tau)}
			{j! \Gamma(n+\frac{1}{2}+\tau)}} \hat c_{0,2j} .\tag*{\qed}
		\end{gather*}
\renewcommand{\qed}{}
\end{proof}

\begin{Lemma}\label{l: hat c_k,ell, case k = 2m > 0, ell = 2n > 0 and sigma = theta ne tau}
	If $k=2m>0$ and $\ell=2n>0$, then $\hat c_{k,\ell}=\sqrt{n/m} \hat c_{k-1,\ell-1}$.
\end{Lemma}

\begin{proof}
	By~\eqref{phi_k},~\eqref{B phi_k} and~\eqref{B_sigma},
		\begin{gather*}
			\hat c_{k,\ell}
			=\frac{1}{\sqrt{2ks}}
			\langle B'_\sigma\phi_{\sigma,k-1},\phi_{\tau,\ell}\rangle_\sigma
			=\frac{1}{\sqrt{2ks}}
			\langle\phi_{\sigma,k-1},B_\tau\phi_{\tau,\ell}\rangle_\sigma
			=\sqrt{\frac{n}{m}}\hat c_{k-1,\ell-1}.\tag*{\qed}
		\end{gather*}
\renewcommand{\qed}{}
\end{proof}

\begin{Lemma}\label{l: hat d_k,ell, case k = 2m+1, ell = 2n+1, and sigma = theta ne tau}
	If $k=2m+1$ and $\ell=2n+1$, then
		\begin{gather*}
			\hat c_{k,\ell}=\frac{n+\frac{1}{2}+\sigma}
			{\sqrt{(m+\frac{1}{2}+\sigma)(n+\frac{1}{2}+\tau)}} \hat c_{k-1,\ell-1} \\
\hphantom{\hat c_{k,\ell}=}{}
 -\frac{v}{\sqrt{m+\frac{1}{2}+\sigma}}
			\sum_{j=0}^{n-1}(-1)^{n-j}\sqrt{\frac{n! \Gamma(j+\frac{1}{2}+\tau)}
			{j! \Gamma(n+\frac{3}{2}+\tau)}} \hat c_{k-1,2j} .
		\end{gather*}
\end{Lemma}

\begin{proof}
	By~\eqref{phi_k},~\eqref{B phi_k},~\eqref{x^-1 p_k} and~\eqref{B_sigma},
		\begin{gather*}
			\hat c_{k,\ell}
			 =\frac{1}{\sqrt{2(k+2\sigma)s}}
			\langle B'_\sigma\phi_{\sigma,k-1},\phi_{\tau,\ell}\rangle_\sigma
			 =\frac{1}{\sqrt{2(k+2\sigma)s}}
			\big\langle\phi_{\sigma,k-1},\big(B_\tau-2vx^{-1}\big)\phi_{\tau,\ell}\big\rangle_\sigma\\
\hphantom{\hat c_{k,\ell}}{}
			 =\sqrt{\frac{n+\frac{1}{2}+\tau}{m+\frac{1}{2}+\sigma}} \hat c_{k-1,\ell-1}
			-\frac{v}{\sqrt{m+\frac{1}{2}+\sigma}}
			\sum_{j=0}^n(-1)^{n-j}\sqrt{\frac{n! \Gamma(j+\frac{1}{2}+\tau)}
			{j! \Gamma(n+\frac{3}{2}+\tau)}} \hat c_{k-1,2j}\\
\hphantom{\hat c_{k,\ell}}{}
			=\frac{n+\frac{1}{2}+\sigma}
			{\sqrt{(m+\frac{1}{2}+\sigma)(n+\frac{1}{2}+\tau)}}\hat c_{k-1,\ell-1}\\
\hphantom{\hat c_{k,\ell}=}{}
-\frac{v}{\sqrt{m+\frac{1}{2}+\sigma}}
			\sum_{j=0}^{n-1}(-1)^{n-j}\sqrt{\frac{n!\Gamma(j+\frac{1}{2}+\tau)}
			{j!\Gamma(n+\frac{3}{2}+\tau)}}\hat c_{k-1,2j}.\tag*{\qed}
		\end{gather*}
\renewcommand{\qed}{}
\end{proof}

\begin{Corollary}\label{c: hat c_k,ell = 0 if k >ell, case where sigma = theta ne tau}
	If $k>\ell$, then $\hat c_{k,\ell}=0$.
\end{Corollary}

\begin{proof}
	This follows by induction on $\ell$ using Lemmas~\ref{l: hat c_k,0, case where sigma = theta ne tau},~\ref{l: hat c_k,ell, case k = 2m > 0, ell = 2n > 0 and sigma = theta ne tau} and~\ref{l: hat d_k,ell, case k = 2m+1, ell = 2n+1, and sigma = theta ne tau}.
\end{proof}

\begin{Remark}\label{r: hat c_k,ell, case k = 2m+1 and ell = 2n+1}
	By Corollary~\ref{c: hat c_k,ell = 0 if k >ell, case where sigma = theta ne tau}, in Lemma~\ref{l: hat d_k,ell, case k = 2m+1, ell = 2n+1, and sigma = theta ne tau}, it is enough to consider the sum with $j$ running from $m$ to $n-1$.
\end{Remark}

\begin{Proposition}\label{p: hat c_k,ell, case where sigma = theta ne tau}
	If $k=2m\le\ell=2n$, then		
		\begin{gather*}
			\hat c_{k,\ell}=(-1)^{m+n}s^{\frac{v}{2}}
			\sqrt{\frac{n! \Gamma(m+\frac{1}{2}+\sigma)}{m! \Gamma(n+\frac{1}{2}+\tau)}}
			\frac{\Gamma(n-m+v)}{(n-m)! \Gamma(v)} ,
		\end{gather*}
	and, if $k=2m+1\le\ell=2n+1$, then
		\begin{gather*}
			\hat c_{k,\ell}=(-1)^{m+n}s^{\frac{v}{2}}
			\sqrt{\frac{n! \Gamma(m+\frac{3}{2}+\sigma)}{m! \Gamma(n+\frac{3}{2}+\tau)}}
			\frac{\Gamma(n-m+v)}{(n-m)! \Gamma(v)} .
		\end{gather*}
\end{Proposition}

\begin{proof}
	This is proved by induction on $k$. In turn, the case $k=0$,
		\begin{gather}\label{hat c_0,ell}
			\hat c_{0,\ell}=(-1)^ns^{\frac{v}{2}}
			\sqrt{\frac{\Gamma(\frac{1}{2}+\sigma)}{n! \Gamma(n+\frac{1}{2}+\tau)}}
			\frac{\Gamma(n+v)}{\Gamma(v)} ,
		\end{gather}
	is proved by induction on $\ell$. If $k=\ell=0$,~\eqref{hat c_0,ell} is~\eqref{hat c_0,0, case sigma = theta ne tau}. Given $\ell=2n>0$, assume that the result holds for $k=0$ and all $\ell'=2n'<\ell$. Then, by Lemma~\ref{l: hat c_0,ell, 1st expression, case where sigma = theta ne tau},
		\begin{gather*}
			\hat c_{0,\ell} =\frac{v}{\sqrt{n}}
			\sum_{j=0}^{n-1}(-1)^{n-j}\sqrt{\frac{(n-1)! \Gamma(j+\frac{1}{2}+\tau)}
			{j! \Gamma(n+\frac{1}{2}+\tau)}}
			(-1)^js^{\frac{v}{2}}\sqrt{\frac{\Gamma(\frac{1}{2}+\sigma)}{j! \Gamma(j+\frac{1}{2}+\tau)}}
			\frac{\Gamma(j+v)}{\Gamma(v)}\\
\hphantom{\hat c_{0,\ell}}{}
			 =(-1)^ns^{\frac{v}{2}}\sqrt{\frac{(n-1)! \Gamma(\frac{1}{2}+\sigma)}{n \Gamma(n+\frac{1}{2}+\tau)}}
			\frac{v}{\Gamma(v)}\sum_{j=0}^{n-1}\frac{\Gamma(j+v)}{j!} ,
		\end{gather*}
	obtaining~\eqref{hat c_0,ell} because
		\begin{gather}\label{v sum_i=0^p frac Gamma(i+v) i!}
			\frac{\Gamma(p+1+t)}{p!}=t\sum_{i=0}^p\frac{\Gamma(i+t)}{i!}
		\end{gather}
	for all $p\in{\mathbb{N}}$ and $t\in{\mathbb{R}}\setminus(-{\mathbb{N}})$, as can be easily checked by induction on~$p$.
	
	Given $k>0$, assume that the result holds for all $k'<k$. If $k$ is even, the statement follows directly from Lemma~\ref{l: hat c_k,ell, case k = 2m > 0, ell = 2n > 0 and sigma = theta ne tau}. If $k$ is odd, by Lemma~\ref{l: hat d_k,ell, case k = 2m+1, ell = 2n+1, and sigma = theta ne tau}, Remark~\ref{r: hat c_k,ell, case k = 2m+1 and ell = 2n+1} and~\eqref{v sum_i=0^p frac Gamma(i+v) i!},
		\begin{gather*}
			\hat c_{k,\ell} =\frac{n+\frac{1}{2}+\sigma}
			{\sqrt{(m+\frac{1}{2}+\sigma)(n+\frac{1}{2}+\tau)}}
			(-1)^{m+n}s^{\frac{v}{2}}
			\sqrt{\frac{n! \Gamma(m+\frac{1}{2}+\sigma)}{m! \Gamma(n+\frac{1}{2}+\tau)}}
			\frac{\Gamma(n-m+v)}{(n-m)! \Gamma(v)}\\
\hphantom{\hat c_{k,\ell} =}{}
-\frac{v}{\sqrt{m+\frac{1}{2}+\sigma}}
			\sum_{j=m}^{n-1}(-1)^{n-j}\sqrt{\frac{n! \Gamma(j+\frac{1}{2}+\tau)}
			{j! \Gamma(n+\frac{3}{2}+\tau)}}\\
\hphantom{\hat c_{k,\ell} =}{}
\times(-1)^{m+j}s^{\frac{v}{2}}
			\sqrt{\frac{j! \Gamma(m+\frac{1}{2}+\sigma)}{m! \Gamma(j+\frac{1}{2}+\tau)}}
			\frac{\Gamma(j-m+v)}{(j-m)! \Gamma(v)}\\
\hphantom{\hat c_{k,\ell}}{}
=(-1)^{m+n}s^{\frac{v}{2}}\sqrt{\frac{n!\Gamma(m+\frac{1}{2}+\sigma)}
			{(m+\frac{1}{2}+\sigma)m!\Gamma(n+\frac{3}{2}+\tau)}}\frac{1}{\Gamma(v)}\\
\hphantom{\hat c_{k,\ell} =}{}
\times\left(\frac{\Gamma(n-m+v) (n+\frac{1}{2}+\sigma)}{(n-m)!}
			-v\sum_{i=0}^{n-m-1}\frac{\Gamma(i+v)}{i!}\right)\\
\hphantom{\hat c_{k,\ell}}{}
			=(-1)^{m+n}s^{\frac{v}{2}}\sqrt{\frac{n!\Gamma(m+\frac{3}{2}+\sigma)}
			{m!\Gamma(n+\frac{3}{2}+\tau)}}\frac{\Gamma(n-m+v)}{(n-m)!\Gamma(v)}.\tag*{\qed}
		\end{gather*}
\renewcommand{\qed}{}
\end{proof}

\begin{Remark}\label{r: hat c_k,ell, case where sigma = theta ne tau}
	By~\eqref{hat c_sigma,tau,xi,k,ell = hat c_sigma+1,tau+1,xi+1,k-1,ell-1}, if $k$ and $\ell$ are odd, then Corollary~\ref{c: hat c_k,ell = 0 if k >ell, case where sigma = theta ne tau} and Proposition~\ref{p: hat c_k,ell, case where sigma = theta ne tau} also hold when $\sigma,\tau>-\frac{3}{2}$.
\end{Remark}

\subsection[Case where $\sigma\ne\theta\ne\tau$ and $\sigma-\theta,\tau-\theta\not\in-{\mathbb{N}}$]{Case where $\boldsymbol{\sigma\ne\theta\ne\tau}$ and $\boldsymbol{\sigma-\theta,\tau-\theta\not\in-{\mathbb{N}}}$}\label{ss: mixed, case where sigma ne theta ne tau}

By~\eqref{phi_0} and~\eqref{int_-infty^infty e^-sx^2 |x|^2 kappa dx},
	\begin{gather}\label{hat c_0,0, case where sigma ne theta ne tau}
		\hat c_{0,0}=s^{\frac{v}{2}}\Gamma\big(\sigma+\tfrac{1}{2}\big)^{-\frac{1}{2}}\Gamma\big(\tau+\tfrac{1}{2}\big)^{-\frac{1}{2}}\Gamma(\theta+\tfrac{1}{2}).
	\end{gather}

\begin{Lemma}\label{l: hat c_k,0, case where sigma ne theta ne tau}
	If $k=2m>0$, then
		\begin{gather*}
			\hat c_{k,0}=\frac{\sigma-\theta}{\sqrt{m}}\sum_{i=0}^{m-1}(-1)^{m-i}
			\sqrt{\frac{(m-1)!\Gamma(i+\frac{1}{2}+\sigma)}{i!\Gamma(m+\frac{1}{2}+\sigma)}}
			\hat c_{2i,0}.
		\end{gather*}
\end{Lemma}

\begin{proof}
	By~\eqref{phi_k} and~\eqref{B'_sigma},
		\begin{gather*}
			\hat c_{k,0}=\frac{1}{\sqrt{2ks}}\langle B'_\sigma\phi_{\sigma,k-1},\phi_{\tau,0}\rangle_\theta
			=\frac{1}{2\sqrt{ms}}\langle B'_\theta\phi_{\sigma,k-1},\phi_{\tau,0}\rangle_\theta
			+\frac{\theta-\sigma}{\sqrt{ms}}\big\langle x^{-1}\phi_{\sigma,k-1},\phi_{\tau,0}\big\rangle_\theta.
		\end{gather*}
	Here, by~\eqref{B phi_0},~\eqref{x^-1 p_k} and~\eqref{B_sigma},
		\begin{gather*}
			\langle B'_\theta\phi_{\sigma,k-1},\phi_{\tau,0}\rangle_\theta
			=\langle\phi_{\sigma,k-1},B_\theta\phi_{\tau,0}\rangle_\theta
			=\langle\phi_{\sigma,k-1},B_\tau\phi_{\tau,0}\rangle_\theta=0,\\
			\langle x^{-1}\phi_{\sigma,k-1},\phi_{\tau,0}\rangle_\theta
			=-\sum_{i=0}^{m-1}(-1)^{m-i}
			\sqrt{\frac{(m-1)!\Gamma(i+\frac{1}{2}+\sigma)s}{i!\Gamma(m+\frac{1}{2}+\sigma)}}\hat c_{2i,0}.\tag*{\qed}
		\end{gather*}
\renewcommand{\qed}{}
\end{proof}

\begin{Lemma}\label{l: hat c_k,ell, case k = 2m > 0 and ell = 2n > 0 and sigma ne theta ne tau}
	If $k=2m>0$ and $\ell=2n>0$, then
		\begin{gather*}
			\hat c_{k,\ell}=\sqrt{\frac{n}{m}}\hat c_{k-1,\ell-1}+\frac{\sigma-\theta}{m}
			\sum_{i=0}^{m-1}(-1)^{m-i}
			\sqrt{\frac{m!\Gamma(i+\frac{1}{2}+\sigma)}{i!\Gamma(m+\frac{1}{2}+\sigma)}}\hat c_{2i,\ell}.
		\end{gather*}
\end{Lemma}

\begin{proof}
	Like in the proof of Lemma~\ref{l: hat c_k,0, case where sigma ne theta ne tau},
		\begin{gather*}
			\hat c_{k,\ell}
			=\frac{1}{2\sqrt{ms}}\langle B'_\theta\phi_{\sigma,k-1},\phi_{\tau,\ell}\rangle_\theta
			+\frac{\theta-\sigma}{\sqrt{ms}}\big\langle x^{-1}\phi_{\sigma,k-1},\phi_{\tau,\ell}\big\rangle_\theta .
		\end{gather*}
	Now, by~\eqref{B phi_k},~\eqref{x^-1 p_k} and~\eqref{B_sigma},
		\begin{gather*}
			\langle B'_\theta\phi_{\sigma,k-1},\phi_{\tau,\ell}\rangle_\theta
			=\langle\phi_{\sigma,k-1},B_\theta\phi_{\tau,\ell}\rangle_\theta
			=\langle\phi_{\sigma,k-1},B_\tau\phi_{\tau,\ell}\rangle_\theta
			=2\sqrt{ns} \hat c_{k-1,\ell-1} ,\\
			\big\langle x^{-1}\phi_{\sigma,k-1},\phi_{\tau,\ell}\big\rangle_\theta
			=-\sum_{i=0}^{m-1}(-1)^{m-i}
			\sqrt{\frac{(m-1)! \Gamma(i+\frac{1}{2}+\sigma)s}{i! \Gamma(m+\frac{1}{2}+\sigma)}}
			\hat c_{2i,\ell} .\tag*{\qed}
		\end{gather*}
\renewcommand{\qed}{}
\end{proof}

\begin{Lemma}\label{l: hat c_k,ell, case k = 2m+1 and ell = 2n+1 and sigma ne theta ne tau}
	If $k=2m+1$ and $\ell=2n+1$, then
		\begin{gather*}
			\hat c_{k,\ell}=\frac{m+\frac{1}{2}+\theta}{\sqrt{(m+\frac{1}{2}+\sigma)(n+\frac{1}{2}+\tau)}}
			\hat c_{k-1,\ell-1}\\
\hphantom{\hat c_{k,\ell}=}{}
-\frac{\sigma-\theta}{\sqrt{n+\frac{1}{2}+\tau}}
			\sum_{i=0}^{m-1}(-1)^{m-i}\sqrt{\frac{m! \Gamma(i+\frac{1}{2}+\sigma)}
			{i! \Gamma(m+\frac{3}{2}+\sigma)}} \hat c_{2i,\ell-1} .
		\end{gather*}
\end{Lemma}

\begin{proof}
	By~\eqref{phi_k},
		\begin{gather*}
			\hat c_{k,\ell}
			=\frac{1}{2\sqrt{(n+\frac{1}{2}+\tau)s}}
			\langle\phi_{\sigma,k},B'_\tau\phi_{\tau,\ell-1}\rangle_\theta,
		\end{gather*}
	where, by~\eqref{B'_sigma},
		\begin{gather*}
			\langle\phi_{\sigma,k},B'_\tau\phi_{\tau,\ell-1}\rangle_\theta
			=\langle\phi_{\sigma,k},B'_\theta\phi_{\tau,\ell-1}\rangle_\theta
			=\langle B_\theta\phi_{\sigma,k},\phi_{\tau,\ell-1}\rangle_\theta\\
\hphantom{\langle\phi_{\sigma,k},B'_\tau\phi_{\tau,\ell-1}\rangle_\theta}{}
			=\langle B_\sigma\phi_{\sigma,k},\phi_{\tau,\ell-1}\rangle_\theta
			+2(\theta-\sigma)\langle x^{-1}\phi_{\sigma,k},\phi_{\tau,\ell-1}\rangle_\theta .
		\end{gather*}
	Hence, by~\eqref{B phi_k} and~\eqref{x^-1 p_k},
		\begin{gather*}
			\hat c_{k,\ell} =\sqrt{\frac{m+\frac{1}{2}+\sigma}{n+\frac{1}{2}+\tau}} \hat c_{k-1,\ell-1}
			-\frac{\sigma-\theta}{\sqrt{n+\frac{1}{2}+\tau}}
			\sum_{i=0}^m(-1)^{m-i}\sqrt{\frac{m! \Gamma(i+\frac{1}{2}+\sigma)}
			{i! \Gamma(m+\frac{3}{2}+\sigma)}} \hat c_{2i,\ell-1}\\
\hphantom{\hat c_{k,\ell}}{}
			 =\frac{m+\frac{1}{2}+\theta}{\sqrt{(n+\frac{1}{2}+\tau)(m+\frac{1}{2}+\sigma)}} \hat c_{k-1,\ell-1}\\
\hphantom{\hat c_{k,\ell}=}{}
-\frac{\sigma-\theta}{\sqrt{n+\frac{1}{2}+\tau}}
			\sum_{i=0}^{m-1}(-1)^{m-i}\sqrt{\frac{m! \Gamma(i+\frac{1}{2}+\sigma)}
			{i! \Gamma(m+\frac{3}{2}+\sigma)}} \hat c_{2i,\ell-1} .\tag*{\qed}
\end{gather*}
\renewcommand{\qed}{}
\end{proof}

\begin{Proposition}\label{p: hat c_k,ell, case where sigma ne theta ne tau}
	If $k=2m$ and $\ell=2n$, then		
		\begin{gather*}
			\hat c_{k,\ell}=(-1)^{m+n}s^{\frac{v}{2}}
			\sqrt{\frac{m! n!}{\Gamma(m+\frac{1}{2}+\sigma)\Gamma(n+\frac{1}{2}+\tau)}}\\
\hphantom{\hat c_{k,\ell}=}{} \times\sum_{p=0}^{\min\{m,n\}}
			\frac{\Gamma(p+\frac{1}{2}+\theta)\Gamma(m-p+\sigma-\theta)\Gamma(n-p+\tau-\theta)}
			{p!(m-p)! (n-p)! \Gamma(\sigma-\theta)\Gamma(\tau-\theta)} ,
		\end{gather*}
	and, if $k=2m+1$ and $\ell=2n+1$, then
		\begin{gather*}
			\hat c_{k,\ell}=(-1)^{m+n}s^{\frac{v}{2}}
			\sqrt{\frac{m! n!}{\Gamma(m+\frac{3}{2}+\sigma)\Gamma(n+\frac{3}{2}+\tau)}}\\
\hphantom{\hat c_{k,\ell}=}{}
 \times\sum_{p=0}^{\min\{m,n\}}
			\frac{\Gamma(p+\frac{3}{2}+\theta)\Gamma(m-p+\sigma-\theta)\Gamma(n-p+\tau-\theta)}
			{p!(m-p)! (n-p)! \Gamma(\sigma-\theta)\Gamma(\tau-\theta)} .
		\end{gather*}
\end{Proposition}

\begin{proof}
	The result is proved by induction on $k$ and $\ell$. First, consider the case $\ell=0$. When $k=\ell=0$, the result is given by~\eqref{hat c_0,0, case where sigma ne theta ne tau}. Now, take any $k=2m>0$, and assume that the result holds for all $\hat c_{k',0}$ with $k'=2m'<k$. Then, by Lemma~\ref{l: hat c_k,0, case where sigma ne theta ne tau} and~\eqref{v sum_i=0^p frac Gamma(i+v) i!},
		\begin{gather*}
			\hat c_{k,0} =\frac{\sigma-\theta}{\sqrt{m}}\sum_{i=0}^{m-1}(-1)^{m-i}
			\sqrt{\frac{(m-1)! \Gamma(i+\frac{1}{2}+\sigma)}{i! \Gamma(m+\frac{1}{2}+\sigma)}}\\
\hphantom{\hat c_{k,0} =}{}
\times
			(-1)^is^{\frac{v}{2}}
			\sqrt{\frac{1}{i!\Gamma(i+\frac{1}{2}+\sigma)\Gamma(\frac{1}{2}+\tau)}}
			\frac{\Gamma(\frac{1}{2}+\theta)\Gamma(i+\sigma-\theta)}{\Gamma(\sigma-\theta)}\\
\hphantom{\hat c_{k,0}}{}
			=(-1)^ms^{\frac{v}{2}}
			\sqrt{\frac{m!}{\Gamma(m+\frac{1}{2}+\sigma)\Gamma(\frac{1}{2}+\tau)}}
			\frac{\Gamma(\frac{1}{2}+\theta)(\sigma-\theta)}{m}
			\sum_{i=0}^{m-1}\frac{\Gamma(i+\sigma-\theta)}{i! \Gamma(\sigma-\theta)}\\
\hphantom{\hat c_{k,0}}{}
			=(-1)^ms^{\frac{v}{2}}
			\sqrt{\frac{1}
			{m!\Gamma(m+\frac{1}{2}+\sigma)\Gamma(\frac{1}{2}+\tau)}}
			\frac{\Gamma(\frac{1}{2}+\theta)\Gamma(m+\sigma-\theta)}{\Gamma(\sigma-\theta)}.
		\end{gather*}
	
	From the case $\ell=0$, the result also follows for the case $k=0$ by~\eqref{hat c_sigma,tau,xi,k,ell = hat c_tau,sigma,xi,ell,k}.
	
	Now, take $k=2m>0$ and $\ell=2n>0$, and assume that the result holds for all $\hat c_{k',\ell'}$ with $k'<k$ and $\ell'\le\ell$. By Lemma~\ref{l: hat c_k,ell, case k = 2m > 0 and ell = 2n > 0 and sigma ne theta ne tau},
		\begin{gather*}
			\hat c_{k,\ell} =\sqrt{\frac{n}{m}} (-1)^{m+n-2}s^{\frac{v}{2}}
			\sqrt{\frac{(m-1)! (n-1)!}
			{\Gamma(m+\frac{1}{2}+\sigma)\Gamma(n+\frac{1}{2}+\tau)}}\\
\hphantom{\hat c_{k,\ell} =}{}
\times\sum_{q=0}^{\min\{m-1,n-1\}}
			\frac{\Gamma(q+\frac{3}{2}+\theta)\Gamma(m-1-q+\sigma-\theta)\Gamma(n-1-q+\tau-\theta)}
			{q!(m-1-q)! (n-1-q)! \Gamma(\sigma-\theta)\Gamma(\tau-\theta)}\\
\hphantom{\hat c_{k,\ell} =}{}
+\frac{\sigma-\theta}{m}
			\sum_{i=0}^{m-1}(-1)^{m-i}
			\sqrt{\frac{m! \Gamma(i+\frac{1}{2}+\sigma)}{i! \Gamma(m+\frac{1}{2}+\sigma)}}
(-1)^{i+n}s^{\frac{v}{2}}
			\sqrt{\frac{i! n!}
			{\Gamma(i+\frac{1}{2}+\sigma)\Gamma(n+\frac{1}{2}+\tau)}}\\
\hphantom{\hat c_{k,\ell} =}{}
			\times\sum_{p=0}^{\min\{i,n\}}
			\frac{\Gamma(p+\frac{1}{2}+\theta)\Gamma(i-p+\sigma-\theta)\Gamma(n-p+\tau-\theta)}
			{p!(i-p)!(n-p)!\Gamma(\sigma-\theta)\Gamma(\tau-\theta)}\\
\hphantom{\hat c_{k,\ell}}{}
			=(-1)^{m+n}s^{\frac{v}{2}}\frac{1}{m}
			\sqrt{\frac{m!n!}
			{\Gamma(m+\frac{1}{2}+\sigma)\Gamma(n+\frac{1}{2}+\tau)}}\\
\hphantom{\hat c_{k,\ell} =}{}
			 \times\Bigg(\sum_{q=0}^{\min\{m-1,n-1\}}
			\frac{\Gamma(q+\frac{3}{2}+\theta)\Gamma(m-1-q+\sigma-\theta)\Gamma(n-1-q+\tau-\theta)}
			{q!(m-1-q)!(n-1-q)!\Gamma(\sigma-\theta)\Gamma(\tau-\theta)}\\
\hphantom{\hat c_{k,\ell} =}{}
+(\sigma-\theta)\sum_{i=0}^{m-1}
			\sum_{p=0}^{\min\{i,n\}}\frac{\Gamma(p+\frac{1}{2}+\theta)\Gamma(i-p+\sigma-\theta)
			\Gamma(n-p+\tau-\theta)}{p!(i-p)! (n-p)!
			\Gamma(\sigma-\theta)\Gamma(\tau-\theta)}\Bigg) .
		\end{gather*}
	Then the desired expression for $\hat c_{k,\ell}$ follows because
		\begin{gather*}
			\sum_{q=0}^{\min\{m-1,n-1\}}
			\frac{\Gamma(q+\frac{3}{2}+\theta)\Gamma(m-1-q+\sigma-\theta)\Gamma(n-1-q+\tau-\theta)}
			{q!(m-1-q)!(n-1-q)! \Gamma(\sigma-\theta)\Gamma(\tau-\theta)}\\
\qquad{}
			=\sum_{p=0}^{\min\{m,n\}}
			\frac{p\Gamma(p+\frac{1}{2}+\theta)\Gamma(m-p+\sigma-\theta)\Gamma(n-p+\tau-\theta)}
			{p!(m-p)!(n-p)! \Gamma(\sigma-\theta)\Gamma(\tau-\theta)} ,
		\end{gather*}
	and, by~\eqref{v sum_i=0^p frac Gamma(i+v) i!},
		\begin{gather}
			(\sigma-\theta)\sum_{i=0}^{m-1}
			\sum_{p=0}^{\min\{i,n\}}\frac{\Gamma(p+\frac{1}{2}+\theta)\Gamma(i-p+\sigma-\theta)
			\Gamma(n-p+\tau-\theta)}{p!(i-p)! (n-p)!
			\Gamma(\sigma-\theta)\Gamma(\tau-\theta)}\nonumber\\
\qquad{}
			=(\sigma-\theta)\sum_{p=0}^{\min\{m-1,n\}}\sum_{j=0}^{m-1-p}
			\frac{\Gamma(p+\frac{1}{2}+\theta)\Gamma(j+\sigma-\theta)
			\Gamma(n-p+\tau-\theta)}{p!j! (n-p)!
			\Gamma(\sigma-\theta)\Gamma(\tau-\theta)}\nonumber\\
\qquad{}
			=\sum_{p=0}^{\min\{m,n\}}\frac{\Gamma(p+\frac{1}{2}+\theta)(m-p)\Gamma(m-p+\sigma-\theta)
			\Gamma(n-p+\tau-\theta)}{p!(m-p)! (n-p)!
			\Gamma(\sigma-\theta)\Gamma(\tau-\theta)} .\label{(sigma - 1/2) sum_i=0^m-1 ...}
		\end{gather}
	
	Finally, take $k=2m+1$ and $\ell=2n+1$, and assume that the result holds for all $\hat c_{k',\ell'}$ with $k'<k$ and $\ell'<\ell$. By Lemma~\ref{l: hat c_k,ell, case k = 2m+1 and ell = 2n+1 and sigma ne theta ne tau},
		\begin{gather*}
			\hat c_{k,\ell} =\frac{(m+\frac{1}{2}+\theta)(-1)^{m+n}s^{\frac{v}{2}}}
			{\sqrt{(m+\frac{1}{2}+\sigma)(n+\frac{1}{2}+\tau)}}
			\sqrt{\frac{m! n!}{\Gamma(m+\frac{1}{2}+\sigma)\Gamma(n+\frac{1}{2}+\tau)}}\\
\hphantom{\hat c_{k,\ell} =}{}
\times\sum_{p=0}^{\min\{m,n\}}
			\frac{\Gamma(p+\frac{1}{2}+\theta)\Gamma(m-p+\sigma-\theta)\Gamma(n-p+\tau-\theta)}
			{p!(m-p)! (n-p)! \Gamma(\sigma-\theta)\Gamma(\tau-\theta)}\\
\hphantom{\hat c_{k,\ell} =}{}
 -\frac{\sigma-\theta}{\sqrt{n+\frac{1}{2}+\tau}}
			\sum_{i=0}^{m-1}(-1)^{m-i}
			\sqrt{\frac{m! \Gamma(i+\frac{1}{2}+\sigma)}{i! \Gamma(m+\frac{3}{2}+\sigma)}}\\
\hphantom{\hat c_{k,\ell} =}{}\times
(-1)^{i+n}s^{\frac{v}{2}}
			\sqrt{\frac{i! n!}{\Gamma(i+\frac{1}{2}+\sigma)\Gamma(n+\frac{1}{2}+\tau)}}\\
\hphantom{\hat c_{k,\ell} =}{}\times
\sum_{p=0}^{\min\{i,n\}}
			\frac{\Gamma(p+\frac{1}{2}+\theta)\Gamma(i-p+\sigma-\theta)\Gamma(n-p+\tau-\theta)}
			{p!(i-p)!(n-p)!\Gamma(\sigma-\theta)\Gamma(\tau-\theta)}\\
\hphantom{\hat c_{k,\ell}}{}
			=(-1)^{m+n}s^{\frac{v}{2}}
			\sqrt{\frac{m! n!}{\Gamma(m+\frac{3}{2}+\sigma)\Gamma(n+\frac{3}{2}+\tau)}}\\
\hphantom{\hat c_{k,\ell} =}{}
\times\Bigg(\sum_{p=0}^{\min\{m,n\}}
			\frac{(m+\frac{1}{2}+\theta)\Gamma(p+\frac{1}{2}+\theta)\Gamma(m-p+\sigma-\theta)\Gamma(n-p+\tau-\theta)}
			{p!(m-p)!(n-p)! \Gamma(\sigma-\theta)\Gamma(\tau-\theta)}\\
\hphantom{\hat c_{k,\ell} =}{}
-(\sigma-\theta)\sum_{i=0}^{m-1}
			\sum_{p=0}^{\min\{i,n\}}\frac{\Gamma(p+\frac{1}{2}+\theta)\Gamma(i-p+\sigma-\theta)
			\Gamma(n-p+\tau-\theta)}{p!(i-p)! (n-p)!
			\Gamma(\sigma-\theta)\Gamma(\tau-\theta)}\Bigg) .
\end{gather*}
	Then we get the stated expression for $\hat c_{k,\ell}$ using~\eqref{(sigma - 1/2) sum_i=0^m-1 ...} again.
\end{proof}

\begin{Remark}\label{r: hat c_k,ell, case where sigma ne theta ne tau}
	By~\eqref{hat c_sigma,tau,xi,k,ell = hat c_sigma+1,tau+1,xi+1,k-1,ell-1}, if $k$ and $\ell$ are odd, then Proposition~\ref{p: hat c_k,ell, case where sigma ne theta ne tau} also holds when $\sigma,\tau>-\frac{3}{2}$.
\end{Remark}

\section[The sesquilinear form ${\mathfrak{t}}'$]{The sesquilinear form $\boldsymbol{{\mathfrak{t}}'}$}\label{s: ft'}

Consider the notation of Section~\ref{s: mixed}. Since $x^{-1}{\mathcal{S}}_{\text{\rm odd}}={\mathcal{S}}_{\text{\rm ev}}$, a sesquilinear form ${\mathfrak{t}}'$ in $L^2_{\sigma,\tau}$, with ${\mathsf{D}}({\mathfrak{t}}')={\mathcal{S}}$, is def\/ined by
	\begin{gather*}
		{\mathfrak{t}}'(\phi,\psi)=\big\langle\phi_{\text{\rm ev}},x^{-1}\psi_{\text{\rm odd}}\big\rangle_\theta
		=\langle x\phi_{\text{\rm ev}},\psi_{\text{\rm odd}}\rangle_{\theta-1}.
	\end{gather*}
Note that ${\mathfrak{t}}'$ is neither symmetric nor bounded from the left. The goal of this section is to study~${\mathfrak{t}}'$, and use it to prove Theorem~\ref{t: VV}.

Let $c'_{k,\ell}={\mathfrak{t}}'(\phi_{\sigma,k},\phi_{\tau,\ell})$. Clearly, $c'_{k,\ell}=0$ if $k$ is odd or $\ell$ is even.

\subsection[Case where $\sigma=\theta=\tau$]{Case where $\boldsymbol{\sigma=\theta=\tau}$}\label{ss: ft', case where sigma = theta = tau}

In this case, we have $v=0$.

\begin{Proposition}\label{p: c'_k,ell, case where sigma = theta = tau}
	For $k=2m$ and $\ell=2n+1$, if $k>\ell$ {\rm(}$m>n${\rm)}, then $c'_{k,\ell}=0$, and, if $k<\ell$ {\rm(}$m\le n${\rm)}, then
		\begin{gather*}
			c'_{k,\ell}=(-1)^{n-m}s^{\frac{1}{2}}
			\sqrt{\frac{n! \Gamma(m+\frac{1}{2}+\sigma)}{m! \Gamma(n+\frac{3}{2}+\sigma)}} .
		\end{gather*}
\end{Proposition}

\begin{proof}
	This follows from~\eqref{x^-1 p_k} since $\hat c_{k,\ell}=\delta_{k,\ell}$ in this case.
\end{proof}

\begin{Proposition}\label{p: |c'_k,ell|, case where sigma = theta = tau}
	There is some $\omega=\omega(\sigma)>0$ so that, for $k=2m$ and $\ell=2n+1$,
		\begin{gather*}
			|c'_{k,\ell}|\preccurlyeq s^{\frac{1}{2}}(m+1)^{-\omega}(n+1)^{-\omega}.
		\end{gather*}
\end{Proposition}

\begin{proof}
	We can assume that $m\le n$ according to Proposition~\ref{p: c'_k,ell, case where sigma = theta = tau}. Moreover
		\begin{gather*}
			|c'_{k,\ell}|\preccurlyeq s^{\frac{1}{2}}(m+1)^{\frac{\sigma}{2}-\frac{1}{4}}(n+1)^{-\frac{\sigma}{2}-\frac{1}{4}}
		\end{gather*}
	for all $m\le n$ by Proposition~\ref{p: c'_k,ell, case where sigma = theta = tau} and Lemma~\ref{l: Gautschi}. Therefore the result follows using Lemma~\ref{l: (m+1)^alpha (n+1)^beta (m-n+1)^gamma}, reversing the roles of $m$ and $n$, because $-\frac{\sigma}{2}-\frac{1}{4}<-\frac{u}{2}<0$.
\end{proof}

\subsection[Case where $\sigma=\theta\ne\tau$ and $\tau-\sigma\not\in-{\mathbb{N}}$]{Case where $\boldsymbol{\sigma=\theta\ne\tau}$ and $\boldsymbol{\tau-\sigma\not\in-{\mathbb{N}}}$}\label{ss: ft', case where sigma = theta ne tau}

Recall that $v=\tau-\sigma$ in this case. Moreover $c'_{k,\ell}=0$ if $k>\ell$ by~\eqref{x^-1 p_k} and Corollary~\ref{c: hat c_k,ell = 0 if k >ell, case where sigma = theta ne tau}.

\begin{Proposition}\label{p: c'_k,ell, case where sigma = theta ne tau}
	For $k=2m<\ell=2n+1$ {\rm(}$m\le n${\rm)},
		\begin{gather*}
			c'_{k,\ell}=(-1)^{m+n}s^{\frac{1+v}{2}}
			\sqrt{\frac{n!\Gamma(m+\frac{1}{2}+\sigma)}{m!\Gamma(n+\frac{3}{2}+\tau)}}
			\frac{\Gamma(n-m+1+v)}{(n-m)!\Gamma(1+v)}.
		\end{gather*}
\end{Proposition}

\begin{proof}
	By~\eqref{x^-1 p_k}, Corollary~\ref{c: hat c_k,ell = 0 if k >ell, case where sigma = theta ne tau}, Proposition~\ref{p: hat c_k,ell, case where sigma = theta ne tau} and~\eqref{v sum_i=0^p frac Gamma(i+v) i!},
		\begin{gather*}
			c'_{k,\ell} =s^{\frac{1}{2}}\sum_{j=m}^n(-1)^{n-j}
 			\sqrt{\frac{n! \Gamma(j+\frac{1}{2}+\tau)}{j! \Gamma(n+\frac{3}{2}+\tau)}}(-1)^{m+j}s^{\frac{v}{2}}
			\sqrt{\frac{j! \Gamma(m+\frac{1}{2}+\sigma)}{m! \Gamma(j+\frac{1}{2}+\tau)}}
			\frac{\Gamma(j-m+v)}{(j-m)! \Gamma(v)}\\
\hphantom{c'_{k,\ell}}{}
			=(-1)^{m+n}s^{\frac{1+v}{2}}\sqrt{\frac{n! \Gamma(m+\frac{1}{2}+\sigma)}
			{m! \Gamma(n+\frac{3}{2}+\tau)}}\frac{1}{\Gamma(v)}\sum_{i=0}^{n-m}\frac{\Gamma(i+v)}{i!}\\
\hphantom{c'_{k,\ell}}{}
			=(-1)^{m+n}s^{\frac{1+v}{2}}
			\sqrt{\frac{n! \Gamma(m+\frac{1}{2}+\sigma)}{m! \Gamma(n+\frac{3}{2}+\tau)}}
			\frac{\Gamma(n-m+1+v)}{(n-m)! \Gamma(1+v)}.\tag*{\qed}
		\end{gather*}
\renewcommand{\qed}{}
\end{proof}

\begin{Proposition}\label{p: |c'_k,ell|, case where sigma = theta ne tau}
	If $(\sigma,\tau)$ satisfies~\eqref{VV, a}, then there is some $\omega=\omega(\sigma,\tau)>0$ so that, for $k=2m<\ell=2n+1$,
		\begin{gather*}
			|c'_{k,\ell}|\preccurlyeq s^{\frac{1+v}{2}}(m+1)^{-\omega}(n+1)^{-\omega}.
		\end{gather*}
\end{Proposition}

\begin{proof}
	By Proposition~\ref{p: c'_k,ell, case where sigma = theta ne tau} and Lemma~\ref{l: Gautschi},
		\begin{gather*}
			|c'_{k,\ell}|\preccurlyeq s^{\frac{1+v}{2}}(m+1)^{\frac{\sigma}{2}-\frac{1}{4}}(n+1)^{-\frac{\tau}{2}-\frac{1}{4}}(n-m+1)^v.
		\end{gather*}
	Then the result follows by Lemma~\ref{l: (m+1)^alpha (n+1)^beta (m-n+1)^gamma}, interchanging the roles of $m$ and $n$, using the condition of Theorem~\ref{t: VV}(a).
\end{proof}

\subsection[Case where $\sigma\ne\theta=\tau$ and $\sigma-\theta\not\in -{\mathbb{N}}$]{Case where $\boldsymbol{\sigma\ne\theta=\tau}$ and $\boldsymbol{\sigma-\theta\not\in -{\mathbb{N}}}$}\label{ss: ft', case where sigma ne theta = tau}

Recall that $v=\sigma-\tau$ in this case.

\begin{Proposition}\label{p: c'_k,ell, case where sigma ne theta = tau}
	For $k=2m$ and $\ell=2n+1$,
		\begin{gather*}
			c'_{k,\ell}=(-1)^{m+n}s^{\frac{1+v}{2}}
			\sqrt{\frac{m! n!}{\Gamma(m+\frac{1}{2}+\sigma)\Gamma(n+\frac{3}{2}+\tau)}}
			\sum_{j=0}^{\min\{m,n\}}\frac{\Gamma(j+\frac{1}{2}+\tau)\Gamma(m-j+v)}{j! (m-j)! \Gamma(v)} .
		\end{gather*}
\end{Proposition}

\begin{proof}
	Let $j$ run from $0$ to $\min\{m,n\}$. By~\eqref{x^-1 p_k}, Corollary~\ref{c: hat c_k,ell = 0 if k >ell, case where sigma = theta ne tau}, Proposition~\ref{p: hat c_k,ell, case where sigma = theta ne tau} and~\eqref{hat c_sigma,tau,xi,k,ell = hat c_tau,sigma,xi,ell,k},
		\begin{gather*}
			c'_{k,\ell} =s^{\frac{1}{2}}\sum_j(-1)^{n-j}
 			\sqrt{\frac{n! \Gamma(j+\frac{1}{2}+\tau)}{j! \Gamma(n+\frac{3}{2}+\tau)}}
			(-1)^{j+m}s^{\frac{v}{2}}
			\sqrt{\frac{m! \Gamma(j+\frac{1}{2}+\tau)}
			{j! \Gamma(m+\frac{1}{2}+\sigma)}}\frac{\Gamma(m-j+v)}{(m-j)! \Gamma(v)}\\
\hphantom{c'_{k,\ell}}{}
=(-1)^{m+n}s^{\frac{1+v}{2}}
			\sqrt{\frac{m! n!}{\Gamma(m+\frac{1}{2}+\sigma)\Gamma(n+\frac{3}{2}+\tau)}}
			\sum_j\frac{\Gamma(j+\frac{1}{2}+\tau)\Gamma(m-j+v)}
			{j! (m-j)! \Gamma(v)}. \tag*{\qed}
		\end{gather*}
\renewcommand{\qed}{}
\end{proof}

\begin{Proposition}\label{p: |c'_k,ell|, case where sigma ne theta = tau}
	If $(\sigma,\tau)$ satisfies~\eqref{VV, b}, then there is some $\omega=\omega(\sigma,\tau)>0$ so that, for $k=2m$ and $\ell=2n+1$,
		\begin{gather*}
			|c'_{k,\ell}|\preccurlyeq s^{\frac{1+v}{2}}(m+1)^{-\omega}(n+1)^{-\omega}.
		\end{gather*}
\end{Proposition}

\begin{proof}
	By Proposition~\ref{p: c'_k,ell, case where sigma ne theta = tau} and Lemma~\ref{l: Gautschi},
		\begin{gather*}
			|c'_{k,\ell}|\preccurlyeq s^{\frac{1+v}{2}}(m+1)^{\frac{1}{4}-\frac{\sigma}{2}}(n+1)^{-\frac{1}{4}-\frac{\tau}{2}}
			\sum_{j=0}^{\min\{m,n\}}(m-j+1)^{v-1}(j+1)^{\tau-\frac{1}{2}}.
		\end{gather*}
	Then the result follows by Corollary~\ref{c: check fJ}, proved in Section~\ref{s: main estimates}, since $(\sigma,\tau)$ satisf\/ies~\eqref{VV, b}.
\end{proof}

\subsection[Case where $\sigma\ne\theta=\tau+1$ and $\sigma-\tau-1\not\in-{\mathbb{N}}$]{Case where $\boldsymbol{\sigma\ne\theta=\tau+1}$ and $\boldsymbol{\sigma-\tau-1\not\in-{\mathbb{N}}}$}\label{ss: sigma ne theta = tau+1}

Note that $v=\sigma-\tau-2$ in this case. Moreover
	\begin{gather}\label{c'_k,ell, case sigma ne theta = tau+1, 1st}
		c'_{k,\ell}=\big\langle\phi_{\sigma,k},x^{-1}\phi_{\tau,\ell}\big\rangle_{\tau+1}
		=\langle x\phi_{\sigma,k},\phi_{\tau,\ell}\rangle_\tau=\langle\phi_{\tau,\ell},x\phi_{\sigma,k}\rangle_\tau
	\end{gather}
for $k=2m$ and $\ell=2n+1$ (Remark~\ref{r: VV}(ii)).

\begin{Proposition}\label{p: c'_k,ell, case sigma ne theta = tau+1}
	Let $k=2m$ and $\ell=2n+1$. If $k+1<\ell$ {\rm(}$m<n${\rm)}, then $c'_{k,\ell}=0$. If $k+1\ge\ell$ {\rm(}$m\ge n${\rm)}, then
		\begin{gather*}
			c'_{k,\ell}=(-1)^{m+n}s^{\frac{v+1}{2}}
			\sqrt{\frac{m!\Gamma(n+\frac{3}{2}+\tau)}{n!\Gamma(m+\frac{1}{2}+\sigma)}}
			\frac{\Gamma(m-n+v+1)}{(m-n)!\Gamma(v+1)}.
		\end{gather*}
\end{Proposition}

\begin{proof}
	By~\eqref{recursion} and~\eqref{c'_k,ell, case sigma ne theta = tau+1, 1st},
		\begin{gather}\label{c'_k,ell, case sigma ne theta = tau+1, 2nd}
			c'_{k,\ell}=\sqrt{\frac{m+\frac{1}{2}+\sigma}{s}}\hat c_{\tau,\sigma,\tau,\ell,k+1}
			+\sqrt{\frac{m}{s}}\hat c_{\tau,\sigma,\tau,\ell,k-1}.
		\end{gather}
	So $c'_{k,\ell}=0$ if $k+1<\ell$ by Corollary~\ref{c: hat c_k,ell = 0 if k >ell, case where sigma = theta ne tau}. When $k+1=\ell$ ($m=n$), by~\eqref{c'_k,ell, case sigma ne theta = tau+1, 2nd} and Proposition~\ref{p: hat c_k,ell, case where sigma = theta ne tau},
		\begin{gather*}
			c'_{k,\ell}=\sqrt{\frac{m+\frac{1}{2}+\sigma}{s}}s^{\frac{v+2}{2}}
			\sqrt{\frac{\Gamma(n+\frac{3}{2}+\tau)}{\Gamma(m+\frac{3}{2}+\sigma)}}
			=s^{\frac{v+1}{2}}\sqrt{\frac{\Gamma(n+\frac{3}{2}+\tau)}{\Gamma(m+\frac{1}{2}+\sigma)}}.
		\end{gather*}
	When $k-1\ge\ell$ ($m>n$), by~\eqref{c'_k,ell, case sigma ne theta = tau+1, 2nd} and Proposition~\ref{p: hat c_k,ell, case where sigma = theta ne tau},
		\begin{gather*}
			c'_{k,\ell} =\sqrt{\frac{m+\frac{1}{2}+\sigma}{s}}(-1)^{m+n}s^{\frac{v+2}{2}}
			\sqrt{\frac{m! \Gamma(n+\frac{3}{2}+\tau)}{n! \Gamma(m+\frac{3}{2}+\sigma)}}
			\frac{\Gamma(m-n+v+2)}{(m-n)! \Gamma(v+2)}\\
\hphantom{c'_{k,\ell} =}{}
 +\sqrt{\frac{m}{s}}(-1)^{m+n-1}s^{\frac{v+2}{2}}
			\sqrt{\frac{(m-1)! \Gamma(n+\frac{3}{2}+\tau)}{n! \Gamma(m+\frac{1}{2}+\sigma)}}
			\frac{\Gamma(m-n+v+1)}{(m-1-n)! \Gamma(v+2)}\\
\hphantom{c'_{k,\ell} }{}			
=(-1)^{m+n}s^{\frac{v+1}{2}}\sqrt{\frac{m!\Gamma(n+\frac{3}{2}+\tau)}
			{n!\Gamma(m+\frac{1}{2}+\sigma)}}
			\frac{\Gamma(m-n+v+1)}{(m-1-n)!\Gamma(v+2)}\left(\frac{m-n+v+1}{m-n}-1\right)\\
\hphantom{c'_{k,\ell} }{}	
			=(-1)^{m+n}s^{\frac{v+1}{2}}\sqrt{\frac{m!\Gamma(n+\frac{3}{2}+\tau)}
			{n!\Gamma(m+\frac{1}{2}+\sigma)}}
			\frac{\Gamma(m-n+v+1)}{(m-n)!\Gamma(v+1)}. \tag*{\qed}
		\end{gather*}
\renewcommand{\qed}{}
\end{proof}

\begin{Proposition}\label{p: |c'_k,ell|, case sigma ne theta = tau+1}
	If $(\sigma,\tau)$ satisfies~\eqref{VV, c}, then there is some $\omega=\omega(\sigma,\tau)>0$ so that, for $k=2m$ and $\ell=2n+1$,
		\begin{gather*}
			|c'_{k,\ell}|\preccurlyeq s^{\frac{v+1}{2}}(m+1)^{-\omega}(n+1)^{-\omega} .
		\end{gather*}
\end{Proposition}

\begin{proof}
	By Proposition~\ref{p: c'_k,ell, case sigma ne theta = tau+1}, we can assume that $k+1\ge\ell$ ($m\ge n$), and, in this case, using also Lemma~\ref{l: Gautschi}, we get
		\begin{gather*}
			|c'_{k,\ell}|\preccurlyeq s^{\frac{v+1}{2}}(m+1)^{\frac{1}{4}-\frac{\sigma}{2}}(n+1)^{\frac{1}{4}+\frac{\tau}{2}}(m-n+1)^v .
		\end{gather*}
	Then the result follows using Lemma~\ref{l: (m+1)^alpha (n+1)^beta (m-n+1)^gamma}.
\end{proof}

\subsection[Case where $\sigma\ne\theta\ne\tau$ and $\sigma-\theta,\tau-\theta\not\in-{\mathbb{N}}$]{Case where $\boldsymbol{\sigma\ne\theta\ne\tau}$ and $\boldsymbol{\sigma-\theta,\tau-\theta\not\in-{\mathbb{N}}}$}
\label{ss: ft', case where sigma ne theta ne tau}

\begin{Proposition}\label{p: c'_k,ell, case where sigma ne theta ne tau}
	For $k=2m$ and $\ell=2n+1$,
		\begin{gather*}
			c'_{k,\ell}=(-1)^{m+n}s^{\frac{1+v}{2}}
			\sqrt{\frac{m! n!}{\Gamma(m+\frac{1}{2}+\sigma)\Gamma(n+\frac{3}{2}+\tau)}}\\
\hphantom{c'_{k,\ell}=}{}
\times\sum_{p=0}^{\min\{m,n\}}
			\frac{\Gamma(p+\frac{1}{2}+\theta)\Gamma(m-p+\sigma-\theta)\Gamma(n-p+1+\tau-\theta)}
			{p!(m-p)! (n-p)! \Gamma(\sigma-\theta)\Gamma(1+\tau-\theta)} .
		\end{gather*}
\end{Proposition}

\begin{proof}
	By~\eqref{x^-1 p_k} and Proposition~\ref{p: hat c_k,ell, case where sigma ne theta ne tau},
		\begin{gather*}
			c'_{k,\ell} =s^{\frac{1}{2}}\sum_{j=0}^n(-1)^{n-j}
 			\sqrt{\frac{n! \Gamma(j+\frac{1}{2}+\tau)}{j! \Gamma(n+\frac{3}{2}+\tau)}}
			(-1)^{m+j}s^{\frac{v}{2}}
			\sqrt{\frac{m! j!}
			{\Gamma(m+\frac{1}{2}+\sigma)\Gamma(j+\frac{1}{2}+\tau)}}\\
\hphantom{c'_{k,\ell} =}{}
\times\sum_{p=0}^{\min\{m,j\}}
			\frac{\Gamma(p+\frac{1}{2}+\theta)\Gamma(m-p+\sigma-\theta)\Gamma(j-p+\tau-\theta)}
			{p!(m-p)! (j-p)! \Gamma(\sigma-\theta)\Gamma(\tau-\theta)}\\
\hphantom{c'_{k,\ell}}{}
			=(-1)^{m+n}s^{\frac{1+v}{2}}\sqrt{\frac{m!n!}
			{\Gamma(m+\frac{1}{2}+\sigma)\Gamma(n+\frac{3}{2}+\tau)}}\\
\hphantom{c'_{k,\ell} =}{}
\times\sum_{j=0}^n\sum_{p=0}^{\min\{m,j\}}
			\frac{\Gamma(p+\frac{1}{2}+\theta)\Gamma(m-p+\sigma-\theta)\Gamma(j-p+\tau-\theta)}
			{p!(m-p)!(j-p)!\Gamma(\sigma-\theta)\Gamma(\tau-\theta)}.
		\end{gather*}
	But, by~\eqref{v sum_i=0^p frac Gamma(i+v) i!},
		\begin{gather*}
			\sum_{j=0}^n\sum_{p=0}^{\min\{m,j\}}
			\frac{\Gamma(m-p+\sigma-\theta)\Gamma(j-p+\tau-\theta)}
			{(m-p)! (j-p)! \Gamma(\sigma-\theta)\Gamma(\tau-\theta)}\\
					\qquad {} =\sum_{p=0}^{\min\{m,n\}}\sum_{j=p}^n
					\frac{\Gamma(m-p+\sigma-\theta)\Gamma(j-p+\tau-\theta)}
					{(m-p)! (j-p)! \Gamma(\sigma-\theta)\Gamma(\tau-\theta)} \\
\qquad{}
					 =\sum_{p=0}^{\min\{m,n\}}\sum_{i=0}^{n-p}
					\frac{\Gamma(m-p+\sigma-\theta)\Gamma(i+\tau-\theta)}
					{(m-p)! i! \Gamma(\sigma-\theta)\Gamma(\tau-\theta)}\\
\qquad{}
					 =\sum_{p=0}^{\min\{m,n\}}
					\frac{\Gamma(m-p+\sigma-\theta)\Gamma(n-p+1+\tau-\theta)}
					{(m-p)! (n-p)! \Gamma(\sigma-\theta)\Gamma(1+\tau-\theta)} . \tag*{\qed}
\end{gather*}
\renewcommand{\qed}{}
\end{proof}

\begin{Proposition}\label{p: |c'_k,ell|, case where sigma ne theta ne tau}
	If $(\sigma,\tau,\theta)$ satisfies~\eqref{VV, d}, then there is some $\omega=\omega(\sigma,\tau,\theta)>0$ so that, for $k=2m$ and $\ell=2n+1$,
		\begin{gather*}
			|c'_{k,\ell}|\preccurlyeq s^{\frac{1+v}{2}}(m+1)^{-\omega}(n+1)^{-\omega} .
		\end{gather*}
\end{Proposition}

\begin{proof}
	Let $p$ run from $0$ to $\min\{m,n\}$. By Proposition~\ref{p: c'_k,ell, case where sigma ne theta ne tau} and Lemma~\ref{l: Gautschi},
		\begin{gather*}
			|c'_{k,\ell}|\preccurlyeq s^{\frac{1+v}{2}}(m+1)^{\frac{1}{4}-\frac{\sigma}{2}}(n+1)^{-\frac{1}{4}-\frac{\tau}{2}}
			\sum_p(m-p+1)^{\sigma-\theta-1}(n-p+1)^{\tau-\theta}(p+1)^{\theta-\frac{1}{2}}.
		\end{gather*}
	Then the result follows by Corollary~\ref{c: fK_conv cap fK'_conv}, proved in Section~\ref{s: main estimates}, since $(\sigma,\tau,\theta)$ satisf\/ies~\eqref{VV, d}.
\end{proof}

\subsection{Proof of Theorem~\ref{t: VV}}\label{ss: proof main}

Assume the conditions of Theorem~\ref{t: VV}. Let ${\mathfrak{j}}_{\sigma,\tau}$ be the positive def\/inite symmetric sesquilinear form in $L^2_{\sigma,\tau}$, with domain ${\mathcal{S}}$, def\/ined by ${\mathfrak{j}}_{\sigma,\tau}(\phi,\psi)=\langle J_{\sigma,\tau}\phi,\psi\rangle_{\sigma,\tau}$.

\begin{Proposition}\label{p: |ft'(phi)|}
	For any $\epsilon>0$, there is some $E=E(\epsilon,\sigma,\tau,\theta)>0$ such that, for all $\phi\in{\mathcal{S}}$,
		\begin{gather*}
			|{\mathfrak{t}}'(\phi)|\le\epsilon s^{\frac{v-1}{2}}{\mathfrak{j}}_{\sigma,\tau}(\phi)+Es^{\frac{1+v}{2}}\|\phi\|_{\sigma,\tau}^2.
		\end{gather*}
\end{Proposition}

\begin{proof}
	This follows from Propositions~\ref{p: |c'_k,ell|, case where sigma = theta = tau},~\ref{p: |c'_k,ell|, case where sigma = theta ne tau},~\ref{p: |c'_k,ell|, case where sigma ne theta = tau},~\ref{p: |c'_k,ell|, case sigma ne theta = tau+1} and~\ref{p: |c'_k,ell|, case where sigma ne theta ne tau} using the arguments of the proof of Proposition~\ref{p: ft le epsilon s^u-1 fl(phi)+C s^u |phi|_sigma^2}.
\end{proof}

\begin{proof}[Proof of Theorem~\ref{t: VV}]
	This is analogous to the proof of Theorem~\ref{t: UU}. Thus some details and the bibliographic references are omitted.
	
	Let ${\mathfrak{t}}_{\sigma,\tau}$ be the positive def\/inite symmetric sesquilinear form in $L^2_{\sigma,\tau}$, with ${\mathsf{D}}({\mathfrak{t}}_{\sigma,\tau})={\mathcal{S}}$, def\/ined by ${\mathfrak{t}}_\sigma$ on ${\mathcal{S}}_{\text{\rm ev}}$ and ${\mathfrak{t}}_\tau$ on ${\mathcal{S}}_{\text{\rm odd}}$, and vanishing on ${\mathcal{S}}_{\text{\rm ev}}\times{\mathcal{S}}_{\text{\rm odd}}$. Let ${\mathfrak{s}}$ be the symmetric sesquilinear form in $L^2_{\sigma,\tau}$, with ${\mathsf{D}}({\mathfrak{s}})={\mathcal{S}}$, def\/ined by ${\mathfrak{s}}(\phi,\psi)={\mathfrak{t}}'(\phi,\psi)+\overline{{\mathfrak{t}}'(\psi,\phi)}$. Then the symmetric sesquilinear form ${\mathfrak{v}}={\mathfrak{j}}_{\sigma,\tau}+\xi{\mathfrak{t}}_{\sigma,\tau}+\eta{\mathfrak{s}}$ in $L^2_{\sigma,\tau}$, with ${\mathsf{D}}({\mathfrak{v}})={\mathcal{S}}$, is given by the right hand side of~\eqref{langle VV^1/2 phi, VV^1/2 psi rangle_sigma,tau}. Using Propositions~\ref{p: ft le epsilon s^u-1 fl(phi)+C s^u |phi|_sigma^2} and~\ref{p: |ft'(phi)|}, for any $\epsilon>0$, there are some $C=C(\epsilon,\sigma,\tau,u)>0$ and $E=E(\epsilon,\sigma,\tau,\theta)>0$ such that, for all $\phi\in{\mathcal{S}}$,
		\begin{gather}\label{(xi ft_sigma,tau + eta fs)(phi)|}
			|(\xi{\mathfrak{t}}_{\sigma,\tau}+\eta{\mathfrak{s}}(\phi)|
			\le\epsilon\big(\xi s^{u-1}+2|\eta|s^{\frac{v-1}{2}}\big){\mathfrak{j}}_{\sigma,\tau}(\phi)
			+\big(\xi Cs^u+2|\eta|Es^{\frac{1+v}{2}}\big)\|\phi\|_{\sigma,\tau}^2.
		\end{gather}
	Then, taking $\epsilon$ so that $\epsilon\big(\xi s^{u-1}+2|\eta|s^{\frac{v-1}{2}}\big)<1$, since ${\mathfrak{j}}_{\sigma,\tau}$ is closable and positive def\/inite, it follows that ${\mathfrak{v}}$ is sectorial and closable, and ${\mathsf{D}}(\bar{\mathfrak{v}})={\mathsf{D}}(\overline{{\mathfrak{j}}_{\sigma,\tau}})$; in particular, ${\mathfrak{v}}$ is bounded from below because it is also symmetric. So $\bar{\mathfrak{v}}$ is induced by a self-adjoint operator ${\mathcal{V}}$ in $L^2_{\sigma,\tau}$ with ${\mathsf{D}}({\mathcal{V}}^{1/2})={\mathsf{D}}(\bar{\mathfrak{v}})$. Thus ${\mathcal{S}}$ is a core of $\bar{\mathfrak{v}}$ and ${\mathcal{V}}^{1/2}$.
	
	For all $\phi\in{\mathcal{S}}$,
		\begin{gather}\label{|fv(phi)| ge ...}
			{\mathfrak{v}}(\phi)\ge{\mathfrak{j}}_{\sigma,\tau}(\phi)+\xi{\mathfrak{t}}_{\sigma,\tau}(\phi)-|\eta|{\mathfrak{s}}(\phi)\ge{\mathfrak{j}}_{\sigma,\tau}(\phi)
			+\xi{\mathfrak{t}}_{\sigma,\tau}(\phi)-2|\eta||{\mathfrak{t}}'(\phi)|.
		\end{gather}
	Since ${\mathcal{S}}$ is a core of $\bar{\mathfrak{v}}$ and $\overline{{\mathfrak{j}}_{\sigma,\tau}}$, using Propositions~\ref{p: ft(phi) ge D s^u/2 (k+1)^-u |phi|_sigma^2} and~\ref{p: |ft'(phi)|} like in the proof of Theorem~\ref{t: UU}, it follows from~\eqref{|fv(phi)| ge ...} that ${\mathcal{V}}$ has a discrete spectrum, which consists of two groups of eigenvalues, $\lambda_0\le\lambda_2\le\cdots$ and $\lambda_1\le\lambda_3\le\cdots$, repeated according to their multiplicity, satisfying~\eqref{lambda_k ge ..., case of VV}. On the other hand, by~\eqref{(xi ft_sigma,tau + eta fs)(phi)|}, for all $\phi\in{\mathcal{S}}$,
		\begin{gather}\label{|fv(phi)| le ...}
			{\mathfrak{v}}(\phi)\le\big(1+\epsilon\big(\xi s^{u-1}+2|\eta|s^{\frac{v-1}{2}}\big)\big){\mathfrak{j}}_{\sigma,\tau}(\phi)
			+\big(\xi Cs^u+2|\eta|Es^{\frac{1+v}{2}}\big)\|\phi\|_{\sigma,\tau}^2,
		\end{gather}
	obtaining~\eqref{lambda_k le ..., case of VV} because ${\mathcal{S}}$ is a core of $\bar{\mathfrak{v}}$ and $\overline{{\mathfrak{j}}_{\sigma,\tau}}$.
	
	With the notation of~(iii), let $\tilde{\mathfrak{t}}_\sigma$ (respectively, $\tilde{\mathfrak{t}}_\tau$) be the symmetric sesquilinear form in~$L^2_\sigma$ (respectively, $L^2_\tau$), with ${\mathsf{D}}(\tilde{\mathfrak{t}}_\sigma)={\mathcal{S}}$ (respectively, ${\mathsf{D}}(\tilde{\mathfrak{t}}_\tau)={\mathcal{S}}$), def\/ined like ${\mathfrak{t}}_\sigma$ (respectively,~${\mathfrak{t}}_\tau$), using $\tilde u$ (respectively, $v-\tilde u+1$) instead of $u$. Let $\tilde{\mathfrak{t}}_{\sigma,\tau}$ be the positive def\/inite symmetric sesquilinear form in $L^2_{\sigma,\tau}$, with ${\mathsf{D}}(\tilde{\mathfrak{t}}_{\sigma,\tau})={\mathcal{S}}$, def\/ined by $\tilde{\mathfrak{t}}_\sigma$ on ${\mathcal{S}}_{\text{\rm ev}}$ and $\tilde{\mathfrak{t}}_\tau$ on ${\mathcal{S}}_{\text{\rm odd}}$, and vanishing on ${\mathcal{S}}_{\text{\rm ev}}\times{\mathcal{S}}_{\text{\rm odd}}$. By the Schwartz inequality,
		\begin{gather}
			2|{\mathfrak{t}}'(\phi)|=2\big|\big\langle\phi_{\text{\rm ev}}|x|^{-\tilde u+\sigma-\theta},x^{-1}\phi_{\text{\rm odd}}|x|^{\tilde u-\sigma+\theta}\big\rangle_\theta\big|\notag\\
			\hphantom{2|{\mathfrak{t}}'(\phi)|}{}
			\le2\big\|\phi_{\text{\rm ev}}|x|^{-\tilde u+\sigma-\theta}\big\|_\theta\cdot\big\|\phi_{\text{\rm odd}}|x|^{\tilde u-\sigma+\theta-1}\big\|_\theta\notag\\
			\hphantom{2|{\mathfrak{t}}'(\phi)|}{}
			=2\big\|\phi_{\text{\rm ev}}|x|^{-\tilde u}\big\|_\sigma\cdot\big\|\phi_{\text{\rm odd}}|x|^{\tilde u-v-1}\big\|_\tau\notag\\
			\hphantom{2|{\mathfrak{t}}'(\phi)|}{}
			\le\big\|\phi_{\text{\rm ev}}|x|^{-\tilde u}\big\|_\sigma^2+\big\|\phi_{\text{\rm odd}}|x|^{\tilde u-v-1}\big\|_\tau^2
			=\tilde{\mathfrak{t}}_{\sigma,\tau}(\phi).\label{2 |ft'(phi)| le tilde ft_sigma,tau(phi)}
		\end{gather}
	Hence~\eqref{lambda_k ge ..., case of VV with tilde u} follows like in the proof of Theorem~\ref{t: UU}, using Propositions~\ref{p: ft le epsilon s^u-1 fl(phi)+C s^u |phi|_sigma^2} and~\ref{p: ft(phi) ge D s^u/2 (k+1)^-u |phi|_sigma^2}.
	
	If $u=\frac{v+1}{2}$, then we can take $\tilde u=u=v-\tilde u+1$ in~(iii), yielding
		\begin{gather*}
			{\mathfrak{v}}(\phi)\ge{\mathfrak{j}}_{\sigma,\tau}(\phi)+(\xi-|\eta|){\mathfrak{t}}_{\sigma,\tau}(\phi)
		\end{gather*}
	by~\eqref{|fv(phi)| ge ...} and~\eqref{2 |ft'(phi)| le tilde ft_sigma,tau(phi)}. Thus~\eqref{lambda_k ge ..., case of VV with u = (v+1)/2, xi ge |eta|} follows if moreover $|\eta|\le\xi$, like in the proof of Theorem~\ref{t: UU}, using Proposition~\ref{p: ft(phi) ge D s^u/2 (k+1)^-u |phi|_sigma^2}.
	
	If we add the term $\xi'\langle\phi_{\text{\rm ev}},\psi_{\text{\rm ev}}\rangle_\sigma+\xi''\langle\phi_{\text{\rm odd}},\psi_{\text{\rm odd}}\rangle_\tau$ to the right hand side of~\eqref{langle VV^1/2 phi, VV^1/2 psi rangle_sigma,tau}, for some $\xi',\xi''\in{\mathbb{R}}$, then the same argument can be used by adding the term $\xi'\|\phi_{\text{\rm ev}}\|_\sigma^2+\xi''\|\phi_{\text{\rm odd}}\|_\tau^2$ to ${\mathfrak{j}}_{\sigma,\tau}$, obtaining~(v).	
\end{proof}

\section{A preliminary estimate}\label{s: preliminary estimate}

\subsection{Statement}\label{ss: statement of the preliminary estimate}

The standard coordinates of ${\mathbb{R}}^5$ are denoted by $(\alpha,\beta,\gamma,\delta,\varkappa)$. Consider the partition of ${\mathbb{R}}$ into the following intervals: $I_1=(-\infty,-1]$, $I_2=(-1,-\frac{1}{2}]$, $I_3=(-\frac{1}{2},-\frac{1}{3}]$, $I_4=(-\frac{1}{3},0)$ and $I_5=[0,\infty)$. Let $Q_{ijk}=I_i\times I_j\times I_k$, and consider the following subsets of ${\mathbb{R}}^5$:
	\begin{description}
	
		\item[${\mathfrak{S}}_{515}$:] This is the subset of ${\mathbb{R}}^2\times Q_{515}$ def\/ined by
			\begin{gather}\label{alpha + gamma, alpha + beta + gamma + varkappa < 0}
				\alpha+\gamma,\alpha+\beta+\gamma+\varkappa<0.
			\end{gather}
			
		\item[${\mathfrak{S}}_{522}$:] This is the subset of ${\mathbb{R}}^2\times Q_{522}$ def\/ined by
			\begin{gather}\label{alpha + gamma, alpha + beta + gamma < 0}
				\alpha+\gamma,\alpha+\beta+\gamma<0.
			\end{gather}
			
		\item[${\mathfrak{S}}_{252}$:] This is the subset of ${\mathbb{R}}^2\times Q_{252}$ def\/ined by
			\begin{gather}
				0\le\gamma+\delta\Rightarrow\alpha+\gamma,\alpha+\beta+\gamma+\delta+\varkappa+1<0,
				\label{0 le gamma + delta => alpha + gamma, alpha + beta + gamma + delta + varkappa + 1 < 0}\\
				\gamma+\delta<0\Rightarrow\left\{\negmedspace\negthickspace
					\begin{array}{l}
						\alpha+\gamma+\delta,\alpha+\beta+\varkappa+1<0,\ \text{or}\\[2pt]
						\alpha+\gamma+\frac{1}{2},\alpha+\beta+\delta<0,\ \text{or}\\[2pt]
						\alpha+\gamma+\frac{1}{3},\alpha+\beta+\delta+\frac{1}{3}<0,\ \text{or}\\[2pt]
						\alpha+\gamma,\alpha+\beta+\delta+\varkappa+1<0.
					\end{array}\right.
				\label{gamma + delta < 0 => 4 cases}
			\end{gather}
			
		\item[${\mathfrak{S}}_{155}$:] This is the subset of ${\mathbb{R}}^2\times Q_{155}$ def\/ined by~\eqref{0 le gamma + delta => alpha + gamma, alpha + beta + gamma + delta + varkappa + 1 < 0} and
			\begin{gather}\label{gamma + delta < 0 => 5 cases}
				\gamma+\delta<0\Rightarrow\left\{\negmedspace\negthickspace
					\begin{array}{l}
						\alpha+\gamma+\delta,\alpha+\beta+\varkappa+1<0,\ \text{or}\\
						\alpha+\gamma+1,\alpha+\beta+\delta+\varkappa<0,\ \text{or}\\[2pt]
						\alpha+\gamma+\frac{1}{2},\alpha+\beta+\delta+\varkappa+\frac{1}{2}<0,\ \text{or}\\[2pt]
						\alpha+\gamma+\frac{1}{3},\alpha+\beta+\delta+\varkappa+\frac{2}{3}<0,\ \text{or}\\[2pt]
						\alpha+\gamma,\alpha+\beta+\delta+\varkappa+1<0.
					\end{array}\right.
			\end{gather}
			
		\item[${\mathfrak{S}}_{212}$:] This is the subset of ${\mathbb{R}}^2\times Q_{212}$ def\/ined by
			\begin{gather}
				\gamma\le\varkappa\Rightarrow\alpha+\gamma,\alpha+\beta+\varkappa<0,
				\label{gamma le varkappa => alpha + gamma, alpha + beta + varkappa < 0}\\
				\varkappa\le\gamma\Rightarrow\alpha+\gamma,\alpha+\beta+\gamma<0.
				\label{varkappa le gamma => alpha + gamma, alpha + beta + gamma < 0}
			\end{gather}
	
	\end{description}
Let $\check{\mathfrak{S}}={\mathfrak{S}}_{515}\cup{\mathfrak{S}}_{522}\cup{\mathfrak{S}}_{252}\cup{\mathfrak{S}}_{155}\cup{\mathfrak{S}}_{212}$. On the other hand, consider the linear isomorphism of ${\mathbb{R}}^5$ def\/ined by
	\begin{gather}\label{R^5 to R^5}
		(\alpha,\beta,\gamma,\delta,\varkappa)\mapsto(\beta,\alpha,\delta,\gamma,\varkappa).
	\end{gather}
This is the ref\/lection with respect to the linear subspace def\/ined by $\alpha=\beta$ and $\gamma=\delta$. The image of any subset $X\subset{\mathbb{R}}^5$ by the mapping~\eqref{R^5 to R^5} is denoted by $X'$, and let $X_{\text{\rm conv}}$ be the convex hull of $X$. Thus $X'_{\text{\rm conv}}:=(X')_{\text{\rm conv}}=(X_{\text{\rm conv}})'$.

\begin{Lemma}\label{l: preliminary estimate}
	If $(\alpha,\beta,\gamma,\delta,\varkappa)\in\check{\mathfrak{S}}_{\text{\rm conv}}\cap\check{\mathfrak{S}}'_{\text{\rm conv}}$, then there is some $\omega>0$ such that, for all $m,n\in{\mathbb{N}}$,
		\begin{gather}\label{preliminary estimate}
			(m+1)^\alpha(n+1)^\beta\sum_{p=0}^{\min\{m,n\}}(m-p+1)^\gamma(n-p+1)^\delta(p+1)^\varkappa
			\preccurlyeq(m+1)^{-\omega}(n+1)^{-\omega}.
		\end{gather}
\end{Lemma}

\subsection{Proof of Lemma~\ref{l: preliminary estimate}}\label{ss: proof of the preliminary estimate}

Since the roles of $m$ and $n$ in Lemma~\ref{l: preliminary estimate} are interchanged by the mapping~\eqref{R^5 to R^5}, we can assume that $m\ge n$. Then Lemma~\ref{l: (m+1)^alpha (n+1)^beta (m-n+1)^gamma} gives~\eqref{preliminary estimate} once
	\begin{gather*}
		\sum_{p=0}^n(m-p+1)^\gamma(n-p+1)^\delta(p+1)^\varkappa
	\end{gather*}
is appropriately estimated. Estimates of this expression are achieved with several strategies explained in Sections~\ref{sss: A_1}--\ref{sss: A_7}, giving rise to several lists of conditions that guarantee~\eqref{preliminary estimate} when $m\ge n$. Then, for the chosen subindices $ijk$ equal to 515, 522, 252, 155, 212, every ${\mathfrak{S}}_{ijk}$ is def\/ined by the most general of those conditions on ${\mathbb{R}}^2\times Q_{ijk}$. This will show that~\eqref{preliminary estimate} holds for $m\ge n$ and $(\alpha,\beta,\gamma,\delta,\varkappa)\in\check{\mathfrak{S}}$. In Section~\ref{sss: convex}, it will be shown that this property can be extended to the convex hull $\check{\mathfrak{S}}_{\text{\rm conv}}$, completing the proof of Lemma~\ref{l: preliminary estimate}.

\subsubsection{First list of conditions}\label{sss: A_1}

For all $\epsilon>0$,
		\begin{gather}
			\sum_{p=0}^n(p+1)^\varkappa=\sum_{q=1}^{n+1}q^\varkappa\le
				\begin{cases}
					\displaystyle \int_1^{n+2}x^\varkappa dx & \text{if $\varkappa\ge0$,}\\
					\displaystyle 1+\int_1^{n+1}x^\varkappa dx & \text{if $\varkappa<0$}
				\end{cases}\notag\\
			\hphantom{\sum_{p=0}^n(p+1)^\varkappa}{}
			\preccurlyeq
				\begin{cases}
					(n+1)^{\varkappa+1} & \text{if $\varkappa>-1$,}\\
					1+\ln(n+1) & \text{if $\varkappa=-1$,}\\
					1 & \text{if $\varkappa<-1$}
				\end{cases}
			\preccurlyeq
				\begin{cases}
					(n+1)^{\varkappa+1} & \text{if $\varkappa>-1$,}\\
					(n+1)^\epsilon & \text{if $\varkappa=-1$,}\\
					1 & \text{if $\varkappa<-1$}.
				\end{cases}\label{A_1, sum_p=0^n (p+1)^varkappa}
		\end{gather}
On the other hand, we claim that
	\begin{gather}\label{A_1, (m-p+1)^gamma (n-p+1)^delta preccurlyeq ...}
		(m-p+1)^\gamma(n-p+1)^\delta\preccurlyeq
			\begin{cases}
				(m+1)^\gamma(n+1)^\delta & \text{if $\delta\ge-\gamma,0$,}\\
				\left.\!
					\begin{array}{@{}l@{}}
						(m-n+1)^{\gamma+\delta}\ \text{and}\\
						(m-n+1)^\gamma(n+1)^\delta
					\end{array}
				\right\} & \text{if $0\le\delta<-\gamma$,}\\
				(m-n+1)^\gamma & \text{if $\delta\le-\gamma,0$,}\\
				\left.\!
					\begin{array}{@{}l@{}}
						(m-n+1)^\gamma\ \text{or}\\
						(m+1)^\gamma(n+1)^\delta
					\end{array}
				\right\} & \text{if $-\gamma<\delta<0$,}
			\end{cases}
	\end{gather}
for all $p=0,\dots,n$. Combining~\eqref{A_1, sum_p=0^n (p+1)^varkappa} and~\eqref{A_1, (m-p+1)^gamma (n-p+1)^delta preccurlyeq ...}, it follows that
	\begin{gather}\label{A_1}
		(m+1)^\alpha(n+1)^\beta\sum_{p=0}^n(m-p+1)^\gamma(n-p+1)^\delta(p+1)^\varkappa\preccurlyeq A_1,
	\end{gather}
where $A_1=A_1(m,n,\alpha,\beta,\gamma,\delta,\varkappa)$ can be taken to be equal to
	\begin{alignat*}{2}
		&(m+1)^{\alpha+\gamma}(n+1)^{\beta+\delta+\varkappa+1} &\quad&\text{if $-\gamma,0\le\delta$, $-1<\varkappa$},\\[2pt]
		&\left.\!\!
			\begin{array}{@{}l@{}}
				(m+1)^\alpha(n+1)^{\beta+\varkappa+1}(m-n+1)^{\gamma+\delta}\ \text{and}\\[2pt]
				(m+1)^\alpha(n+1)^{\beta+\delta+\varkappa+1}(m-n+1)^\gamma
			\end{array}
		\right\} &\quad&\text{if $0\le\delta<-\gamma$, $-1<\varkappa$}.
	\end{alignat*}
In the def\/inition of $A_1$, the other cases of $\gamma,\delta,\varkappa$ are omitted because they will not be used. We will continue omitting such cases, often without further comment. Resulting tautologies will be also removed without further comment. By~\eqref{A_1}, applying Lemma~\ref{l: (m+1)^alpha (n+1)^beta (m-n+1)^gamma} to the above list, we get the f\/irst list of conditions that guarantee~\eqref{preliminary estimate} when $m\ge n$:
	\begin{gather}
		-\gamma,0\le\delta,\ -1<\varkappa\Rightarrow\alpha+\gamma,\alpha+\beta+\gamma+\delta+\varkappa+1<0,
		\label{A_1, -gamma,0 le delta, -1 < varkappa}\\
		0\le\delta<-\gamma,\ -1<\varkappa\Rightarrow\left\{
			\begin{array}{@{}l@{}}
				\alpha+\gamma+\delta,\alpha+\beta+\varkappa+1<0,\ \text{or}\\
				\alpha+\gamma,\alpha+\beta+\delta+\varkappa+1<0.
			\end{array}
		\right.\label{A_1, 0 le delta < -gamma, -1 < varkappa}
	\end{gather}

To prove~\eqref{A_1, (m-p+1)^gamma (n-p+1)^delta preccurlyeq ...}, it is enough to study the maximum of the $C^\infty$ function
	\begin{gather*}
		f(x)=(m-x+1)^\gamma(n-x+1)^\delta
	\end{gather*}
on $[0,n]$ (the natural domain of $f$ contains $(-\infty,n+1)$). We have
	\begin{gather*}
		f'(x)=(m-x+1)^{\gamma-1}(n-x+1)^{\delta-1}h(x),
	\end{gather*}
where
	\begin{gather*}
		h(x)=(\gamma+\delta)x-\gamma(n+1)-\delta(m+1).
	\end{gather*}
Observe that this expression is valid even when $\gamma=0$ or $\delta=0$. Since $f'$ and $h$ have the same zero set on $[0,n]$, and they have the same sign on the complement of the zero set in $[0,n]$, it is enough to analyze $h$ to know where $f$ reaches its maximum on $[0,n]$. We consider several cases.

{\em Case where $\gamma+\delta=0$.\/} Then $h\equiv\gamma(m-n)$. If $m>n$ and $\gamma\ne0$, then $h\ne0$ and $\operatorname{sign} h=\operatorname{sign} \gamma$. If $m=n$ or $\gamma=0$, then $h\equiv0$. Hence:
	\begin{gather}\label{A_1, max_0 le x le n f(x) if gamma + delta = 0}
		\max_{0\le x\le n}f(x)=
			\begin{cases}
				f(n)=(m-n+1)^\gamma & \text{if $\gamma=-\delta\ge0$},\\
				f(0)=(m+1)^\gamma(n+1)^\delta & \text{if $\gamma=-\delta\le0$}.
			\end{cases}
	\end{gather}
	
{\em Case where $\gamma+\delta\ne0$.\/} Then $h$ vanishes just at the point
	\begin{gather*}
		x_0:=\frac{\gamma(n+1)+\delta(m+1)}{\gamma+\delta}.
	\end{gather*}

{\em Case where $\gamma+\delta<0$.\/} We have $h>0$ on $(-\infty,x_0)$ and $h<0$ on $(x_0,\infty)$, yielding
	\begin{gather}\label{A_1, max_0 le x le n f(x) if gamma + delta < 0}
		\max_{0\le x\le n}f(x)=
			\begin{cases}
				f(0)=(m+1)^\gamma(n+1)^\delta & \text{if $x_0\le0$},\\
				f(x_0) & \text{if $0\le x_0\le n$},\\
				f(n)=(m-n+1)^\gamma & \text{if $x_0\ge n$}.
			\end{cases}
	\end{gather}

{\em Case where $\gamma+\delta<0$ and $\delta\le0$.\/} Then $x_0\ge n+1$, and therefore, by~\eqref{A_1, max_0 le x le n f(x) if gamma + delta < 0},
	\begin{gather}\label{A_1, max_0 le x le n f(x) if gamma + delta < 0 and delta le 0}
		\gamma+\delta<0,\ \delta\le0\Rightarrow
		\max_{0\le x\le n}f(x)=(m-n+1)^\gamma.
	\end{gather}

{\em Case where $\gamma+\delta<0$ and $\delta>0$; i.e., $0<\delta<-\gamma$.\/} We may have $x_0\le0$, $0\le x_0\le n$ or $n\le x_0$. Moreover
	\begin{gather*}
		f(x_0)=\frac{(-\gamma)^\gamma\delta^\delta}{(-\gamma-\delta)^{\gamma+\delta}}(m-n)^{\gamma+\delta}\preccurlyeq(m-n+1)^{\gamma+\delta}.
	\end{gather*}
Therefore
	\begin{gather}\label{1st estimate of max_0 le x le n f(x) if 0 < delta < -gamma}
		0<\delta<-\gamma\Rightarrow
		\max_{0\le x\le n}f(x)\preccurlyeq(m-n+1)^{\gamma+\delta}
	\end{gather}
 by~\eqref{A_1, max_0 le x le n f(x) if gamma + delta < 0} and since
	\begin{gather*}
		(m-n+1)^\gamma,(m+1)^\gamma(n+1)^\delta<(m-n+1)^{\gamma+\delta},
	\end{gather*}
which follows using that $\gamma<\gamma+\delta$ and
	\begin{gather*}
		n\ge\frac{m}{2}\Rightarrow\frac{m-n+1}{n+1}\le1\Rightarrow\frac{m-n+1}{m+1}<\frac{m-n+1}{n+1}
		\le\left(\frac{m-n+1}{n+1}\right)^{-\frac{\delta}{\gamma}},\\
		n<\frac{m}{2}\Rightarrow\frac{m-n+1}{n+1}>1\Rightarrow\frac{m-n+1}{m+1}\le1
		<\left(\frac{m-n+1}{n+1}\right)^{-\frac{\delta}{\gamma}},
	\end{gather*}
because $0<-\frac{\delta}{\gamma}<1$. On the other hand, in this case,
	\begin{gather*}
		\max_{0\le x\le n}f(x)\le\max_{0\le x\le n}(m-x+1)^\gamma\max_{0\le y\le n}(n-y+1)^\delta\\
		\hphantom{\max_{0\le x\le n}f(x)}{}
		=(m-n+1)^\gamma(n+1)^\delta<(m-n+1)^{\gamma+\delta}
	\end{gather*}
if $n<\frac{m}{2}$. So~\eqref{1st estimate of max_0 le x le n f(x) if 0 < delta < -gamma} can be improved by
	\begin{gather}\label{A_1, max_0 le x le n f(x) if 0 < delta < -gamma}
		0<\delta<-\gamma\Rightarrow
		\max_{0\le x\le n}f(x)\preccurlyeq
			\begin{cases}
				(m-n+1)^{\gamma+\delta}\ \text{and}\\
				(m-n+1)^\gamma(n+1)^\delta.
			\end{cases}
	\end{gather}
	
{\em Case where $\gamma+\delta>0$.\/} We have $h<0$ on $(-\infty,x_0)$ and $h>0$ on $(x_0,\infty)$, yielding
	\begin{gather}\label{max_0 le x le n f(x) if gamma + delta > 0}
		\max_{0\le x\le n}f(x)=
			\begin{cases}
				f(n)=(m-n+1)^\gamma & \text{if $x_0\le0$},\\
				\max\{f(0),f(n)\} & \text{if $0\le x_0\le n$},\\
				f(0)=(m+1)^\gamma(n+1)^\delta & \text{if $x_0\ge n$}.
			\end{cases}
	\end{gather}
	
{\em Case where $\gamma+\delta>0$ and $\delta\ge0$.\/} We get $x_0\ge n+1$, and therefore, by~\eqref{max_0 le x le n f(x) if gamma + delta > 0},
	\begin{gather}\label{A_1, max_0 le x le n f(x) if gamma + delta > 0 and delta ge 0}
		\gamma+\delta>0,\ \delta\ge0\Rightarrow
		\max_{0\le x\le n}f(x)=(m+1)^\gamma(n+1)^\delta.
	\end{gather}
		
Gathering together~\eqref{A_1, max_0 le x le n f(x) if gamma + delta = 0},~\eqref{A_1, max_0 le x le n f(x) if gamma + delta < 0 and delta le 0},~\eqref{A_1, max_0 le x le n f(x) if 0 < delta < -gamma} and~\eqref{A_1, max_0 le x le n f(x) if gamma + delta > 0 and delta ge 0}, we get the f\/irst two cases of~\eqref{A_1, (m-p+1)^gamma (n-p+1)^delta preccurlyeq ...}. The other cases will not be used, and they follow with similar arguments.

\subsubsection{Second list of conditions}\label{sss: A_2}

According to~\eqref{A_1, sum_p=0^n (p+1)^varkappa}, for all $\epsilon>0$,
	\begin{gather}\label{A_2, sum_p=0^n (n-p+1)^delta}
		\sum_{p=0}^n(n-p+1)^\delta=\sum_{q=1}^{n+1}q^\delta\preccurlyeq
			\begin{cases}
				(n+1)^{\delta+1} & \text{if $\delta>-1$},\\
				(n+1)^\epsilon & \text{if $\delta=-1$},\\
				1 & \text{if $\delta<-1$}.
			\end{cases}
	\end{gather}
On the other hand, we claim that
	\begin{gather}\label{A_2, (m-p+1)^gamma (p+1)^varkappa preccurlyeq ...}
		(m-p+1)^\gamma(p+1)^\varkappa\le
			\begin{cases}
				(m+1)^\gamma(n+1)^\varkappa & \text{if $\gamma,\varkappa\ge0$},\\
				(m-n+1)^\gamma(n+1)^\varkappa & \text{if $\gamma\le0\le\varkappa$},\\
				\left.\!
					\begin{array}{@{}l@{}}
						(m-n+1)^\gamma(n+1)^\varkappa\ \text{or}\\
						(m+1)^\gamma
					\end{array}
				\right\} &\text{if $\gamma<\varkappa<0$},\\
				(m+1)^\gamma & \text{if $\varkappa\le\gamma,0$,}
			\end{cases}
	\end{gather}
for all $p=0,\dots,n$. Combining~\eqref{A_2, sum_p=0^n (n-p+1)^delta} and~\eqref{A_2, (m-p+1)^gamma (p+1)^varkappa preccurlyeq ...}, it follows that, for all $\epsilon>0$,
	\begin{gather}\label{A_2}
		(m+1)^\alpha(n+1)^\beta\sum_{p=0}^n(m-p+1)^\gamma(n-p+1)^\delta(p+1)^\varkappa\preccurlyeq A_2,
	\end{gather}
where $A_2=A_2(m,n,\alpha,\beta,\gamma,\delta,\varkappa,\epsilon)$ can be taken to be equal to
	\begin{alignat*}{2}
		&(m+1)^{\alpha+\gamma}(n+1)^{\beta+\varkappa+\epsilon} &\quad&\text{if $0\le\gamma,\varkappa$, $\delta=-1$},\\
		&(m+1)^{\alpha+\gamma}(n+1)^{\beta+\varkappa} &\quad&\text{if $0\le\gamma,\varkappa$, $\delta<-1$},\\
		&(m+1)^\alpha(n+1)^{\beta+\varkappa+\epsilon}(m-n+1)^\gamma&\quad&\text{if $\gamma\le0\le\varkappa$, $\delta=-1$},\\
		&(m+1)^\alpha(n+1)^{\beta+\varkappa}(m-n+1)^\gamma&\quad&\text{if $\gamma\le0\le\varkappa$, $\delta<-1$},\\
		&\left.\!\!
			\begin{array}{@{}l@{}}
				(m+1)^\alpha(n+1)^{\beta+\varkappa+\epsilon}(m-n+1)^\gamma\ \text{or}\\[2pt]
				(m+1)^{\alpha+\gamma}(n+1)^{\beta+\epsilon}
			\end{array}
		\right\} &\quad&\text{if $\gamma<\varkappa<0$, $\delta=-1$},\\
		&\left.\!\!
			\begin{array}{@{}l@{}}
				(m+1)^\alpha(n+1)^{\beta+\varkappa}(m-n+1)^\gamma\ \text{or}\\[2pt]
				(m+1)^{\alpha+\gamma}(n+1)^\beta
			\end{array}
		\right\} &\quad&\text{if $\gamma<\varkappa<0$, $\delta<-1$},\\
		&(m+1)^{\alpha+\gamma}(n+1)^{\beta+\epsilon} &\quad&\text{if $\varkappa\le\gamma,0$, $\delta=-1$},\\
		&(m+1)^{\alpha+\gamma}(n+1)^\beta &\quad&\text{if $\varkappa\le\gamma,0$, $\delta<-1$}.
	\end{alignat*}
By~\eqref{A_2}, applying Lemma~\ref{l: (m+1)^alpha (n+1)^beta (m-n+1)^gamma} to the above list, we get the second list of conditions that guarantee~\eqref{preliminary estimate} when $m\ge n$:
	\begin{gather}
		0\le\gamma,\varkappa,\ \delta\le-1\Rightarrow\alpha+\gamma,\alpha+\beta+\gamma+\varkappa<0,
		\label{A_2, 0 le gamma,varkappa, delta le -1}\\
		\gamma\le\varkappa,0,\ \delta\le-1\Rightarrow\alpha+\gamma,\alpha+\beta+\varkappa<0,
		\label{A_2, gamma le varkappa,0, delta le -1}\\
		\varkappa\le\gamma,0,\ \delta\le-1\Rightarrow\alpha+\gamma,\alpha+\beta+\gamma<0.
		\label{A_2, varkappa le gamma,0, delta le -1}
	\end{gather}

To prove~\eqref{A_2, (m-p+1)^gamma (p+1)^varkappa preccurlyeq ...}, it is enough to study the maximum of the $C^\infty$ function
	\begin{gather*}
		f(x)=(m-x+1)^\gamma(x+1)^\varkappa
	\end{gather*}
on $[0,n]$ (the natural domain of $f$ contains $(-1,m+1)$). We have
	\begin{gather*}
		f'(x)=(m-x+1)^{\gamma-1}(x+1)^{\varkappa-1}h(x),
	\end{gather*}
where
	\begin{gather*}
		h(x)=-(\gamma+\varkappa)x+\varkappa(m+1)-\gamma.
	\end{gather*}
Observe that this expression is valid even when $\gamma=0$ or $\varkappa=0$. Since $f'$ and $h$ have the same zero set on $[0,n]$, and they have the same sign on the complement of the zero set in $[0,n]$, it is enough to analyze $h$ to know where $f$ reaches its maximum on $[0,n]$. We consider several cases.

{\em Case where $\gamma+\varkappa=0$.\/} Then $h\equiv\varkappa(m+2)$. If $\varkappa\ne0$, then $h\ne0$ and $\operatorname{sign} h=\operatorname{sign}\varkappa$. If $\varkappa=0$, then $h\equiv0$. Hence:
	\begin{gather}\label{A_2, max_0 le x le n f(x) if gamma + varkappa = 0}
		\max_{0\le x\le n}f(x)=
			\begin{cases}
				f(n)=(m-n+1)^\gamma(n+1)^\varkappa & \text{if $\varkappa=-\gamma\ge0$},\\
				f(0)=(m+1)^\gamma & \text{if $\varkappa=-\gamma\le0$}.
			\end{cases}
	\end{gather}
	
{\em Case where $\gamma+\varkappa\ne0$.\/} Then $h$ vanishes just at the point
	\begin{gather*}
		x_0:=\frac{\varkappa(m+1)-\gamma}{\gamma+\varkappa}.
	\end{gather*}

{\em Case where $\gamma+\varkappa<0$.\/} We have $h<0$ on $(-\infty,x_0)$ and $h>0$ on $(x_0,\infty)$, yielding
	\begin{gather}\label{A_2, max_0 le x le n f(x) if gamma + varkappa < 0}
		\max_{0\le x\le n}f(x)=
			\begin{cases}
				f(n)=(m-n+1)^\gamma(n+1)^{\varkappa} & \text{if $x_0\le0$},\\
				\max\{f(0),f(n)\} & \text{if $0\le x_0\le n$},\\
				f(0)=(m+1)^\gamma & \text{if $x_0\ge n$}.
			\end{cases}
	\end{gather}

{\em Case where $\gamma+\varkappa<0$ and $\varkappa\ge0$; i.e., $0\le\varkappa<-\gamma$.\/} Then $x_0<0$, and therefore, by~\eqref{A_2, max_0 le x le n f(x) if gamma + varkappa < 0},
	\begin{gather}\label{A_2, max_0 le x le n f(x) if 0 le varkappa < -gamma}
		0\le\varkappa<-\gamma\Rightarrow\max_{0\le x\le n}f(x)=(m-n+1)^\gamma(n+1)^\varkappa.
	\end{gather}
	
{\em Case where $\gamma+\varkappa<0$ and $\gamma\ge0$; i.e., $0\le\gamma<-\varkappa$.\/} Then $x_0\ge m+1$, and therefore, by~\eqref{A_2, max_0 le x le n f(x) if gamma + varkappa < 0},
	\begin{gather}\label{A_2, max_0 le x le n f(x) if 0 le gamma < -varkappa}
		0\le\gamma<-\varkappa\Rightarrow\max_{0\le x\le n}f(x)=(m+1)^\gamma.
	\end{gather}
	
{\em Case where $\varkappa\le\gamma<0$.\/} Then $x_0\ge\frac{m}{2}$, and we may have $x_0\le n$ or $x_0\ge n$. In any case, by~\eqref{A_2, max_0 le x le n f(x) if gamma + varkappa < 0},
	\begin{gather}\label{A_2, max_0 le x le n f(x) if varkappa le gamma < 0}
		\varkappa\le\gamma<0\Rightarrow\max_{0\le x\le n}f(x)=(m+1)^\gamma.
	\end{gather}
	
{\em Case where $\gamma<\varkappa<0$.\/} Then $x_0<\frac{m}{2}$, and we may have $x_0\le0$, $0\le x_0\le n$ or $x_0\ge n$. In any case, by~\eqref{A_2, max_0 le x le n f(x) if gamma + varkappa < 0},
	\begin{gather}\label{A_2, max_0 le x le n f(x) if gamma < varkappa < 0}
		\gamma<\varkappa<0\Rightarrow\max_{0\le x\le n}f(x)=
			\begin{cases}
				(m-n+1)^\gamma(n+1)^\varkappa\ \text{or}\\
				(m+1)^\gamma .
			\end{cases}
	\end{gather}
	
{\em Case where $\gamma+\varkappa>0$.\/} We have $h>0$ on $(-\infty,x_0)$ and $h<0$ on $(x_0,\infty)$, yielding
	\begin{gather}\label{A_2, max_0 le x le n f(x) if gamma + varkappa > 0}
		\max_{0\le x\le n}f(x)=
			\begin{cases}
				f(0)=(m+1)^\gamma & \text{if $x_0\le0$},\\
				f(x_0) & \text{if $0\le x_0\le n$},\\
				f(n)=(m-n+1)^\gamma(n+1)^\varkappa & \text{if $x_0\ge n$}.
			\end{cases}
	\end{gather}
	
{\em Case where $\gamma+\varkappa>0$ and $\varkappa\le0$; i.e., $-\gamma<\varkappa\le0$.\/} Then $x_0<0$, and therefore, by~\eqref{A_2, max_0 le x le n f(x) if gamma + varkappa > 0},
	\begin{gather}\label{A_2, max_0 le x le n f(x) if -gamma < varkappa le 0}
		-\gamma<\varkappa\le0\Rightarrow\max_{0\le x\le n}f(x)=(m+1)^\gamma.
	\end{gather}
	
{\em Case where $\gamma+\varkappa>0$ and $\gamma\le0$; i.e., $-\varkappa<\gamma\le0$.\/} Then $x_0\ge m+1$, and therefore, by~\eqref{A_2, max_0 le x le n f(x) if gamma + varkappa > 0},
	\begin{gather}\label{A_2, max_0 le x le n f(x) if -varkappa < gamma le 0}
		-\varkappa<\gamma\le0\Rightarrow\max_{0\le x\le n}f(x)=(m-n+1)^\gamma(n+1)^\varkappa.
	\end{gather}
	
{\em Case where $\gamma + \varkappa>0$ and $\gamma, \varkappa\ge 0$.\/} We may have $x_0\le 0$, $0\le x_0\le n$ or $n\le x_0$. Moreover
	\begin{gather*}
		f(x_0)=\frac{\gamma^\gamma\varkappa^\varkappa}{(\gamma+\varkappa)^{\gamma+\varkappa}}(m+2)^{\gamma+\varkappa}\preccurlyeq(m+1)^{\gamma+\varkappa}.
	\end{gather*}
But, in this case,
	\begin{gather*}
		\max_{0\le x\le n}f(x)\le\max_{0\le x\le n}(m-x+1)^\gamma\max_{0\le y\le n}(y+1)^\varkappa
		=(m+1)^\gamma(n+1)^\varkappa\le(m+1)^{\gamma+\varkappa}.
	\end{gather*}	
Therefore, by~\eqref{A_2, max_0 le x le n f(x) if gamma + varkappa > 0},
	\begin{gather}\label{A_2, max_0 le x le n f(x) if gamma,varkappa ge 0}
		\gamma,\varkappa\ge0\Rightarrow\max_{0\le x\le n}f(x)\le(m+1)^\gamma(n+1)^\varkappa.
	\end{gather}

Gathering together~\eqref{A_2, max_0 le x le n f(x) if gamma + varkappa = 0},~\eqref{A_2, max_0 le x le n f(x) if 0 le varkappa < -gamma}--\eqref{A_2, max_0 le x le n f(x) if gamma < varkappa < 0} and~\eqref{A_2, max_0 le x le n f(x) if -gamma < varkappa le 0}--\eqref{A_2, max_0 le x le n f(x) if gamma,varkappa ge 0}, we get~\eqref{A_2, (m-p+1)^gamma (p+1)^varkappa preccurlyeq ...}.
	
\subsubsection{Third list of conditions}\label{sss: A_3}

For all $\epsilon>0$,
		\begin{gather}
			\sum_{p=0}^n(m-p+1)^\gamma=\sum_{q=m-n+1}^{m+1}q^\gamma\le
				\begin{cases}
					\displaystyle \int_{m-n+1}^{m+2}x^\gamma dx & \text{if $\gamma\ge0$},\\
					\displaystyle (m-n+1)^\gamma+\int_{m-n+1}^{m+1}x^\gamma dx & \text{if $\gamma<0$}
				\end{cases}\notag\\
			\hphantom{\sum_{p=0}^n(m-p+1)^\gamma}{}
			\preccurlyeq
				\begin{cases}
					(m+1)^{\gamma+1} & \text{if $\gamma>-1$},\\
					1+\ln(m+1) & \text{if $\gamma=-1$},\\
					(m-n+1)^{\gamma+1} & \text{if $\gamma<-1$}
				\end{cases}
			\preccurlyeq
				\begin{cases}
					(m+1)^{\gamma+1} & \text{if $\gamma>-1$},\\
					(m+1)^\epsilon & \text{if $\gamma=-1$},\\
					(m-n+1)^{\gamma+1} & \text{if $\gamma<-1$}.
				\end{cases}\label{A_3, sum_p=0^n (m-p+1)^gamma, with int}
		\end{gather}
	The following gives a better estimate when $\gamma\ge0$, and an alternative estimate when $\gamma<0$:
		\begin{gather}\label{A_3, sum_p=0^n (m-p+1)^gamma, without int}
			\sum_{p=0}^n(m-p+1)^\gamma\le
				\begin{cases}
					(m+1)^\gamma(n+1) & \text{if $\gamma\ge0$},\\
					(m-n+1)^\gamma(n+1) & \text{if $\gamma<0$}.
				\end{cases}
		\end{gather}
The \looseness=-1 estimate~\eqref{A_3, sum_p=0^n (m-p+1)^gamma, without int} is better than~\eqref{A_3, sum_p=0^n (m-p+1)^gamma, with int} when $\gamma\ge0$ since $m\ge n$, and~\eqref{A_3, sum_p=0^n (m-p+1)^gamma, without int} may be better or worse than~\eqref{A_3, sum_p=0^n (m-p+1)^gamma, with int} when $\gamma<0$, depending on the values of $m$ and $n$. Note that estimates of the type~\eqref{A_3, sum_p=0^n (m-p+1)^gamma, without int} for $\sum\limits_{p=0}^n(p+1)^\varkappa$ and $\sum\limits_{p=0}^n(n-p+1)^\delta$ are worse than~\eqref{A_1, sum_p=0^n (p+1)^varkappa} and~\eqref{A_2, sum_p=0^n (n-p+1)^delta}. Thus it makes no sense to add this kind of estimate in Sections~\ref{sss: A_1} and~\ref{sss: A_2}. On the other hand, we claim that
	\begin{gather}\label{A_3, (n-p+1)^delta (p+1)^varkappa preccurlyeq ...}
		(n-p+1)^\delta(p+1)^\varkappa\le
			\begin{cases}
				(n+1)^\varkappa & \text{if $\delta\le\varkappa,0$},\\
				(n+1)^{\delta+\varkappa} & \text{if $\delta,\varkappa\ge0$},\\
				(n+1)^\delta & \text{if $\varkappa\le\delta,0$}
			\end{cases}
	\end{gather}
for all $p=0,\dots,n$. Combining~\eqref{A_3, sum_p=0^n (m-p+1)^gamma, with int}--\eqref{A_3, (n-p+1)^delta (p+1)^varkappa preccurlyeq ...}, it follows that, for all $\epsilon>0$,
	\begin{gather}\label{A_3}
		(m+1)^\alpha(n+1)^\beta\sum_{p=0}^n(m-p+1)^\gamma(n-p+1)^\delta(p+1)^\varkappa\preccurlyeq A_3,
	\end{gather}
where $A_3=A_3(m,n,\alpha,\beta,\gamma,\delta,\varkappa,\epsilon)$ can be taken to be equal to
	\begin{alignat*}{2}
		&\left.\!\!
			\begin{array}{@{}l@{}}
				(m+1)^{\alpha+\epsilon}(n+1)^{\beta+\delta+\varkappa}\ \text{and}\\[2pt]
				(m+1)^\alpha(n+1)^{\beta+\delta+\varkappa+1}(m-n+1)^\gamma
			\end{array}
		\right\} &\quad&\text{if $\gamma=-1$, $0\le\delta,\varkappa$},\\
		&\left. \!\!
			\begin{array}{@{}l@{}}
				(m+1)^\alpha(n+1)^{\beta+\delta+\varkappa}(m-n+1)^{\gamma+1}\ \text{and}\\[2pt]
				(m+1)^\alpha(n+1)^{\beta+\delta+\varkappa+1}(m-n+1)^\gamma
			\end{array}
		\right\} &\quad&\text{if $\gamma<-1$, $0\le\delta,\varkappa$}.
	\end{alignat*}
By~\eqref{A_3}, applying Lemma~\ref{l: (m+1)^alpha (n+1)^beta (m-n+1)^gamma} to the above list, we get the third list of conditions that guarantee~\eqref{preliminary estimate} when $m\ge n$:
	\begin{gather}\label{A_3, gamma le -1, 0 le delta,varkappa}
		\gamma\le-1,\ 0\le\delta,\varkappa\Rightarrow
			\begin{cases}
				\alpha+\gamma+1,\alpha+\beta+\delta+\varkappa<0,\ \text{or}\\
				\alpha+\gamma,\alpha+\beta+\delta+\varkappa+1<0.
			\end{cases}
	\end{gather}
	
To prove~\eqref{A_3, (n-p+1)^delta (p+1)^varkappa preccurlyeq ...}, it is enough to study the maximum of the $C^\infty$ function
	\begin{gather*}
		f(x)=(n-x+1)^\delta(x+1)^\varkappa
	\end{gather*}
on $[0,n]$ (the natural domain of $f$ contains $(-1,n+1)$). We have
	\begin{gather*}
		f'(x)=(n-x+1)^{\delta-1}(x+1)^{\varkappa-1}h(x),
	\end{gather*}
where
	\begin{gather*}
		h(x)=-(\delta+\varkappa)x+\varkappa(n+1)-\delta.
	\end{gather*}
Observe that this expression is valid even when $\delta=0$ or $\varkappa=0$. Since $f'$ and $h$ have the same zero set on $[0,n]$, and they have the same sign on the complement of the zero set in $[0,n]$, it is enough to analyze $h$ to know where $f$ reaches its maximum on $[0,n]$. We consider several cases.

{\em Case where $\delta+\varkappa=0$.\/} Then $h\equiv\varkappa(n+2)$. If $\varkappa\ne0$, then $h\ne0$ and $\operatorname{sign} h=\operatorname{sign}\varkappa$. If $\varkappa=0$, then $h\equiv0$. Hence:
	\begin{gather}\label{A_3, max_0 le x le n f(x) if delta + varkappa = 0}
		\max_{0\le x\le n}f(x)=
			\begin{cases}
				f(n)=(n+1)^\varkappa & \text{if $\varkappa=-\delta\ge0$},\\
				f(0)=(n+1)^\delta & \text{if $\varkappa=-\delta\le0$}.
			\end{cases}
	\end{gather}
	
{\em Case where $\delta+\varkappa\ne0$.\/} Then $h$ vanishes just at the point
	\begin{gather*}
		x_0:=\frac{\varkappa(n+1)-\delta}{\delta+\varkappa}.
	\end{gather*}

{\em Case where $\delta+\varkappa>0$.\/} We have $h>0$ on $(-\infty,x_0)$ and $h<0$ on $(x_0,\infty)$, yielding
	\begin{gather}\label{A_3, max_0 le x le n f(x) if delta + varkappa > 0}
		\max_{0\le x\le n}f(x)=
			\begin{cases}
				f(0)=(n+1)^\delta & \text{if $x_0\le0$},\\
				f(x_0) & \text{if $0\le x_0\le n$},\\
				f(n)=(n+1)^\varkappa & \text{if $x_0\ge n$}.
			\end{cases}
	\end{gather}
	
{\em Case where $\delta+\varkappa>0$ and $\delta,\varkappa\ge0$.\/} We may have $x_0\le0$, $0\le x_0\le n$ or $n\le x_0$. Moreover
	\begin{gather*}
		f(x_0)=\frac{\delta^\delta\varkappa^\varkappa}{(\delta+\varkappa)^{\delta+\varkappa}}(n+2)^{\delta+\varkappa}\preccurlyeq(n+1)^{\delta+\varkappa}.
	\end{gather*}
But, in this case,
	\begin{gather*}
		\max_{0\le x\le n}f(x)\le\max_{0\le x\le n}(n-x+1)^\delta\max_{0\le y\le n}(y+1)^\varkappa=(n+1)^{\delta+\varkappa}.
	\end{gather*}	
Therefore, by~\eqref{A_3, max_0 le x le n f(x) if delta + varkappa > 0},
	\begin{gather}\label{A_3, max_0 le x le n f(x) if delta,varkappa ge 0}
		\delta,\varkappa\ge0\Rightarrow\max_{0\le x\le n}f(x)\le(n+1)^{\delta+\varkappa}.
	\end{gather}

Gathering together~\eqref{A_3, max_0 le x le n f(x) if delta + varkappa = 0} and~\eqref{A_3, max_0 le x le n f(x) if delta,varkappa ge 0}, we get the second case of~\eqref{A_3, (n-p+1)^delta (p+1)^varkappa preccurlyeq ...}. The other cases will not be used, and they follow with similar arguments.

\subsubsection{Fourth list of conditions}\label{sss: A_4}

We have
	\begin{gather}\label{A_4, (p+1)^varkappa}
		(p+1)^\varkappa\le
			\begin{cases}
				(n+1)^\varkappa & \text{if $\varkappa\ge0$},\\
				1 & \text{if $\varkappa\le0$}
			\end{cases}
	\end{gather}
for $p=0,\dots,n$. Moreover, by~\eqref{A_1, sum_p=0^n (p+1)^varkappa},~\eqref{A_3, sum_p=0^n (m-p+1)^gamma, with int} and~\eqref{A_3, sum_p=0^n (m-p+1)^gamma, without int}, for all $\epsilon>0$,
	\begin{gather}
		\sum_{p=0}^n(n-p+1)^{2\delta}=\sum_{q=1}^{n+1}q^{2\delta}\preccurlyeq
			\begin{cases}
				(n+1)^{2\delta+1} & \text{if $\delta>-\frac{1}{2}$},\\
				(n+1)^\epsilon & \text{if $\delta=-\frac{1}{2}$},\\
				1 & \text{if $\delta<-\frac{1}{2}$},
			\end{cases}
		\label{A_4, sum_p=0^n (n-p+1)^2 delta, with int}\\
		\sum_{p=0}^n(m-p+1)^{2\gamma}\preccurlyeq
			\begin{cases}
				(m+1)^{2\gamma+1} & \text{if $\gamma>-\frac{1}{2}$},\\
				(m+1)^\epsilon & \text{if $\gamma=-\frac{1}{2}$},\\
				(m-n+1)^{2\gamma+1} & \text{if $\gamma<-\frac{1}{2}$},
			\end{cases}
		\label{A_4, sum_p=0^n (m-p+1)^2 gamma, with int}\\
		\sum_{p=0}^n(m-p+1)^{2\gamma}\le
			\begin{cases}
				(m+1)^{2\gamma}(n+1) & \text{if $\gamma\ge0$},\\
				(m-n+1)^{2\gamma}(n+1) & \text{if $\gamma<0$}.
			\end{cases}
		\label{A_4, sum_p=0^n (m-p+1)^2 gamma, without int}
	\end{gather}
The estimate~\eqref{A_4, sum_p=0^n (m-p+1)^2 gamma, without int} is better than~\eqref{A_4, sum_p=0^n (m-p+1)^2 gamma, with int} when $\gamma\ge0$, and it may be better or worse than~\eqref{A_4, sum_p=0^n (m-p+1)^2 gamma, with int} when $\gamma<0$, depending on the values of $m$ and $n$. By the Cauchy--Schwartz inequality,
	\begin{gather*}
		\sum_{p=0}^n(m-p+1)^\gamma(n-p+1)^\delta\le
		\left(\sum_{p=0}^n(m-p+1)^{2\gamma}\right)^{\frac{1}{2}}
		\left(\sum_{p=0}^n(n-p+1)^{2\delta}\right)^{\frac{1}{2}}.
	\end{gather*}
Therefore, by~\eqref{A_4, (p+1)^varkappa}--\eqref{A_4, sum_p=0^n (m-p+1)^2 gamma, without int},
	\begin{gather}\label{A_4}
		(m+1)^\alpha(n+1)^\beta\sum_{p=0}^n(m-p+1)^\gamma(n-p+1)^\delta(p+1)^\varkappa\preccurlyeq A_4,
	\end{gather}
where $A_4=A_4(m,n,\alpha,\beta,\gamma,\delta,\varkappa)$ can be taken to be equal to
	\begin{gather*}
		\left.\!\!
			\begin{array}{@{}l@{}}
				(m+1)^\alpha(n+1)^{\beta+\delta+\varkappa+\frac{1}{2}}(m-n+1)^{\gamma+\frac{1}{2}}\ \text{and}\\[2pt]
				(m+1)^\alpha(n+1)^{\beta+\delta+\varkappa+1}(m-n+1)^\gamma
			\end{array}
		\right\}\quad\text{if $\gamma<-\tfrac{1}{2}<\delta$, $0\le\varkappa$}.
	\end{gather*}
By~\eqref{A_4}, applying Lemma~\ref{l: (m+1)^alpha (n+1)^beta (m-n+1)^gamma} to the above list, we get the fourth list of conditions that guarantee~\eqref{preliminary estimate} when $m\ge n$:
	\begin{gather}\label{A_4, gamma < -1/2 < delta, 0 le varkappa}
		\gamma<- \tfrac{1}{2}<\delta,\ 0\le\varkappa
		\Rightarrow\left\{\negmedspace\negthickspace
			\begin{array}{l}
				\alpha+\gamma+\frac{1}{2},\alpha+\beta+\delta+\varkappa+\frac{1}{2}<0,\ \text{or}\\[2pt]
				\alpha+\gamma,\alpha+\beta+\delta+\varkappa+1<0.
			\end{array}
		\right.
	\end{gather}

\subsubsection{Fifth list of conditions}\label{sss: A_5}

This is analogous to the estimates of Section~\ref{sss: A_4}, interchanging the roles of $\delta$ and $\varkappa$. We have
	\begin{gather}\label{A_5, (n-p+1)^delta}
		(n-p+1)^\delta\le
			\begin{cases}
				(n+1)^\delta & \text{if $\delta\ge0$},\\
				1 & \text{if $\delta\le0$}
			\end{cases}
	\end{gather}
for $p=0,\dots,n$. Moreover, by~\eqref{A_1, sum_p=0^n (p+1)^varkappa}, for all $\epsilon>0$,
	\begin{gather}\label{A_5, sum_p=0^n (p+1)^2 varkappa, with int}
		\sum_{p=0}^n(p+1)^{2\varkappa}\preccurlyeq
			\begin{cases}
				(n+1)^{2\varkappa+1} & \text{if $\varkappa>-\frac{1}{2}$},\\
				(n+1)^\epsilon & \text{if $\varkappa=-\frac{1}{2}$},\\
				1 & \text{if $\varkappa<-\frac{1}{2}$}.
			\end{cases}
	\end{gather}
Applying the Cauchy--Schwartz inequality, we get
	\begin{gather*}
		\sum_{p=0}^n(m-p+1)^\gamma(p+1)^\varkappa\le
		\left(\sum_{p=0}^n(m-p+1)^{2\gamma}\right)^{\frac{1}{2}}
		\left(\sum_{p=0}^n(p+1)^{2\varkappa}\right)^{\frac{1}{2}}.
	\end{gather*}
Therefore, by~\eqref{A_5, (n-p+1)^delta},~\eqref{A_5, sum_p=0^n (p+1)^2 varkappa, with int},~\eqref{A_4, sum_p=0^n (m-p+1)^2 gamma, with int} and~\eqref{A_4, sum_p=0^n (m-p+1)^2 gamma, without int}, for all $\epsilon>0$,
	\begin{gather}\label{A_5}
		(m+1)^\alpha(n+1)^\beta\sum_{p=0}^n(m-p+1)^\gamma(n-p+1)^\delta(p+1)^\varkappa\preccurlyeq A_5,
	\end{gather}
where $A_5=A_5(m,n,\alpha,\beta,\gamma,\delta,\varkappa,\epsilon)$ can be taken to be equal to
	\begin{alignat*}{2}
		&\left.\!\!
			\begin{array}{@{}l@{}}
				(m+1)^{\alpha+\epsilon}(n+1)^{\beta+\delta+\epsilon}\ \text{and}\\[2pt]
				(m+1)^\alpha(n+1)^{\beta+\delta+\frac{1}{2}+\epsilon}(m-n+1)^\gamma
			\end{array}
		\right\} &\quad&\text{if $\gamma=\varkappa=- \tfrac{1}{2}$, $0\le\delta$},\\
		&\left.\!\!
			\begin{array}{@{}l@{}}
				(m+1)^{\alpha+\epsilon}(n+1)^{\beta+\delta}\ \text{and}\\[2pt]
				(m+1)^\alpha(n+1)^{\beta+\delta+\frac{1}{2}}(m-n+1)^\gamma
			\end{array}
		\right\} &\quad&\text{if $\varkappa<\gamma=- \tfrac{1}{2}$, $0\le\delta$},\\
		&\left.\!\!
			\begin{array}{@{}l@{}}
				(m+1)^\alpha(n+1)^{\beta+\delta+\epsilon}(m-n+1)^{\gamma+\frac{1}{2}}\ \text{and}\\[2pt]
				(m+1)^\alpha(n+1)^{\beta+\delta+\frac{1}{2}+\epsilon}(m-n+1)^\gamma
			\end{array}
		\right\} &\quad&\text{if $\gamma<\varkappa=- \tfrac{1}{2}$, $0\le\delta$},\\
		&\left.\!\!
			\begin{array}{@{}l@{}}
				(m+1)^\alpha(n+1)^{\beta+\delta}(m-n+1)^{\gamma+\frac{1}{2}}\ \text{and}\\[2pt]
				(m+1)^\alpha(n+1)^{\beta+\delta+\frac{1}{2}}(m-n+1)^\gamma
			\end{array}
		\right\} &\quad&\text{if $\gamma,\varkappa<- \tfrac{1}{2}$, $0\le\delta$}.
	\end{alignat*}
By~\eqref{A_5}, applying Lemma~\ref{l: (m+1)^alpha (n+1)^beta (m-n+1)^gamma} to the above list, we get the f\/ifth list of conditions that guarantee~\eqref{preliminary estimate} when $m\ge n$:
	\begin{gather}\label{A_5, gamma,varkappa le -1/2, 0 le delta}
		\gamma,\varkappa\le- \tfrac{1}{2},\ 0\le\delta
		\Rightarrow
			\begin{cases}
				\alpha+\gamma+\frac{1}{2},\alpha+\beta+\delta<0,\ \text{or}\\[2pt]
				\alpha+\gamma,\alpha+\beta+\delta+\frac{1}{2}<0.
			\end{cases}
	\end{gather}

\subsubsection{Sixth list of conditions}\label{sss: A_6}

We have
	\begin{gather}\label{A_6, (m-p+1)^gamma}
		(m-p+1)^\gamma\le
			\begin{cases}
				(m+1)^\gamma & \text{if $\gamma\ge0$},\\
				(m-n+1)^\gamma & \text{if $\gamma\le0$}
			\end{cases}
	\end{gather}
for $p=0,\dots,n$. Moreover, by~\eqref{A_1, sum_p=0^n (p+1)^varkappa}, for all $\epsilon>0$,
	\begin{gather}\label{A_6, sum_p=0^n (n-p+1)^2 delta, with int}
		\sum_{p=0}^n(n-p+1)^{2\delta}=\sum_{q=1}^{n+1}q^{2\delta}\preccurlyeq
			\begin{cases}
				(n+1)^{2\delta+1} & \text{if $\delta>-\frac{1}{2}$},\\
				(n+1)^\epsilon & \text{if $\delta=-\frac{1}{2}$},\\
				1 & \text{if $\delta<-\frac{1}{2}$}.
			\end{cases}
	\end{gather}
Applying the Cauchy--Schwartz inequality, we get
	\begin{gather*}
		\sum_{p=0}^n(n-p+1)^\delta(p+1)^\varkappa\le
		\left(\sum_{p=0}^n(n-p+1)^{2\delta}\right)^{\frac{1}{2}}
		\left(\sum_{p=0}^n(p+1)^{2\varkappa}\right)^{\frac{1}{2}}.
	\end{gather*}
Therefore, by~\eqref{A_5, sum_p=0^n (p+1)^2 varkappa, with int},~\eqref{A_6, (m-p+1)^gamma} and~\eqref{A_6, sum_p=0^n (n-p+1)^2 delta, with int}, for all $\epsilon>0$,
	\begin{gather}\label{A_6}
		(m+1)^\alpha(n+1)^\beta\sum_{p=0}^n(m-p+1)^\gamma(n-p+1)^\delta(p+1)^\varkappa\preccurlyeq A_6,
	\end{gather}
where $A_6=A_6(m,n,\alpha,\beta,\gamma,\delta,\varkappa,\epsilon)$ can be taken to be equal to
	\begin{alignat*}{2}
		&(m+1)^{\alpha+\gamma}(n+1)^{\beta+\epsilon}
		&\quad&\text{if}\
			\begin{cases}
				\varkappa\le\delta=- \tfrac{1}{2},\ 0\le\gamma,\ \text{or}\\[2pt]
				\delta\le\varkappa=- \tfrac{1}{2},\ 0\le\gamma,
			\end{cases}\\
		&(m+1)^{\alpha+\gamma}(n+1)^\beta
		&\quad&\text{if $\delta,\varkappa<- \tfrac{1}{2}$, $0\le\gamma$}.
	\end{alignat*}
By~\eqref{A_6}, applying Lemma~\ref{l: (m+1)^alpha (n+1)^beta (m-n+1)^gamma} to the above list, we get the sixth list of conditions that guarantee~\eqref{preliminary estimate} when $m\ge n$:
	\begin{gather}\label{A_6, delta,varkappa le -1/2, 0 le gamma}
		\delta,\varkappa\le- \tfrac{1}{2},\ 0\le\gamma \Rightarrow\alpha+\gamma,\alpha+\beta+\gamma<0.
	\end{gather}

\subsubsection{Seventh list of conditions}\label{sss: A_7}

By~\eqref{A_1, sum_p=0^n (p+1)^varkappa},~\eqref{A_3, sum_p=0^n (m-p+1)^gamma, with int} and~\eqref{A_3, sum_p=0^n (m-p+1)^gamma, without int}, for all $\epsilon>0$,
	\begin{gather}
		\sum_{p=0}^n(n-p+1)^{3\delta}=\sum_{q=1}^{n+1}q^{3\delta}\preccurlyeq
			\begin{cases}
				(n+1)^{3\delta+1} & \text{if $\delta>-\frac{1}{3}$},\\
				(n+1)^\epsilon & \text{if $\delta=-\frac{1}{3}$},\\
				1 & \text{if $\delta<-\frac{1}{3}$},
			\end{cases}\label{A_7, sum_p=0^n (n-p+1)^3 delta, with int}\\
		\sum_{p=0}^n(p+1)^{3\varkappa}\preccurlyeq
			\begin{cases}
				(n+1)^{3\varkappa+1} & \text{if $\varkappa>-\frac{1}{3}$},\\
				(n+1)^\epsilon & \text{if $\varkappa=-\frac{1}{3}$},\\
				1 & \text{if $\varkappa<-\frac{1}{3}$},
			\end{cases}\label{A_7, sum_p=0^n (p+1)^3 varkappa, with int}\\
		\sum_{p=0}^n(m-p+1)^{3\gamma}\preccurlyeq
			\begin{cases}
				(m+1)^{3\gamma+1} & \text{if $\gamma>-\frac{1}{3}$},\\
				(m+1)^\epsilon & \text{if $\gamma=-\frac{1}{3}$},\\
				(m-n+1)^{3\gamma+1} & \text{if $\gamma<-\frac{1}{3}$},
			\end{cases}\label{A_7, sum_p=0^n (m-p+1)^3 gamma, with int}\\
		\sum_{p=0}^n(m-p+1)^{3\gamma}\le
			\begin{cases}
				(m+1)^{3\gamma}(n+1) & \text{if $\gamma\ge0$},\\
				(m-n+1)^{3\gamma}(n+1) & \text{if $\gamma<0$}.
			\end{cases}\label{A_7, sum_p=0^n (m-p+1)^3 gamma, without int}
	\end{gather}
Note that~\eqref{A_7, sum_p=0^n (m-p+1)^3 gamma, without int} is better than~\eqref{A_7, sum_p=0^n (m-p+1)^3 gamma, with int} for $\gamma\ge0$, and it is an alternative estimate for $\gamma<0$. Applying the generalized H\"older inequality \cite{Cheung2001}, we get
	\begin{gather*}
		\sum_{p=0}^n(m-p+1)^\gamma(n-p+1)^\delta(p+1)^\varkappa\\
		\qquad{}\leq \left(\sum_{p=0}^n(m-p+1)^{3\gamma}\right)^{\frac{1}{3}}
		\left(\sum_{p=0}^n(n-p+1)^{3\delta}\right)^{\frac{1}{3}}
		\left(\sum_{p=0}^n(p+1)^{3\varkappa}\right)^{\frac{1}{3}}.
	\end{gather*}
Therefore, by~\eqref{A_7, sum_p=0^n (n-p+1)^3 delta, with int}--\eqref{A_7, sum_p=0^n (m-p+1)^3 gamma, without int},
	\begin{gather}\label{A_7}
		(m+1)^\alpha(n+1)^\beta\sum_{p=0}^n(m-p+1)^\gamma(n-p+1)^\delta(p+1)^\varkappa\preccurlyeq A_7,
	\end{gather}
where $A_7=A_7(m,n,\alpha,\beta,\gamma,\delta,\varkappa)$ can be taken to be equal to
	\begin{gather*}
		\left.\!\!
			\begin{array}{@{}l@{}}
				(m+1)^\alpha(n+1)^{\beta+\delta+\varkappa+\frac{2}{3}}(m-n+1)^{\gamma+\frac{1}{3}}\ \text{and}\\[2pt]
				(m+1)^\alpha(n+1)^{\beta+\delta+\varkappa+1}(m-n+1)^\gamma
			\end{array}\right\}
		 \quad \text{if $\gamma<- \tfrac{1}{3}<\delta,\varkappa$},\\
		 \left.\!\!
			\begin{array}{@{}l@{}}
				(m+1)^\alpha(n+1)^{\beta+\delta+\frac{1}{3}}(m-n+1)^{\gamma+\frac{1}{3}}\ \text{and}\\[2pt]
				(m+1)^\alpha(n+1)^{\beta+\delta+\frac{2}{3}}(m-n+1)^\gamma
			\end{array}\right\}
		 \quad \text{if $\gamma,\varkappa<- \tfrac{1}{3}<\delta$}.
	\end{gather*}
By~\eqref{A_7}, applying Lemma~\ref{l: (m+1)^alpha (n+1)^beta (m-n+1)^gamma} to the above list, we get the seventh list of conditions that guarantee~\eqref{preliminary estimate} when $m\ge n$:
	\begin{gather}
		\gamma<- \tfrac{1}{3}<\delta,\varkappa
		\Rightarrow
			\begin{cases}
				\alpha+\gamma+\frac{1}{3},\alpha+\beta+\delta+\varkappa+ \tfrac{2}{3}<0,\ \text{or}\\[2pt]
				\alpha+\gamma,\alpha+\beta+\delta+\varkappa+1<0,
			\end{cases}\label{A_7, gamma < -1/3 < delta,varkappa}\\
		\gamma,\varkappa<- \tfrac{1}{3}<\delta
		\Rightarrow
			\begin{cases}
				\alpha+\gamma+\frac{1}{3},\alpha+\beta+\delta+ \tfrac{1}{3}<0,\ \text{or}\\[2pt]
				\alpha+\gamma,\alpha+\beta+\delta+ \tfrac{2}{3}<0.
			\end{cases}\label{A_7, gamma,varkappa < -1/3 < delta}
	\end{gather}

\subsubsection[Obtaining the sets ${\mathfrak{S}}_{ijk}$ from the lists of conditions]{Obtaining the sets $\boldsymbol{{\mathfrak{S}}_{ijk}}$ from the lists of conditions}\label{sss: conditions}

The left hand side of the conditions from the lists of Sections~\ref{sss: A_1}--\ref{sss: A_7} involve only $(\gamma,\delta,\varkappa)$. Now, we indicate which of them def\/ine sets covering every $Q_{ijk}$, for the chosen subindices $ijk$ equal to 515, 522, 252, 155, 212. Those conditions will produce the def\/inition of ${\mathfrak{S}}_{ijk}\subset{\mathbb{R}}^2\times Q_{ijk}$ so that~\eqref{preliminary estimate} holds for $m\ge n$.

The set $Q_{515}$ is given by the left hand side of~\eqref{A_2, 0 le gamma,varkappa, delta le -1}, whose right hand side is~\eqref{alpha + gamma, alpha + beta + gamma + varkappa < 0}, def\/i\-ning~${\mathfrak{S}}_{515}$.

Any $(\gamma,\delta,\varkappa)\in Q_{522}$ satisf\/ies the left hand side of~\eqref{A_6, delta,varkappa le -1/2, 0 le gamma}, whose right hand side is~\eqref{alpha + gamma, alpha + beta + gamma < 0}, def\/i\-ning~${\mathfrak{S}}_{522}$.

Any $(\gamma,\delta,\varkappa)\in Q_{252}$ satisf\/ies the left hand side of~\eqref{A_1, -gamma,0 le delta, -1 < varkappa} or~\eqref{A_1, 0 le delta < -gamma, -1 < varkappa}, and satisf\/ies the left hand side of~\eqref{A_5, gamma,varkappa le -1/2, 0 le delta} and~\eqref{A_7, gamma,varkappa < -1/3 < delta}. So, when $m\ge n$, the estimate~\eqref{preliminary estimate} is guaranteed for any $(\alpha,\beta,\gamma,\delta,\varkappa)\in{\mathbb{R}}^2\times Q_{252}$ satisfying both~\eqref{A_1, -gamma,0 le delta, -1 < varkappa} and~\eqref{A_1, 0 le delta < -gamma, -1 < varkappa}, or any of~\eqref{A_5, gamma,varkappa le -1/2, 0 le delta} or~\eqref{A_7, gamma,varkappa < -1/3 < delta}. On ${\mathbb{R}}^2\times Q_{252}$, these conditions mean that $(\alpha,\beta,\gamma,\delta,\varkappa)$ satisf\/ies~\eqref{0 le gamma + delta => alpha + gamma, alpha + beta + gamma + delta + varkappa + 1 < 0} and~\eqref{gamma + delta < 0 => 4 cases}, def\/ining ${\mathfrak{S}}_{252}$.

Any $(\gamma,\delta,\varkappa)\in Q_{155}$ satisf\/ies the left hand side of~\eqref{A_1, -gamma,0 le delta, -1 < varkappa} or~\eqref{A_1, 0 le delta < -gamma, -1 < varkappa}, and satisf\/ies the left hand side of~\eqref{A_4, gamma < -1/2 < delta, 0 le varkappa},~\eqref{A_3, gamma le -1, 0 le delta,varkappa} and~\eqref{A_7, gamma < -1/3 < delta,varkappa}. So, when $m\ge n$, the estimate~\eqref{preliminary estimate} is guaranteed for any $(\alpha,\beta,\gamma,\delta,\varkappa)\in{\mathbb{R}}^2\times Q_{155}$ satisfying both~\eqref{A_1, -gamma,0 le delta, -1 < varkappa} and~\eqref{A_1, 0 le delta < -gamma, -1 < varkappa}, or any of~\eqref{A_4, gamma < -1/2 < delta, 0 le varkappa},~\eqref{A_3, gamma le -1, 0 le delta,varkappa} or~\eqref{A_7, gamma < -1/3 < delta,varkappa}. On ${\mathbb{R}}^2\times Q_{155}$, these conditions mean that $(\alpha,\beta,\gamma,\delta,\varkappa)$ satisf\/ies~\eqref{0 le gamma + delta => alpha + gamma, alpha + beta + gamma + delta + varkappa + 1 < 0} and~\eqref{gamma + delta < 0 => 5 cases}, def\/ining ${\mathfrak{S}}_{155}$.

Any $(\gamma,\delta,\varkappa)\in Q_{212}$ satisf\/ies the left hand side of~\eqref{A_2, gamma le varkappa,0, delta le -1} or~\eqref{A_2, varkappa le gamma,0, delta le -1}. So, when $m\ge n$, the estimate~\eqref{preliminary estimate} is guaranteed for any $(\alpha,\beta,\gamma,\delta,\varkappa)\in{\mathbb{R}}^2\times Q_{212}$ satisfying both~\eqref{A_2, gamma le varkappa,0, delta le -1} and~\eqref{A_2, varkappa le gamma,0, delta le -1}. On ${\mathbb{R}}^2\times Q_{212}$, these conditions become~\eqref{gamma le varkappa => alpha + gamma, alpha + beta + varkappa < 0} and~\eqref{varkappa le gamma => alpha + gamma, alpha + beta + gamma < 0}, def\/ining ${\mathfrak{S}}_{212}$.

\subsubsection{The preliminary estimate is satisf\/ied on a convex set}\label{sss: convex}

Let us show the convexity of the set of elements $x=(\alpha,\beta,\gamma,\delta,\varkappa)\in{\mathbb{R}}^5$ satisfying~\eqref{preliminary estimate} for $m\ge n$, with $\omega=\omega(x)>0$. For $i=0,1$, suppose that $x_i=(\alpha_i,\beta_i,\gamma_i,\delta_i,\varkappa_i)$ satisf\/ies~\eqref{preliminary estimate} for $m\ge n$ with $\omega_i=\omega(x_i)>0$. Recall that the case where $m\le n$ follows from the case where $m\ge n$ by using the mapping~\eqref{R^5 to R^5}.

For $0<t<1$, let
	\begin{gather*}
		x_t=(\alpha_t,\beta_t,\gamma_t,\delta_t,\varkappa_t)=(1-t)x_0+tx_1,\qquad\omega_t=(1-t)\omega_0+t\omega_1>0.
	\end{gather*}
By H\"older inequality, for all $m\ge n$,
	\begin{gather*}
		\sum_{p=0}^n(m-p+1)^{\gamma_t}(n-p+1)^{\delta_t}(p+1)^{\varkappa_t}\\
		\qquad{}
		=\sum_{p=0}^n\!\big((m-p+1)^{\gamma_0}(n-p+1)^{\delta_0}(p+1)^{\varkappa_0}\big)^{1-t}
\big((m-p+1)^{\gamma_1}(n-p+1)^{\delta_1}(p+1)^{\varkappa_1}\big)^t\!\\
		\qquad{}
		\le\left(\sum_{p=0}^n(m-p+1)^{\gamma_0}(n-p+1)^{\delta_0}(p+1)^{\varkappa_0}\right)^{1-t}\\
		\qquad{}
		\phantom{\le\text{}}\text{}\times\left(\sum_{p=0}^n(m-p+1)^{\gamma_1}(n-p+1)^{\delta_1}(p+1)^{\varkappa_1}\right)^t.
	\end{gather*}
So
	\begin{gather*}
		(m+1)^{\alpha_t}(n+1)^{\beta_t}\sum_{p=0}^n(m-p+1)^{\gamma_t}(n-p+1)^{\delta_t}(p+1)^{\varkappa_t}\\
		\qquad{}
		\preccurlyeq\big((m+1)^{-\omega_0}(n+1)^{-\omega_0}\big)^{1-t}\big((m+1)^{-\omega_1}(n+1)^{-\omega_1}\big)^t
		=(m+1)^{-\omega_t}(n+1)^{-\omega_t}.
	\end{gather*}
Thus $x_t$ satisf\/ies~\eqref{preliminary estimate} for $m\ge n$ with $\omega_t$. This completes the proof of Lemma~\ref{l: preliminary estimate}.

\section{The main estimates}\label{s: main estimates}

Here, we show the estimates used in the proofs of Propositions~\ref{p: |c'_k,ell|, case where sigma ne theta = tau} and~\ref{p: |c'_k,ell|, case where sigma ne theta ne tau}. We continue with the notation of Section~\ref{s: preliminary estimate}. Moreover let $(\sigma,\tau,\theta)$ denote the standard coordinates of ${\mathbb{R}}^3$. Consider the af\/f\/ine injection ${\mathbb{R}}^3\to{\mathbb{R}}^5$ and the af\/f\/ine isomorphism of ${\mathbb{R}}^3$ def\/ined by
	\begin{gather}
(\sigma,\tau,\theta)\mapsto \big(\tfrac{1}{4}-\tfrac{\sigma}{2},-\tfrac{1}{4}-\tfrac{\tau}{2},\sigma-\theta-1,\tau-\theta,\theta-\tfrac{1}{2}\big),\label{R^3 to R^5}\\
		(\sigma,\tau,\theta)\mapsto(\tau+1,\sigma-1,\theta).\label{R^3 to R^3}
	\end{gather}
The mapping~\eqref{R^3 to R^3} is the ref\/lection with respect to the plane def\/ined by $\sigma=\tau+1$, and it corresponds to the mapping~\eqref{R^5 to R^5} via~\eqref{R^3 to R^5}. Let $\check{\mathfrak{K}},\check{\mathfrak{K}}'\subset{\mathbb{R}}^3$ be the inverse images of~$\check{\mathfrak{S}}$,~$\check{\mathfrak{S}}'$ by~\eqref{R^3 to R^5}, and let $\check{\mathfrak{K}}_{\text{\rm conv}}$, $\check{\mathfrak{K}}'_{\text{\rm conv}}$ be their convex hulls. So $\check{\mathfrak{K}}_{\text{\rm conv}}$, $\check{\mathfrak{K}}'_{\text{\rm conv}}$ are contained in the inverse images of $\check{\mathfrak{S}}_{\text{\rm conv}}$, $\check{\mathfrak{S}}'_{\text{\rm conv}}$ by~\eqref{R^3 to R^5}, and $\check{\mathfrak{K}}'$, $\check{\mathfrak{K}}'_{\text{\rm conv}}$ are the images of $\check{\mathfrak{K}}$, $\check{\mathfrak{K}}_{\text{\rm conv}}$ by~\eqref{R^3 to R^3}. Thus $\check{\mathfrak{K}}_{\text{\rm conv}}\cap\check{\mathfrak{K}}'_{\text{\rm conv}}$ is symmetric with respect to the plane $\sigma=\tau+1$. We will show the following.

\begin{Lemma}\label{l: check fK_conv cap check fK'_conv}
$\check{\mathfrak{K}}_{\text{\rm conv}}\cap\check{\mathfrak{K}}'_{\text{\rm conv}}$ consists of the elements $(\sigma,\tau,\theta)\in{\mathbb{R}}^3$ that satisfy~\eqref{VV, d}.
\end{Lemma}

The following is a direct consequence of Lemmas~\ref{l: preliminary estimate} and~\ref{l: check fK_conv cap check fK'_conv}.

\begin{Corollary}\label{c: fK_conv cap fK'_conv}
	If $(\sigma,\tau,\theta)\in{\mathbb{R}}^3$ satisfies~\eqref{VV, d}, then there is some $\omega>0$ such that~\eqref{preliminary estimate} holds with the image $(\alpha,\beta,\gamma,\delta,\varkappa)$ of $(\sigma,\tau,\theta)$ by~\eqref{R^3 to R^5}.
\end{Corollary}

Let $\check{\mathfrak{J}}\subset{\mathbb{R}}^2$ be the inverse image of $\check{\mathfrak{K}}_{\text{\rm conv}}\cap\check{\mathfrak{K}}'_{\text{\rm conv}}$ by the af\/f\/ine injection ${\mathbb{R}}^2\to{\mathbb{R}}^3$, $(\sigma,\tau)\mapsto(\sigma,\tau,\tau)$. The following is a direct consequence of Lemma~\ref{l: check fK_conv cap check fK'_conv}.

\begin{Lemma}\label{l: check fJ}
	$\check{\mathfrak{J}}$ consists of the elements $(\sigma,\tau)\in{\mathbb{R}}^2$ that satisfy~\eqref{VV, b}.
\end{Lemma}

Lemma~\ref{l: check fJ} and Corollary~\ref{c: fK_conv cap fK'_conv} have the following direct consequence.

\begin{Corollary}\label{c: check fJ}
	If $(\sigma,\tau)\in{\mathbb{R}}^2$ satisfies~\eqref{VV, b}, then there is some $\omega>0$ such that~\eqref{preliminary estimate} holds with the image $(\alpha,\beta,\gamma,\delta,\varkappa)$ of $(\sigma,\tau,\tau)$ by~\eqref{R^3 to R^5}.
\end{Corollary}

\begin{figure}[h]
\centering
\subfigure[$\check{\mathfrak{K}}$]{\includegraphics[width=5cm]{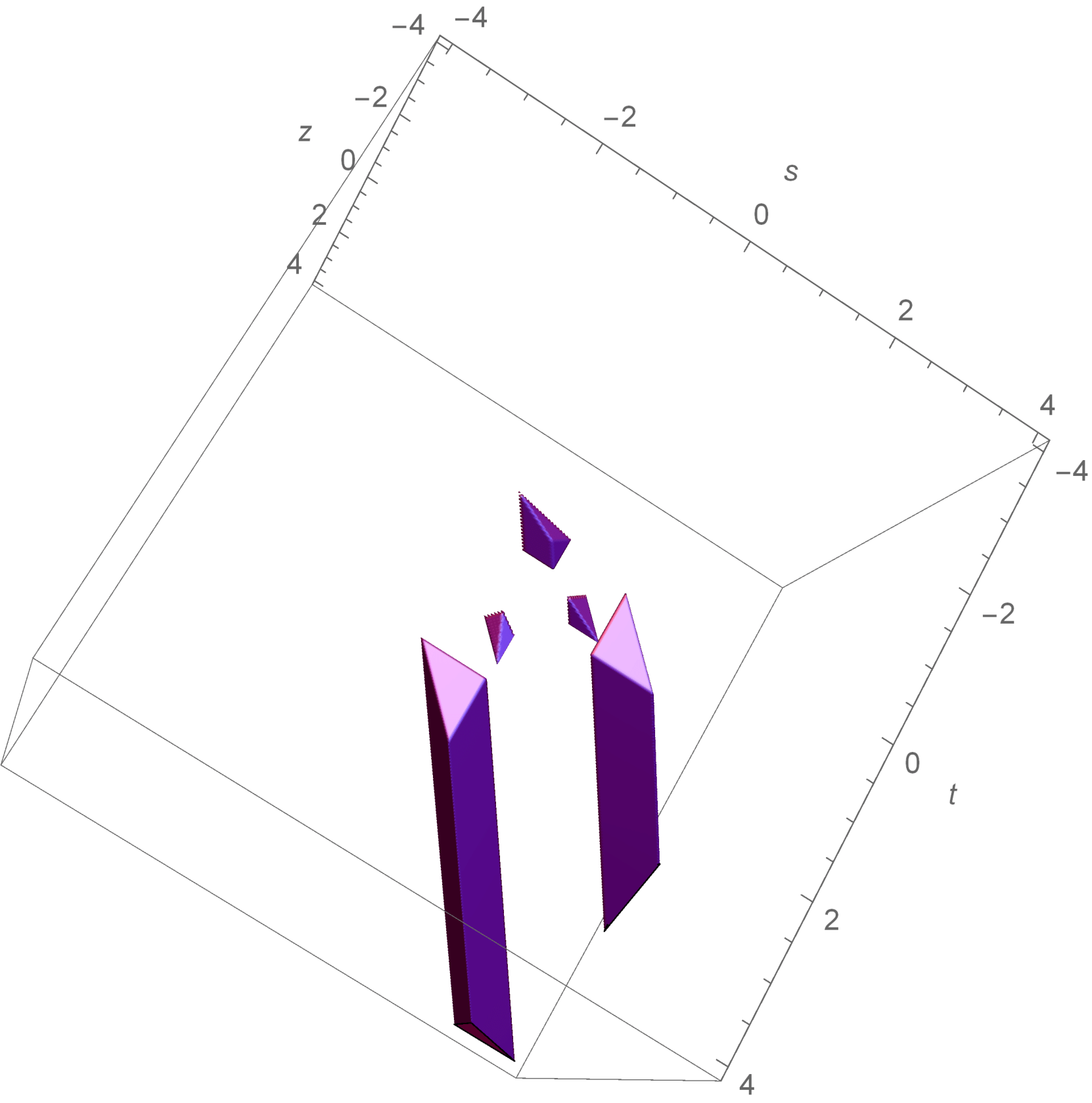}}
\subfigure[$\check{\mathfrak{K}}_{\text{\rm conv}}$]{\includegraphics[width=5cm]{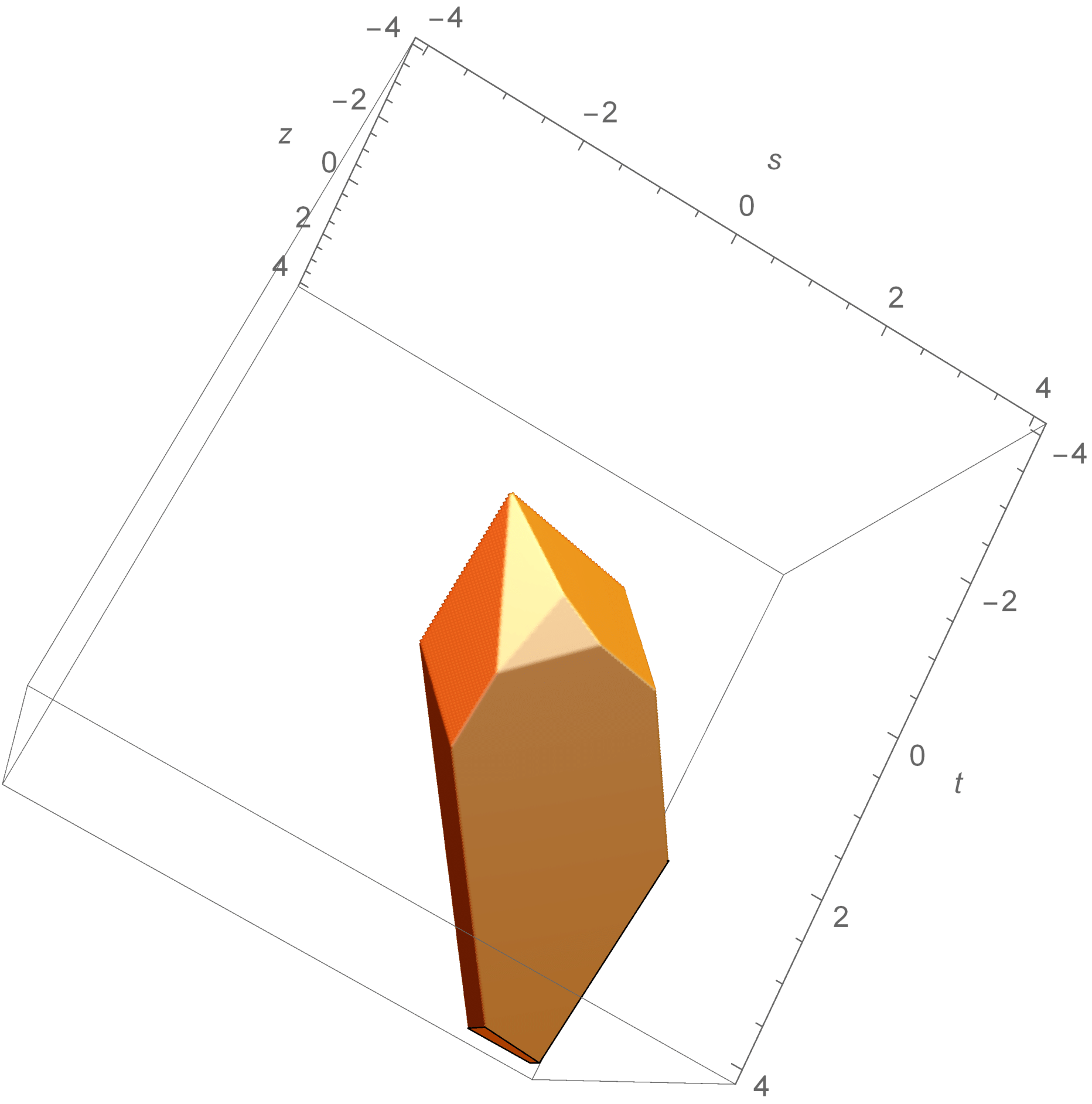}}
\caption{The sets $\check{\mathfrak{K}}$ and $\check{\mathfrak{K}}_{\text{\rm conv}}$.}
\label{fig: check fK & check fK_conv}
\end{figure}

Let us prove Lemma~\ref{l: check fK_conv cap check fK'_conv}. For the subindices $ijk$ equal to 515, 522, 252, 155, 212, let ${\mathfrak{K}}_{ijk}$ and $R_{ijk}$ be the inverse images of~${\mathfrak{S}}_{ijk}$ and~${\mathbb{R}}^2\times Q_{ijk}$ by the mapping~\eqref{R^3 to R^5}. Thus $\check{\mathfrak{K}}={\mathfrak{K}}_{515}\cup{\mathfrak{K}}_{522}\cup{\mathfrak{K}}_{252}\cup{\mathfrak{K}}_{155}\cup{\mathfrak{K}}_{212}$. Moreover, for every $\theta\in{\mathbb{R}}$, let
	\begin{alignat*}{3}
		I_1^1(\theta)&=(-\infty,\theta],&\qquad I_1^2(\theta)&=(-\infty,\theta-1],&\qquad
		I_1^3&=\big({-}\infty,-\tfrac{1}{2}\big],\\
		I_2^1(\theta)&=\big(\theta, \theta+\tfrac{1}{2}\big],&\qquad I_2^2(\theta)&=\big(\theta-1, \theta-\tfrac{1}{2}\big],&\qquad
		I_2^3&=\big({-}\tfrac{1}{2},0\big],\\
		I_3^1(\theta)&=\big( \theta+\tfrac{1}{2},\theta+\tfrac{2}{3}\big],&\qquad I_3^2(\theta)&=\big( \theta-\tfrac{1}{2},\theta-\tfrac{1}{3}\big],&\qquad
		I_3^3&=\big(0,\tfrac{1}{6}\big],\\
		I_4^1(\theta)&=\big(\theta+\tfrac{2}{3},\theta+1\big),&\qquad I_4^2(\theta)&=\big( \theta-\tfrac{1}{3},\theta\big),&\qquad
		I_4^3&=\big( \tfrac{1}{6},\tfrac{1}{2}\big),\\
		I_5^1(\theta)&=[\theta+1,\infty),&\qquad I_5^2(\theta)&=[\theta,\infty),&\qquad I_5^3&=\big[\tfrac{1}{2},\infty\big).
	\end{alignat*}
It can be directly checked that
	\begin{gather*}
		R_{ijk}=\big\{(\sigma,\tau,\theta)\in{\mathbb{R}}^3\,|\,(\sigma,\tau)\in I_i^1(\theta)\times I_j^2(\theta),\ \theta\in I_k^3\big\}.
	\end{gather*}
Simple computations show that, via~\eqref{R^3 to R^5}, the conditions def\/ining the sets ${\mathfrak{S}}_{ijk}$ (Section~\ref{ss: statement of the preliminary estimate}) become the following descriptions of the sets ${\mathfrak{K}}_{ijk}$ (Fig.~\ref{fig: check fK & check fK_conv}(a)):
	\begin{description}
			\item[${\mathfrak{K}}_{515}$:] This is the subset of $R_{515}$ def\/ined by
			\begin{gather*}
				\tfrac{\sigma}{2}-\tfrac{3}{4}<\theta,\qquad\sigma-\tau-3<0.
			\end{gather*}
			
		\item[${\mathfrak{K}}_{522}$:] This is the subset of $R_{522}$ def\/ined by
			\begin{gather*}
				\tfrac{\sigma}{2}-\tfrac{3}{4},\tfrac{\sigma-\tau}{2}-1<\theta.
			\end{gather*}
			
		\item[${\mathfrak{K}}_{252}$:] This is the subset of $R_{252}$ def\/ined by
			\begin{gather}
				\tfrac{\sigma+\tau-1}{2}<\theta,\label{(sigma + tau - 1)/2 < theta}\\
									\begin{cases}
						\frac{\sigma}{2}-\frac{1}{4},\frac{\tau-\sigma}{2}<\theta,\ \text{or}\\[2pt]
						\frac{\sigma}{2}-\frac{5}{12},\frac{\tau-\sigma}{2}+\frac{1}{3}<\theta,\ \text{or}\\[2pt]
						\frac{\sigma}{2}-\frac{3}{4}<\theta,\ \tau-\sigma+1<0.
					\end{cases}\notag
			\end{gather}
			
		\item[${\mathfrak{K}}_{155}$:] This is the subset of $R_{155}$ def\/ined by~\eqref{(sigma + tau - 1)/2 < theta} and
			\begin{gather*}
									\begin{cases}
						\frac{\sigma}{2}+\frac{1}{4}<\theta,\ \tau-\sigma-1<0,\ \text{or}\\[2pt]
						\frac{\sigma}{2}-\frac{1}{4}<\theta,\ \tau-\sigma<0,\ \text{or}\\[2pt]
						\frac{\sigma}{2}-\frac{5}{12}<\theta,\ \tau-\sigma+\frac{1}{3}<0,\ \text{or}\\[2pt]
						\frac{\sigma}{2}-\frac{3}{4}<\theta,\ \tau-\sigma+1<0.
					\end{cases}
			\end{gather*}
			
		\item[${\mathfrak{K}}_{212}$:] This is the subset of $R_{212}$ def\/ined by
			\begin{gather*}
				\tfrac{\sigma}{2}-\tfrac{1}{4}\le\theta \Rightarrow \theta<\tfrac{\sigma+\tau+1}{2},\\
				\theta\le\tfrac{\sigma}{2}-\tfrac{1}{4}\Rightarrow \tfrac{\sigma}{2}-\tfrac{3}{4},\tfrac{\sigma-\tau}{2}-1<\theta.
			\end{gather*}
	
	\end{description}
With tedious computations assisted by graphics produced with Mathematica, it follows that $\check{\mathfrak{K}}_{\text{\rm conv}}$ is the open subset of ${\mathbb{R}}^3$ def\/ined by (Fig.~\ref{fig: check fK & check fK_conv}(b))
	\begin{gather}\label{check fK_conv}
				\tfrac{\sigma-\tau}{2}-1,\tfrac{\tau-\sigma}{2},\tfrac{\sigma+\tau-1}{4},\tfrac{\sigma+3\tau-2}{14},
				\tfrac{\sigma+\tau-1}{2}<\theta<\tfrac{\sigma+\tau+1}{2},\qquad
				\tau-1<\sigma<\tau+3.
	\end{gather}
This is a ``semi-inf\/inite bar'' with 4 lateral faces, and 4 faces at the ``bounded end''.

Applying the af\/f\/ine transformation~\eqref{R^3 to R^3} to this description, we get that $\check{\mathfrak{K}}'_{\text{\rm conv}}$ consists of the triples $(\sigma,\tau,\theta)\in{\mathbb{R}}^3$ satisfying the following conditions:
	\begin{gather}\label{check fK'_conv}
				\tfrac{\sigma-\tau}{2}-1,\tfrac{\tau-\sigma}{2},\tfrac{\sigma+\tau-1}{4},\tfrac{3\sigma+\tau-4}{14},
				\tfrac{\sigma+\tau-1}{2}<\theta<\tfrac{\sigma+\tau+1}{2},\qquad
				\tau-1<\sigma<\tau+3.
	\end{gather}
Combining~\eqref{check fK_conv} and~\eqref{check fK'_conv}, it follows that $\check{\mathfrak{K}}_{\text{\rm conv}}\cap\check{\mathfrak{K}}'_{\text{\rm conv}}$ is given by~\eqref{VV, d} (Fig.~\ref{fig: VV, d}), completing the proof of Lemma~\ref{l: check fK_conv cap check fK'_conv}.

\begin{figure}[h]
\centering
\subfigure[$\bigcup_{i,j}{\mathfrak{K}}_{ij2}$]{\includegraphics[width=5cm]{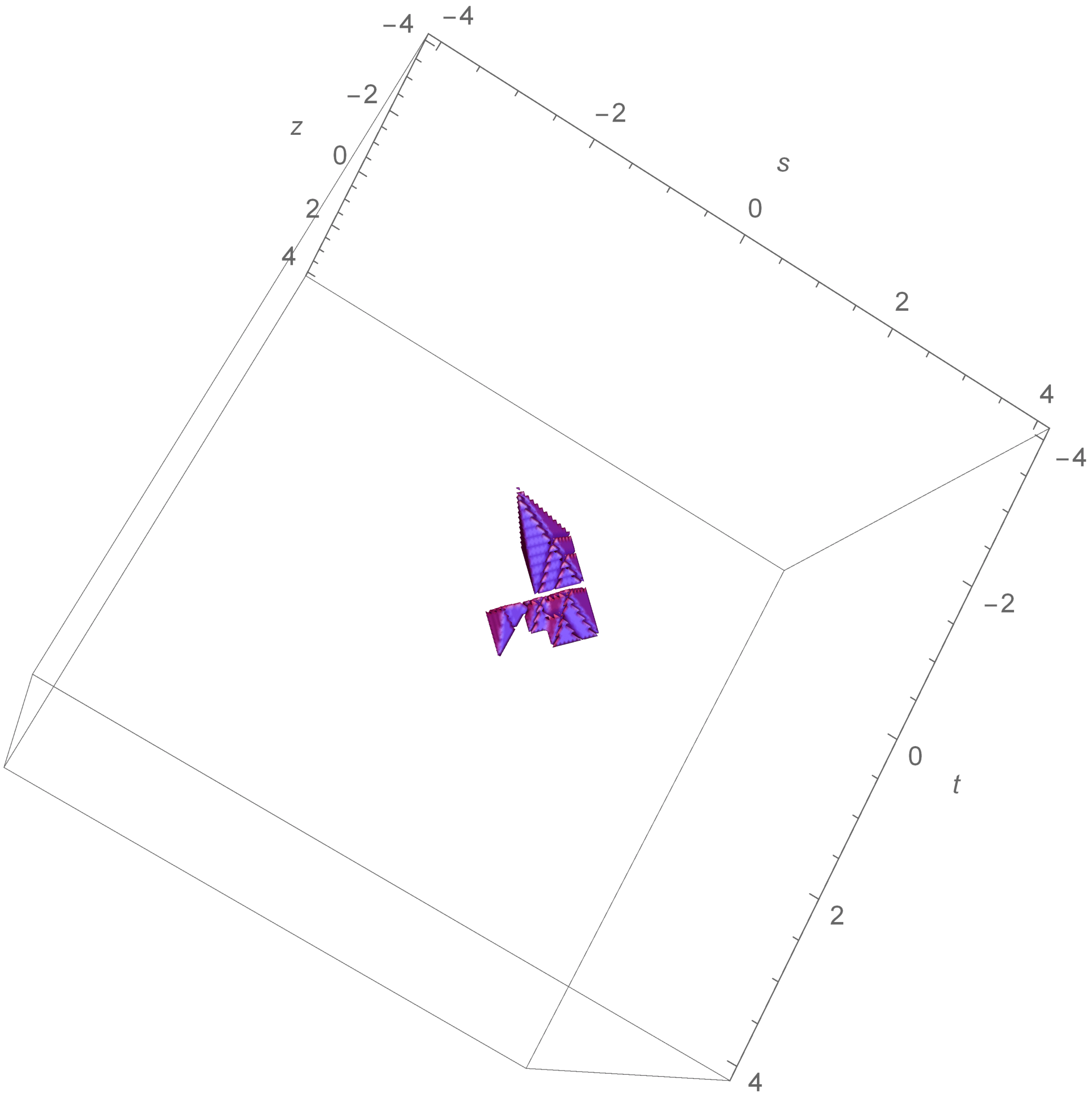}}
\subfigure[$\bigcup_{i,j}{\mathfrak{K}}_{ij3}$]{\includegraphics[width=5cm]{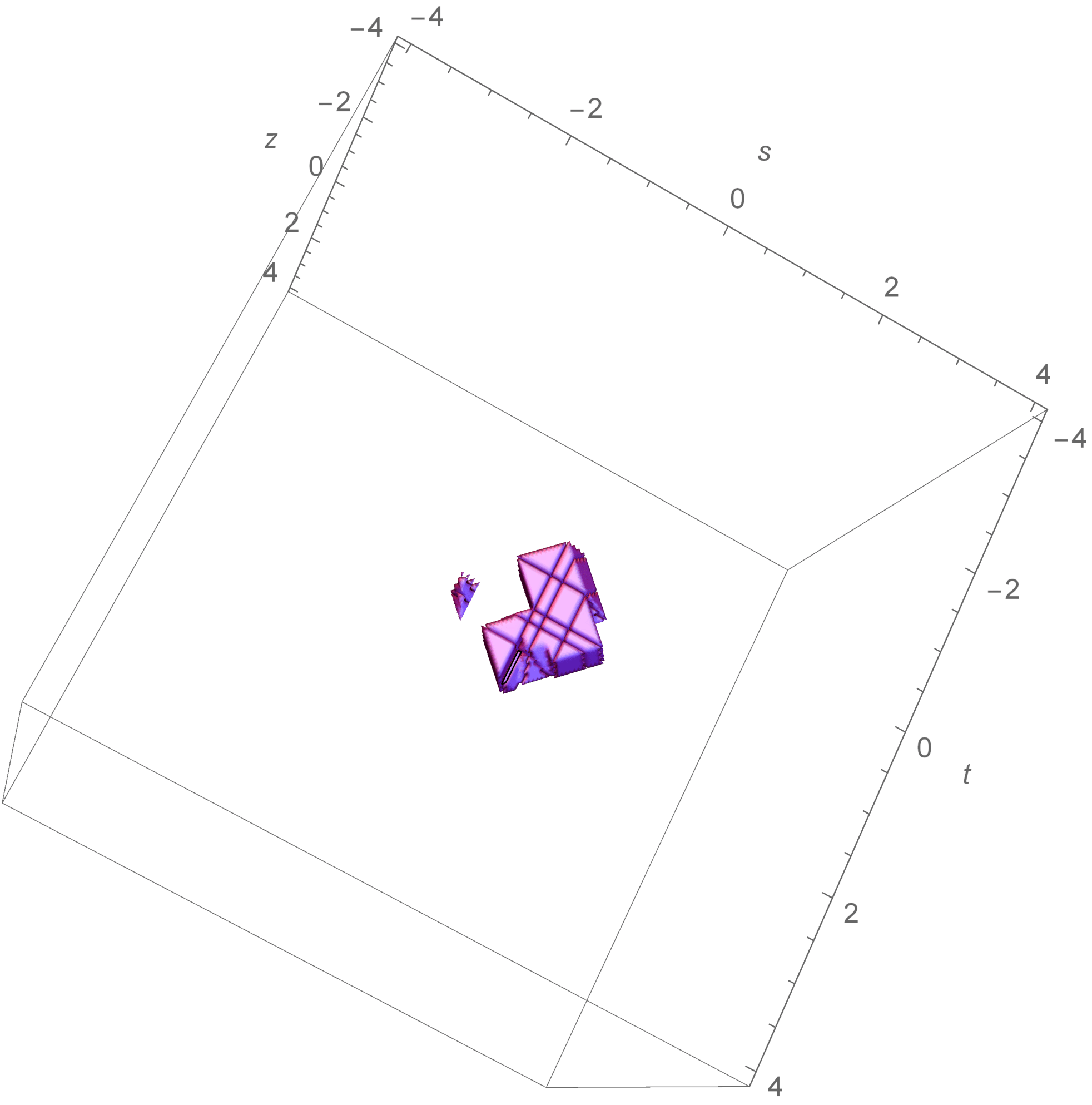}}
\subfigure[$\bigcup_{i,j}{\mathfrak{K}}_{ij4}$]{\includegraphics[width=5cm]{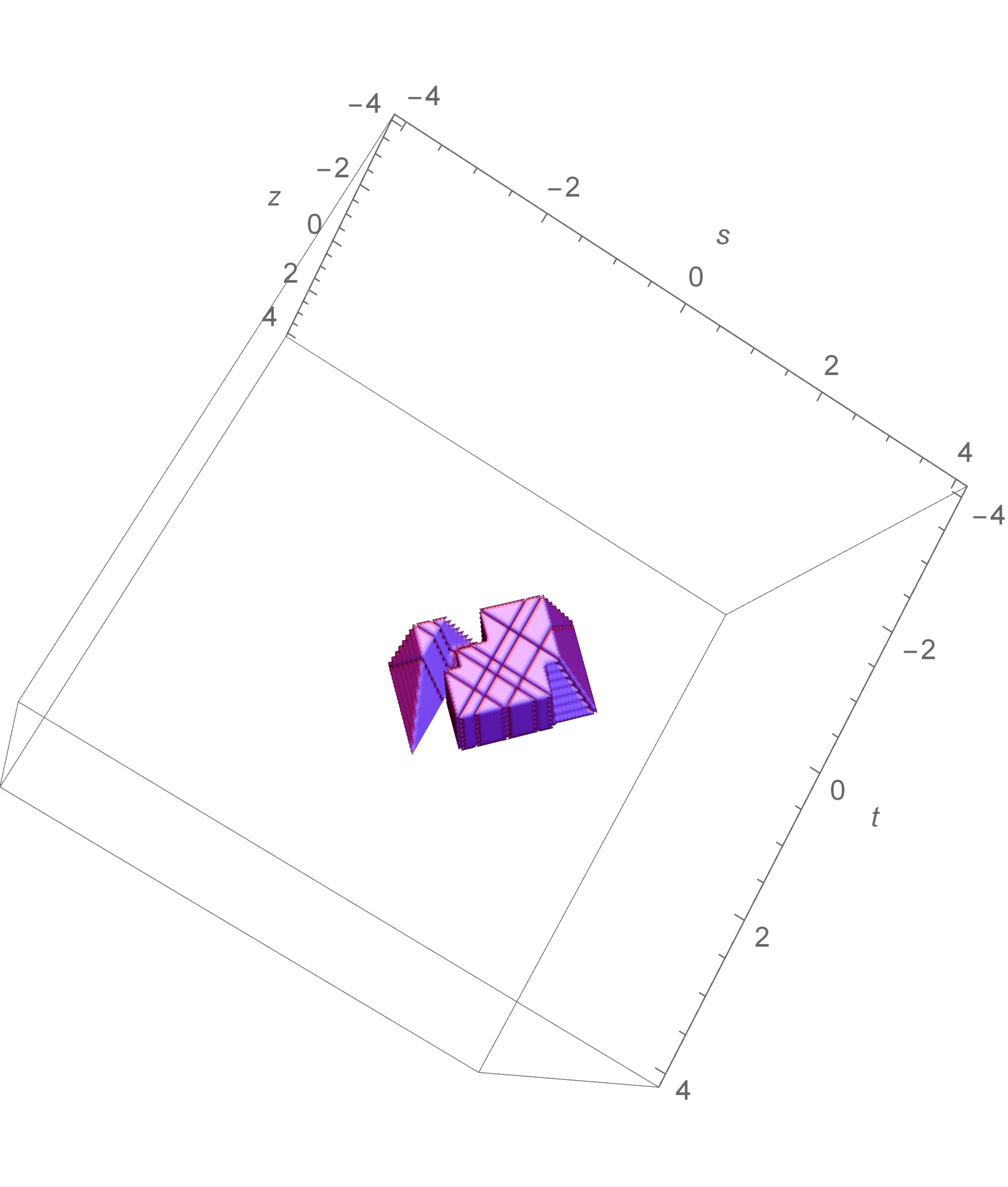}}
\subfigure[$\bigcup_{i,j}{\mathfrak{K}}_{ij5}$]{\includegraphics[width=5cm]{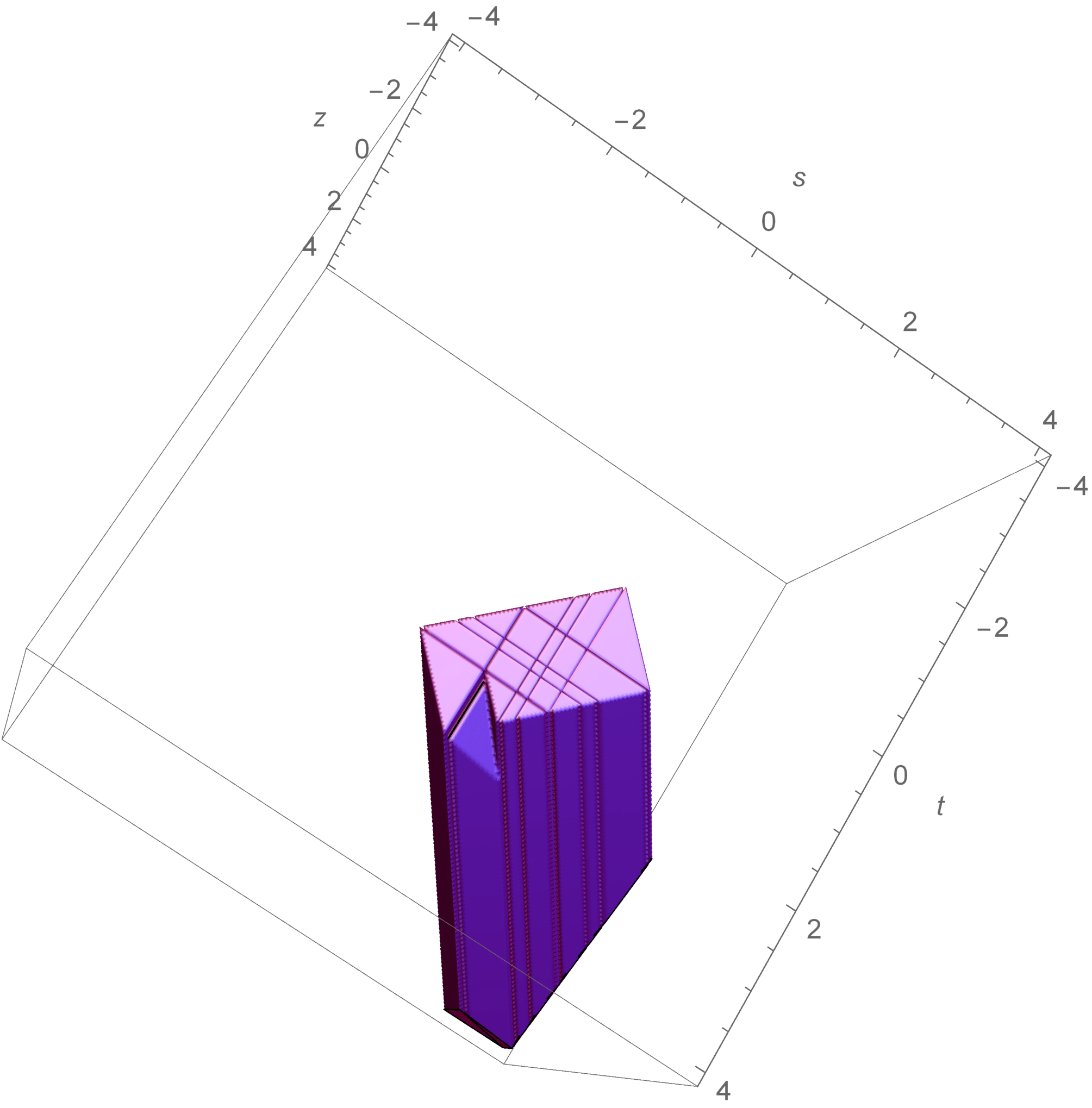}}
\subfigure[${\mathfrak{K}}$]{\includegraphics[width=5cm]{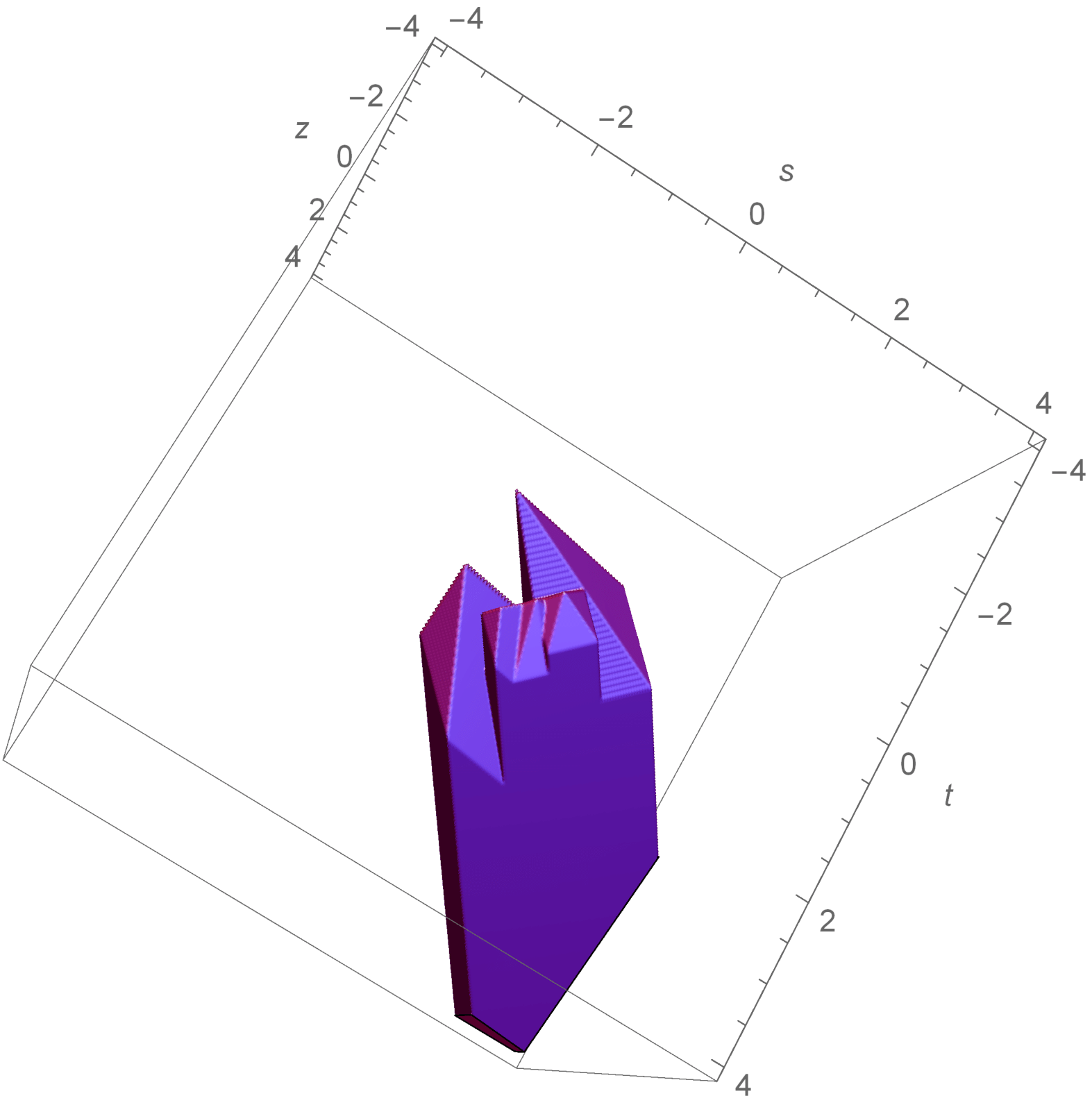}}
\caption{Construction of ${\mathfrak{K}}$.}\label{fig: construction of fK}
\end{figure}

\begin{Remark}\label{r: hidden work}
	In Sections~\ref{sss: A_1}--\ref{sss: A_7}, we have only written the cases that provide the most general conditions to def\/ine ${\mathfrak{S}}_{515},{\mathfrak{S}}_{522},{\mathfrak{S}}_{252},{\mathfrak{S}}_{155},{\mathfrak{S}}_{212}$. But indeed much more hidden work was needed to produce this shorter proof:
		\begin{itemize}\itemsep=0pt
			
			\item We have computed all cases in Sections~\ref{sss: A_1}--\ref{sss: A_7}, giving rise to seven long lists of conditions that guarantee~\eqref{preliminary estimate} when $m\ge n$.
			
			\item We have studied which of those conditions are the most general ones on every ${\mathbb{R}}^2\times Q_{ijk}$, for all $ijk=1,\dots,5$. This produces $125$ sets ${\mathfrak{S}}_{ijk}$, whose inverse images by~\eqref{R^3 to R^5} give 125 sets ${\mathfrak{K}}_{ijk}$. The corresponding unions are denoted by ${\mathfrak{S}}$ and ${\mathfrak{K}}$, and their convex hulls by~${\mathfrak{S}}_{\text{\rm conv}}$ and~${\mathfrak{K}}_{\text{\rm conv}}$.
			
			 \item We got that 41 sets ${\mathfrak{K}}_{ijk}$ are empty, including the 25 sets of the form ${\mathfrak{K}}_{ij1}$, and the remaining 84 sets ${\mathfrak{K}}_{ijk}$ f\/it together forming a ``semi-inf\/inite bar'' (Fig.~\ref{fig: construction of fK}).
			
			\item With tedious computations, we have shown that ${\mathfrak{K}}_{\text{\rm conv}}$ is given by~\eqref{check fK_conv}.
			
			\item We have chosen the most simple family, ${\mathfrak{K}}_{515}$, ${\mathfrak{K}}_{522}$, ${\mathfrak{K}}_{252}$, ${\mathfrak{K}}_{155}$, ${\mathfrak{K}}_{212}$, def\/ining the same convex hull ($\check{\mathfrak{K}}_{\text{\rm conv}}={\mathfrak{K}}_{\text{\rm conv}}$).
			
			\item Finally, we have made some attempts to improve the estimates of Section~\ref{sss: A_7} by using more general versions of the H\"older inequality \cite{Cheung2001}. Some better estimates were obtained in this way, but they produce the same set ${\mathfrak{K}}_{\text{\rm conv}}$ after taking the convex hull.
		
	\end{itemize}
\end{Remark}

\begin{Remark}\label{r: a better version}
	The set ${\mathfrak{S}}_{\text{\rm conv}}$ may have a simple expression, like ${\mathfrak{K}}_{\text{\rm conv}}$, but its computation became too involved. This is the reason we have used ${\mathfrak{K}}_{\text{\rm conv}}$, obtaining the conditions of Theorem~\ref{t: VV}, which are general enough for our applications in~\cite{AlvCalazaFranco:Witten-general}. But, of course, the inverse image of ${\mathfrak{S}}_{\text{\rm conv}}$ by~\eqref{R^3 to R^5} is possibly larger than ${\mathfrak{K}}_{\text{\rm conv}}$. Therefore a simple expression of ${\mathfrak{S}}_{\text{\rm conv}}$ would possibly give a better version of Theorem~\ref{t: VV}. Even a simple expression of $\check{\mathfrak{S}}_{\text{\rm conv}}$ would possibly give a better version of Theorem~\ref{t: VV}.
\end{Remark}

\section[Operators induced on ${\mathbb{R}}_+$]{Operators induced on $\boldsymbol{{\mathbb{R}}_+}$}\label{s: R_+}

Let ${\mathcal{S}}_{\text{\rm ev/odd},+}=\{ \phi|_{{\mathbb{R}}_+}\,|\,\phi\in{\mathcal{S}}_{\text{\rm ev/odd}} \}$. For $c,d>-\frac{1}{2}$, let $L^2_{c,+}=L^2({\mathbb{R}}_+,x^{2c} dx)$ and $L^2_{c,d,+}=L^2_{c,+}\oplus L^2_{d,+}$, whose scalar products are denoted by $\langle\ ,\ \rangle_c$ and $\langle\ ,\ \rangle_{c,d}$, respectively. For $c_1,c_2,d_1,d_2\in{\mathbb{R}}$, let{\samepage	
	\begin{gather*}
		P_0=H-2c_1x^{-1} \frac{d}{dx}+c_2x^{-2} ,\qquad
		Q_0=H-2d_1\frac{d}{dx} x^{-1}+d_2x^{-2} .
	\end{gather*}
Morever let $\xi>0$ and $\eta,\theta\in{\mathbb{R}}$.}

\begin{Corollary}\label{c: PP}
	If $a^2+(2c_1-1)a-c_2=0$, $0<u<1$ and $\sigma:=a+c_1>u-\frac{1}{2}$, then there is a~positive self-adjoint operator ${\mathcal{P}}$ in $L^2_{c_1,+}$ satisfying the following:
 		\begin{itemize}\itemsep=0pt

 			\item[$(i)$] $x^a{\mathcal{S}}_{\text{\rm ev},+}$ is a core of ${{\mathcal{P}}}^{1/2}$ and, for all $\phi,\psi\in x^a{\mathcal{S}}_{\text{\rm ev},+}$,
				\begin{gather*}
					\big\langle{{\mathcal{P}}}^{1/2}\phi,{{\mathcal{P}}}^{1/2}\psi\big\rangle_{c_1}=\langle P_0\phi,\psi\rangle_{c_1}
					+\xi \big\langle x^{-u}\phi,x^{-u}\psi\big\rangle_{c_1} .
				\end{gather*}

 			\item[$(ii)$] ${\mathcal{P}}$ has a discrete spectrum. Let $\lambda_0\le\lambda_2\le\cdots$ be its eigenvalues, repeated according to their multiplicity. There is some $D=D(\sigma,u)>0$, and, for each $\epsilon>0$, there is some $C=C(\epsilon,\sigma,u)>0$ so that~\eqref{lambda_k, case of UU} holds for all $k\in2{\mathbb{N}}$.
 		\end{itemize}
\end{Corollary}

\begin{Corollary}\label{c: QQ}
	If $b^2+(2d_1+1)b-d_2=0$, $0<u<1$ and $\tau:=b+d_1>u-\frac{3}{2}$, then there is a~positive self-adjoint operator ${\mathcal{Q}}$ in $L^2_{d_1,+}$ satisfying the following:
 		\begin{itemize}\itemsep=0pt

 			\item[$(i)$] $x^b{\mathcal{S}}_{\text{\rm odd},+}$ is a core of ${{\mathcal{Q}}}^{1/2}$ and, for all $\phi,\psi\in x^b{\mathcal{S}}_{\text{\rm odd},+}$,
				\begin{gather*}
					\big\langle{{\mathcal{Q}}}^{1/2}\phi,{{\mathcal{Q}}}^{1/2}\psi\big\rangle_{d_1}=\langle Q_0\phi,\psi\rangle_{d_1}
					+\xi\big\langle x^{-u}\phi,x^{-u}\psi\big\rangle_{d_1} .
				\end{gather*}

 			\item[$(ii)$] ${\mathcal{Q}}$ has a discrete spectrum. Let $\lambda_1\le\lambda_3\le\cdots$ be its eigenvalues, repeated according to their multiplicity. There is some $D=D(\tau,u)>0$, and, for each $\epsilon>0$, there is some $C=C(\epsilon,\tau,u)>0$ so that~\eqref{lambda_k, case of UU} holds for all $k\in2{\mathbb{N}}+1$, with $\tau$ instead of $\sigma$.
 		\end{itemize}
\end{Corollary}

\begin{Corollary}\label{c: WW}
	Under the conditions of Corollaries~{\rm \ref{c: PP}} and~{\rm \ref{c: QQ}}, if moreover the conditions of Theorem~{\rm \ref{t: VV}} are satisfied with some $\theta>-\frac{1}{2}$, then there is a positive self-adjoint operator ${\mathcal{W}}$ in $L^2_{c_1,d_1,+}$ satisfying the following:
 		\begin{itemize}\itemsep=0pt

 			\item[$(i)$] $x^a{\mathcal{S}}_{\text{\rm ev},+}\oplus x^b{\mathcal{S}}_{\text{\rm odd},+}$ is a core of ${\mathcal{W}}^{1/2}$, and, for $\phi=(\phi_1,\phi_2)$ and $\psi=(\psi_1,\psi_2)$ in $x^a{\mathcal{S}}_{\text{\rm ev},+}\oplus x^b{\mathcal{S}}_{\text{\rm odd},+}$,
				\begin{gather}
					\big\langle{{\mathcal{W}}}^{1/2}\phi,{{\mathcal{W}}}^{1/2}\psi\big\rangle_{c_1,d_1}
					=\langle(P_0\oplus Q_0)\phi,\psi\rangle_{c_1,d_1}
					+\xi\big\langle x^{-u}\phi,x^{-u}\psi\big\rangle_{c_1,d_1}\nonumber\\
					\hphantom{\langle{{\mathcal{W}}}^{1/2}\phi,{{\mathcal{W}}}^{1/2}\psi\rangle_{c_1,d_1}=}{}
					+\eta\big(\big\langle x^{-a-b-1}\phi_2,\psi_1\big\rangle_\theta
					+\big\langle\phi_1,x^{-a-b-1}\psi_2\big\rangle_\theta\big) .\label{langle WW^1/2 phi, WW^1/2 psi rangle_c_1,d_1}
				\end{gather}

 			\item[$(ii)$] ${\mathcal{W}}$ has a discrete spectrum. Its eigenvalues form two groups, $\lambda_0\le\lambda_2\le\cdots$ and $\lambda_1\le\lambda_3\le\cdots$, repeated according to their multiplicity, such that there is some $D=D(\sigma,\tau,u)>0$ and, for every $\epsilon>0$, there are some $C=C(\epsilon,\sigma,\tau,u)>0$ and $E=E(\epsilon,\sigma,\tau,\theta)>0$ so that~\eqref{lambda_k ge ..., case of VV} and~\eqref{lambda_k le ..., case of VV} hold for all $k\in{\mathbb{N}}$.
			
			\item[$(iii)$] If $\tilde u\in{\mathbb{R}}$ satisf\/ies~\eqref{tilde u}, then there is some $D=D(\sigma,\tau,u)>0$ and, for any $\epsilon>0$, there is some $\widetilde C=\widetilde C(\epsilon,\sigma,\tau,u)>0$ so that~\eqref{lambda_k ge ..., case of VV with tilde u} holds for all $k\in{\mathbb{N}}$.
				
			\item[$(iv)$] If $u=\frac{v+1}{2}$ and $\xi\ge|\eta|$, then there is some $\widetilde D=\widetilde D(\sigma,\tau,u)>0$ so that~\eqref{lambda_k ge ..., case of VV with u = (v+1)/2, xi ge |eta|} holds for all $k\in{\mathbb{N}}$.
							
			\item[$(v)$] If we add the term $\xi'\langle\phi_1,\psi_1\rangle_{c_1}+\xi''\langle\phi_2,\psi_2\rangle_{d_1}$ to the right hand side of~\eqref{langle WW^1/2 phi, WW^1/2 psi rangle_c_1,d_1}, for some $\xi',\xi''\in{\mathbb{R}}$, then the result holds as well with the additional term $\max\{\xi',\xi''\}$ in the right hand side of~\eqref{lambda_k le ..., case of VV}, and the additional term, $\xi'$ for $k\in2{\mathbb{N}}$ and $\xi''$ for $k\in2{\mathbb{N}}+1$, in the right hand sides of~\eqref{lambda_k ge ..., case of VV},~\eqref{lambda_k ge ..., case of VV with tilde u} and~\eqref{lambda_k ge ..., case of VV with u = (v+1)/2, xi ge |eta|}.
			
 		\end{itemize}
\end{Corollary}

These corollaries follow directly from Theorems~\ref{t: UU} and~\ref{t: VV} because the given conditions on~$a$ and~$b$ characterize the cases where $P_0$ and $Q_0$ correspond to $|x|^a U_{\sigma,\text{\rm ev}} |x|^{-a}$ and $|x|^b U_{\tau,\text{\rm odd}} |x|^{-b}$, respectively, via the isomorphisms $|x|^a{\mathcal{S}}_{\text{\rm ev}}\to x^a{\mathcal{S}}_{\text{\rm ev},+}$ and $|x|^b{\mathcal{S}}_{\text{\rm odd}}\to x^b{\mathcal{S}}_{\text{\rm odd},+}$ def\/ined by restriction \cite[Theorem~1.4 and Section~5]{AlvCalaza2014}. In fact, Corollaries~\ref{c: PP} and~\ref{c: QQ} are equivalent because, if $c_1=d_1+1$ and $c_2=d_2$, then $Q_0=xP_0x^{-1}$ and $x\colon L^2_{c_1,+}\to L^2_{d_1,+}$ is a unitary operator.

Remarks~\ref{r: UU}(ii) and~\ref{r: bar fp = bar fl + xi bar ft} have obvious versions for these corollaries. In particular, ${\mathcal{P}}=\overline{P}$, ${\mathcal{Q}}=\overline{Q}$ and ${\mathcal{W}}=\overline{W}$, where $P=P_0+\xi x^{-2u}$, $Q=Q_0+\xi x^{-2u}$ and
	\begin{gather*}
		W =
			\begin{pmatrix}
				P & \eta x^{2(\theta-\sigma)+a-b-1} \\
				\eta x^{2(\theta-\tau)+b-a-1} & Q
			\end{pmatrix}
		 =
			\begin{pmatrix}
				P & \eta x^{2(\theta-c_1)-a-b-1} \\
				\eta x^{2(\theta-d_1)-a-b-1} & Q
			\end{pmatrix} ,
	\end{gather*}
with ${\mathsf{D}}(P)=\bigcap_{m=0}^\infty{\mathsf{D}}({{\mathcal{P}}}^m)$, ${\mathsf{D}}(Q)=\bigcap_{m=0}^\infty{\mathsf{D}}({{\mathcal{Q}}}^m)$ and ${\mathsf{D}}(W)=\bigcap_{m=0}^\infty{\mathsf{D}}({{\mathcal{W}}}^m)$. According to Remark~\ref{r: VV}(ii), we can write~\eqref{langle WW^1/2 phi, WW^1/2 psi rangle_c_1,d_1} as
	\begin{gather*}
		\big\langle{\mathcal{W}}^{1/2}\phi,{\mathcal{W}}^{1/2}\psi\big\rangle_{c_1,d_1}
		=\langle(P_0\oplus Q_0)\phi,\psi\rangle_{c_1,d_1}
		+\xi\big\langle x^{-u}\phi,x^{-u}\psi\big\rangle_{c_1,d_1}\\
		\hphantom{\big\langle{\mathcal{W}}^{1/2}\phi,{\mathcal{W}}^{1/2}\psi\big\rangle_{c_1,d_1}=}{}
		+\eta\big(\big\langle x^{-a-b+1}\phi_2,\psi_1\big\rangle_{\theta'}
		+\big\langle\phi_1,x^{-a-b+1}\psi_2\big\rangle_{\theta'}\big),
	\end{gather*}
and we have
	\begin{gather*}
		W=
			\begin{pmatrix}
				P & \eta x^{2(\theta'-c_1)-a-b+1} \\
				\eta x^{2(\theta'-d_1)-a-b+1} & Q
			\end{pmatrix} .
	\end{gather*}

\section{Application to the Witten's perturbation on strata}\label{s: Witten}

Let $M$ be a Riemannian $n$-manifold. Let $d$, $\delta$ and $\Delta$ denote the de~Rham derivative and coderivative, and the Laplacian, with domain the graded space $\Omega_0(M)$ of compactly supported dif\/ferential forms, and let $L^2\Omega(M)$ be the graded Hilbert space of square integrable dif\/ferential forms. Any closed extension ${\mathbf{d}}$ of $d$ in $L^2\Omega(M)$, def\/ining a complex (${\mathbf{d}}^2=0$), is called an {\em ideal boundary condition} ({\em i.b.c.}) of $d$, which def\/ines a self-adjoint extension $\boldsymbol{\Delta}={\mathbf{d}}^*{\mathbf{d}}+{\mathbf{d}}{\mathbf{d}}^*$ of $\Delta$, called the {\em Laplacian} of ${\mathbf{d}}$. There always exists a minimum/maximum i.b.c., $d_{\min}=\overline d$ and $d_{\max}=\delta^*$, whose Laplacians are denoted by $\Delta_{\min/\max}$. We get corresponding cohomologies $H_{\min/\max}(M)$, and versions of Betti numbers and Euler characteristic, $\beta^i_{\min/\max}$ and $\chi_{\min/\max}$. These are quasi-isometric invariants; in particular, $H_{\max}(M)$ is the usual $L^2$ cohomology. If~$M$ is complete, then there is a unique i.b.c., but these concepts become interesting in the non-complete case. For instance, if $M$ is the interior of a compact manifold with non-empty boundary, then $d_{\min/\max}$ is def\/ined by taking relative/absolute boundary conditions. Given $s>0$ and $f\in C^\infty(M)$, the above ideas can be considered as well for the Witten's perturbations $d_s=e^{-sf}de^{sf}=d+s df\wedge$, with formal adjoint $\delta_s=e^{sf}\delta e^{-sf}=\delta-s df\lrcorner$ and Laplacian $\Delta_s$. In fact, this theory can be considered for any elliptic complex.

On the other hand, let us give a rough idea of the concept of {\em stratified space}. It is a Hausdorf\/f, locally compact and second countable space $A$ with a partition into $C^\infty$ manifolds ({\em strata}) satisfying certain conditions. An order on the family of strata is def\/ined so that $X\le Y$ means that $X\subset\overline Y$. With this order relation, the maximum length of chains of strata is called the {\em depth} of~$A$. Then we continue describing $A$ by induction on $\operatorname{depth} A$, as well as its group $\operatorname{Aut}(A)$ of {\em automorphisms}. If $\operatorname{depth} A=0$, then $A$ is just a $C^\infty$ manifold, whose {\em automorphisms} are its dif\/feomorphisms. Now, assume that $\operatorname{depth} A>0$, and the descriptions are given for lower depth. Then it is required that each stratum $X$ has an open neighborhood~$T$ (a~{\em tube}) that is a~f\/iber bundle whose typical f\/iber is a cone $c(L)=(L\times[0,\infty))/(L\times\{0\})$ and structural group $c(\operatorname{Aut}(L))$, where $L$ is a compact stratif\/ication of lower depth (the {\em link} of~$X$), and $c(\operatorname{Aut}(L))$ consists of the homeomorphisms $c(\phi)$ of $c(L)$ induced by the maps $\phi\times\operatorname{id}$ on $L\times[0,\infty)$ ($\phi\in\operatorname{Aut}(L)$). The point $*=L\times\{0\}\in c(L)$ is called the {\em vertex}. An {\em automorphism} of~$A$ is a homeomorphism that restricts to dif\/feomorphisms between the strata, and whose restrictions to their tubes are f\/iber bundle homomorphisms. This completes the description because the depth is locally f\/inite by the local compactness.

\looseness=1 The local trivializations of the tubes can be considered as ``stratif\/ication charts'', giving a~local description of the form ${\mathbb{R}}^m\times c(L)$. Via these charts, a stratum $M$ of $A$ corresponds, either to ${\mathbb{R}}^m\times\{*\}\equiv{\mathbb{R}}^m$, or to ${\mathbb{R}}^m\times N\times{\mathbb{R}}_+$ for some stratum~$N$ of~$L$. The concept of {\em general adapted metric} on $M$ is def\/ined by induction on the depth. It is any Riemannian metric in the case of depth zero. For positive depth, a Riemannian metric~$g$ on $M$ is called a {\em general adapted metric} if, on each local chart as above, $g$ is quasi-isometric, either to the f\/lat Euclidean metric~$g_0$ if $M$ corresponds to~${\mathbb{R}}^m$, or to $g_0+x^{2u}\tilde g+(dx)^2$ if~$M$ corresponds to ${\mathbb{R}}^m\times N\times{\mathbb{R}}_+$, where~$\tilde g$ is a~general adapted metric on $N$, $x$ is the canonical coordinate of ${\mathbb{R}}_+$, and $u>0$ depends on~$M$ and each stratum $X<M$, whose tube is considered to def\/ine the chart. This assignment $X\mapsto u$ is called the {\em type} of the metric. We omit the term ``general'' when we take $u=1$ for all strata.

Assuming that $A$ is compact, it is proved in \cite{AlvCalazaFranco:Witten-general} that, for any general adapted metric~$g$ on a stratum~$M$ of~$A$ with numbers $u\le1$, the Laplacian $\Delta_{\min/\max}$ has a discrete spectrum, its eigenvalues satisfy a weak version of the Weyl's asymptotic formula, and the method of Witten is extended to get Morse inequalities involving the numbers~$\beta^i_{\min/\max}$ and another numbers $\nu^i_{\min/\max}$ def\/ined by the local data around the ``critical points'' of a version of Morse functions on~$M$; here, the ``critical points'' live in the metric completion of~$M$. This is specially important in the case of a~stratif\/ied pseudo-manifold~$A$ with regular stratum~$M$, where $H_{\max}(M)$ is the intersection homology with perversity depending on the type of the metric \cite{Nagase1983,Nagase1986}. Again, we proceed by induction on the depth to prove these assertions. In the case of depth zero, these properties hold because we are in the case of closed manifolds. Now, assume that the depth is positive, and these properties hold for lower depth. Via a globalization procedure and a version of the K\"unneth formula, the computations boil down to the case of the Witten's perturba\-tion~$d_s$ for a stratum $M=N\times(0,\infty)$ of a cone $c(L)$ with an adapted metric $g=x^{2u}\tilde g+(dx)^2$, where we consider the ``Morse function'' $f=\pm x^2/2$.

Let $\tilde d_{\min/\max}$, $\tilde\delta_{\min/\max}$ and $\widetilde\Delta_{\min/\max}$ denote the operators def\/ined as above for~$N$ with~$\tilde g$. Take dif\/ferential forms $0\ne\gamma\in\ker\widetilde\Delta_{\min/\max}$, of degree~$r$, and $0\ne\alpha,\beta\in{\mathsf{D}}(\widetilde\Delta_{\min/\max})$, of degrees~$r$ and~$r-1$, with $\tilde d_{\min/\max}\beta=\mu\alpha$ and $\tilde\delta_{\min/\max}\alpha=\mu\beta$ for some $\mu>0$. Since $\widetilde\Delta_{\min/\max}$ is assumed to have a discrete spectrum, $L^2\Omega(N)$ has a complete orthonormal system consisting of forms of these types. Correspondingly, there is a ``direct sum splitting'' of $d_s$ into the following two types of subcomplexes:
	\begin{gather*}
		\begin{CD}
			C^\infty_0({\mathbb{R}}_+)\,\gamma @>{d_{s,r}}>> C^\infty_0({\mathbb{R}}_+)\,\gamma\wedge dx ,
		\end{CD}\\
		\begin{CD}
			C^\infty_0({\mathbb{R}}_+)\,\beta @>{d_{s,r-1}}>>
			C^\infty_0({\mathbb{R}}_+)\,\alpha+C^\infty_0({\mathbb{R}}_+)\,\beta\wedge dx
			@>{d_{s,r}}>> C^\infty_0({\mathbb{R}}_+)\,\alpha\wedge dx .
		\end{CD}
	\end{gather*}
Forgetting the dif\/ferential form part, they can be considered as two types of simple elliptic complexes of lengths one and two,
	\begin{gather*}
		\begin{CD}
			C^\infty_0({\mathbb{R}}_+) @>{d_{s,r}}>> C^\infty_0({\mathbb{R}}_+) ,
		\end{CD}\\
		\begin{CD}
			C^\infty_0({\mathbb{R}}_+) @>{d_{s,r-1}}>> C^\infty_0({\mathbb{R}}_+)\oplus C^\infty_0({\mathbb{R}}_+)
			@>{d_{s,r}}>> C^\infty_0({\mathbb{R}}_+) .
		\end{CD}
	\end{gather*}
Let $\kappa=(n-2r-1)\frac{u}{2}$. In the complex of length one, $d_{s,r}$ is a densely def\/ined operator of $L^2_{\kappa,+}$ to $L^2_{\kappa,+}$, we have
	\begin{gather*}
 		d_{s,r}=\frac{d}{dx}\pm sx ,\qquad
 		\delta_{s,r}=-\frac{d}{dx}-2\kappa x^{-1}\pm sx ,
 	\end{gather*}
and the corresponding components of the Laplacian are
	\begin{gather*}
 		\Delta_{s,r} =H-2\kappa x^{-1} \frac{d}{dx}\mp s(1+2\kappa) ,\qquad
 		\Delta_{s,r+1} =H-2\kappa\frac{d}{dx} x^{-1}\mp s(-1+2\kappa) .
 	\end{gather*}
Up to the constant terms, these operators are of the form already considered in~\cite{AlvCalaza2014}, without the term with $x^{-2u}$, and the spectrum of $\Delta_{s,\min/\max,r}$ and $\Delta_{s,\min/\max,r+1}$ is well known.

\begin{table}[t]\centering\caption{Self-adjoint extensions of $\Delta_{s,r-1}$ and $\Delta_{s,r+1}$.}
\label{table: D_s,r-1 and D_s,r+1}
\vspace{1mm}
\renewcommand{\arraystretch}{1.3}
\begin{tabular}{|c|c|l||c|c|l|}
\hline
$a$ & $\sigma$ & condition & $b$ & $\tau$ & condition\\
\hline
$0$ & $\kappa+u$ & $\kappa>-\frac{1}{2}$ & $0$ & $\kappa$ & $\kappa>u-\frac{3}{2}$ \\
\hline
$1-2(\kappa+u)$ & $1-\kappa-u$ & $\kappa<\frac{3}{2}-2u$ & $1-2\kappa$ & $-1-\kappa$ & $\kappa<\frac{1}{2}-u$ \\
\hline
\end{tabular}
\end{table}

\begin{table}[t]\centering
\caption{Self-adjoint extensions of $\Delta_{s,r}$.}
\label{table: D_s,r}
\vspace{1mm}
\renewcommand{\arraystretch}{1.3}
\begin{tabular}{|c|c|c|c|c|c|}
\hline
$a$ & $b$ & $\sigma$ & $\tau$ & $\theta$ & Condition \\
\hline
$0$ & $0$ & $\kappa$ & $\kappa+u$ & $\kappa$ & $\kappa>u-\frac{1}{2}$ \\
\hline
$1-2\kappa$ & $-1-2(\kappa+u)$ & $1-\kappa$ & $-1-\kappa-u$ & $-\kappa-u$ & $\kappa<\frac{1}{2}-2u$ \\
\hline
$0$ & $-1-2(\kappa+u)$ & $\kappa$ & $-1-\kappa-u$ & $-\frac{1}{2}-u$ & Impossible \\
\hline
$1-2\kappa$ & $0$ & $1-\kappa$ & $\kappa+u$ & $\frac{1}{2}$ & $-1-\frac{u}{2}<\kappa<1-\frac{u}{2}$ \\
\hline
\end{tabular}
\end{table}

In the complex of length two, $d_{s,r-1}$ is a densely def\/ined operator of $L^2_{\kappa+u,+}$ to $L^2_{\kappa,+}\oplus L^2_{\kappa+u,+}$, $d_{s,r}$ is a densely def\/ined operator of $L^2_{\kappa,+}\oplus L^2_{\kappa+u,+}$ to $L^2_{\kappa,+}$, we have
	 \begin{gather*}
 d_{s,r-1} =
 \begin{pmatrix}
 \mu\\
 \frac{d}{dx}\pm sx
 \end{pmatrix} ,\qquad
 \delta_{s,r-1} =
 \begin{pmatrix}
 \mu x^{-2u} & -\frac{d}{d x}-2(\kappa+u) x^{-1}\pm s x
 \end{pmatrix} ,\\
 d_{s,r} =
 \begin{pmatrix}
 \frac{d}{d x}\pm s x & -\mu
 \end{pmatrix} ,\qquad
 \delta_{s,r} =
 \begin{pmatrix}
 -\frac{d}{d x}-2\kappa x^{-1}\pm s x\\
 -\mu x^{-2u}
 \end{pmatrix} ,
 \end{gather*}
and the corresponding components of the Laplacian are
	\begin{gather*}
		\Delta_{s,r-1} =H-2(\kappa+u)x^{-1} \frac{d}{dx}+\mu^2x^{-2u}\mp s(1+2(\kappa+u)) ,\\
		\Delta_{s,r+1} = H-2\kappa\frac{d}{dx} x^{-1}+\mu^2x^{-2u}\mp s(-1+2\kappa),\\
 		\Delta_{s,r} =
 			\begin{pmatrix}
 				A & -2\mu ux^{-1} \\
 				-2\mu ux^{-2u-1} & B
 			\end{pmatrix} ,
 	\end{gather*}
where
	\begin{gather*}
		A =H-2\kappa x^{-1} \frac{d}{dx}+\mu^2x^{-2u}\mp s(1+2\kappa) ,\\
 		B =H-2(\kappa+u)\frac{d}{dx} x^{-1}+\mu^2x^{-2u}\mp s(-1+2(\kappa+u)) .
 	\end{gather*}
Up to the constant terms, $\Delta_{s,r-1}$ and $A$ are of the form of $P$, and $\Delta_{s,r+1}$ and $B$ are of the form of $Q$, in Section~\ref{s: R_+}. In the case $u=1$, these operators were studied in \cite{AlvCalaza2014}. Thus assume that $u<1$. Then, according to Corollaries~\ref{c: PP}--\ref{c: WW}, we get self-adjoint extensions of $\Delta_{s,r-1}$, $\Delta_{s,r+1}$ and $\Delta_{s,r}$ as indicated in Tables~\ref{table: D_s,r-1 and D_s,r+1} and~\ref{table: D_s,r}, where the conditions are determined by the hypotheses; indeed most possibilities of the hypothesis are needed. With further analysis~\cite{AlvCalazaFranco:Witten-general}, the maximum and minimum Laplacians can be given by appropriate choices of these operators, depending on the values of~$\kappa$. Moreover the eigenvalue estimates of these corollaries play a key role in this research.

If $A$ is a stratif\/ied pseudo-manifold, our restrictions on $u$ allow to get enough metrics to represent all intersection cohomologies of $A$ with perversity less or equal than the lower middle perversity, according to \cite{Nagase1983,Nagase1986}.

\subsection*{Acknowledgements}

The f\/irst author was partially supported by MICINN, Grants MTM2011-25656 and MTM2014-56950-P, and by Xunta de Galicia, Grant Consolidaci\'on e estructuraci\'on 2015 GPC GI-1574. The third author has received f\/inancial support from the Xunta de Galicia and the European Union (European Social Fund - ESF).

\pdfbookmark[1]{References}{ref}
\LastPageEnding

\end{document}